\documentclass[12pt]{amsart}

\usepackage{amsthm}

\usepackage{graphicx}
\usepackage{epsfig}
\usepackage{epstopdf}
\usepackage[totalwidth=17cm, totalheight=25cm]{geometry}

\usepackage{latexsym}
\usepackage{amsmath}
\usepackage{amssymb}
\usepackage{amsfonts}

\theoremstyle{definition}

\newtheorem{dfn}{Definition}[section]
\newtheorem{defn}[dfn]{Definition}
\newtheorem{definition}[dfn]{Definition}

\newtheorem{example}[dfn]{Example}
\newtheorem{bexample}[dfn]{Basic Example}
\newtheorem{rem}[dfn]{Remark}
\newtheorem{remark}[dfn]{Remark}
\theoremstyle{plain}

\newtheorem{thm}[dfn]{Theorem}
\newtheorem{theorem}[dfn]{Theorem}
\newtheorem{ff}[dfn]{Fundamental Fact}

\newtheorem{lem}[dfn]{Lemma}
\newtheorem{lemma}[dfn]{Lemma}
\newtheorem{klemma}[dfn]{Key Lemma}

\newtheorem{prop}[dfn]{Proposition}

\newtheorem{cor}[dfn]{Corollary}
\newtheorem{corollary}[dfn]{Corollary}

\newtheorem{question}[dfn]{Question}

\numberwithin{equation}{section}

\def\proof{\par\medskip\noindent{\it Proof. }}

\def\R{{\mathbb R}}
\def\C{{\mathbb C}}
\def\Z{{\mathbb Z}}
\def\Q{{\mathbb Q}}
\def\H{{\mathbb H}}
\def\N{{\mathbb N}}

\def\eps{\epsilon}
\def\al{\alpha}
\def\be{\beta}
\def\ga{\gamma}
\def\Ga{\Gamma}
\def\de{\delta}

\def\De{\Delta}

\def\Si{\Sigma}
\def\si{\sigma}

\def\la{\lambda}
\def\La{\Lambda}
\def\om{\omega}
\def\Om{\Omega}

\def\acts{\curvearrowright}
\def\amod{a_{mod}}
\def\CF{C_{F\ddot u}}
\def\D{\partial}
\def\DF{\partial_{F\ddot u}}
\def\diamo{\diamondsuit}
\def\Flag{\operatorname{Flag}}
\def\Flagit{\Flag_{\iota\tau_{mod}}}
\def\Flagt{\operatorname{Flag_{\tau_{mod}}}}
\def\Flags{\operatorname{Flag_{\si_{mod}}}}
\def\Fmod{F_{mod}}
\def\geo{\partial_{\infty}}

\def\embed{\hookrightarrow}
\def\Hom{\operatorname{Hom}}
\def\Homeo{\operatorname{Homeo}}
\def\id{\operatorname{id}}
\def\8{\infty}
\def\inte{\operatorname{int}}
\def\Isom{\operatorname{Isom}}
\def\Las{\La_{\si_{mod}}}
\def\Lat{\La_{\tau_{mod}}}
\def\lra{\longrightarrow}
\def\ol{\overline}

\def\Ra{\Rightarrow}
\def\rank{\mathop{\hbox{rank}}}

\def\simod{\si_{mod}}
\def\st{\operatorname{st}}
\def\stF{\operatorname{st}_{F\ddot u}}
\def\taumod{\tau_{mod}}
\def\Th{\mathop{\hbox{Th}}\nolimits}
\def\tits{\partial_{Tits}}
\def\Wt{W_{\tau_{mod}}}

\def\<{\langle}
\def\>{\rangle}

\begin{document}

\title{	Discrete isometry groups of symmetric spaces}
\author{Michael Kapovich, Bernhard Leeb}
\date{March 3, 2017}
\maketitle

\tableofcontents  

This survey
is based on a series of lectures that we gave at MSRI in Spring 2015 and on a series of papers, 
mostly written jointly with Joan Porti \cite{coco13, coco15, morse, mlem, anolec,bordif}. 
A shorter survey of our work appeared in \cite{anosov}. 
Our goal here is to:

\smallskip
1. Describe a class of discrete subgroups $\Ga<G$ of higher rank semisimple Lie groups, 
to be called RA (regular antipodal) subgroups, 
which exhibit some ``rank 1 behavior''. 

2. Give different characterizations of the subclass of {\em Anosov subgroups},
which generalize convex-cocompact subgroups of rank 1 Lie groups,
in terms of various equivalent dynamical and geometric properties 
(such as asymptotically embedded, RCA, Morse, URU).

3. Discuss the topological dynamics of discrete subgroups $\Ga$
on flag manifolds associated to $G$ 
and Finsler compactifications of associated symmetric spaces $X=G/K$.
Find domains of proper discontinuity 
and use them to construct natural bordifications and compactifications of the locally symmetric spaces $X/\Gamma$. 

\smallskip 
The ultimate goal of this project is to find a class of {\em geometrically finite} discrete subgroups of higher rank Lie groups which  includes some natural classes of discrete subgroups such as 
Anosov subgroups and arithmetic groups. One reason for looking for such a class is that one should be able to prove structural theorems about such groups and associated quotient spaces, analogous to theorems in the case of geometrically finite isometry groups of rank 1 symmetric spaces. 
This is somewhat analogous to the subclass of hyperbolic groups among all finitely presented groups: While there are very few general results about finitely presented groups, there are many interesting results about hyperbolic groups. 
Our work is guided by the theory of Kleinian groups and, more generally, the theory of discrete subgroups of rank 1 Lie groups and quasiconvex subgroups of Gromov-hyperbolic groups. However, there are instances when the ``hyperbolic intuition'' leads one astray and naive generalizations from rank 1 fail in the higher rank setting: One of the earliest examples of this failure is the fact that {\em convex cocompactness} does not have a straightforward generalization 
to higher rank and simply repeating the rank 1 definitions results in a too narrow   class of groups, see \cite{convcoco} and \cite{Quint}. 
Many proofs and constructions in the theory of Anosov subgroups are more complex than the ones for convex cocompact subgroups of rank 1 Lie groups. Partly, this is due to the fact that in higher rank one has to take into account the rich combinatorial structure of Tits boundaries of symmetric spaces (the latter is trivial in rank 1).

\smallskip 
{\bf Organization of the paper.}  In section \ref{ch:1}, we review the theory of discrete isometry groups of negatively curved symmetric spaces with main emphasis on geometrically finite groups and their properties. Section \ref{sec:2} covers basics of nonpositively curved symmetric spaces $X=G/K$, including visual boundaries, horofunction compactification (with respect to the Riemannian as well as Finsler metrics), regularity of sequences, 
higher rank convergence dynamics 
as well as the 
Higher Rank Morse Lemma for 
regular quasigeodesics in $X$. 
Some of this material is covered in more detail in the four appendices to the paper. 
In section \ref{sec:3}, 
we start the discussion of discrete subgroups $\Ga<G$.
We define the class of regular subgroups
in terms of the asymptotic geometry of their orbits in $X$
and give a dynamical characterization as generalized convergence groups
in terms of their action on the full flag manifold $G/B$.
We then impose further restrictions 
(on their dynamics at infinity, the asymptotic geometry of their orbits,
their coarse geometry)
and state an equivalence theorem for the resulting classes of groups.\footnote{These classes of groups are defined in general 
with respect to faces $\taumod$ of the model spherical Weyl chamber $\simod$,
equivalently, with respect to conjugacy classes of parabolic subgroups $P<G$.
In this survey, for simplicity, we limit ourselves to the case $\taumod=\simod$, equivalently,
when the parabolic subgroups are minimal parabolics (minimal parabolic in the case of complex Lie groups).}
We further study the coarse geometry of Anosov subgroups:
We state a local-to-global principle for the Morse property 
which leads to a new proof for structural stability.
A form of the local-to-global principle 
is used to construct basic examples of Anosov subgroups,
namely free discrete subgroups with controlled coarse extrinsic geometry
(Morse-Schottky subgroups). 
Section \ref{sec:4} deals with 
domains of discontinuity of actions of Anosov subgroups $\Gamma< G$ on the full flag manifold $G/B$ 
as well as on the regular Finsler compactification of $X$. 
Our construction of domains of discontinuity is motivated by Mumford's Geometric Invariant Theory.
One of the applications 
is the existence of a certain manifold with corners compactification of the locally symmetric space  $X/\Gamma$ 
obtained from the Finsler compactification of $X$. 
The existence of such a compactification 
is yet another characterization of Anosov subgroups among uniformly regular subgroups. 
In section \ref{sec:5}, we discuss some potential future directions for the study of geometrically finite isometry groups of higher rank symmetric spaces.

\smallskip
{\bf Acknowledgements.} The first author was partly supported by the NSF grants  DMS-12-05312 and DMS-16-04241. 
He is also thankful to KIAS (the Korea Institute for Advanced Study) for its hospitality. 
This paper grew out of lecture notes 
written during the MSRI program ``Dynamics on Moduli Spaces of Geometric Structures'' in Spring 2015; 
we are grateful to MSRI for conducting this program.
We are also thankful to the referee for useful suggestions.

\section{Kleinian groups -- discrete isometry groups of negatively curved symmetric spaces}\label{ch:1}

\subsection{Basics of negatively curved symmetric spaces} \label{sec:rk1}

In this section we review basic properties of negatively curved symmetric spaces; later on, we will compare and contrast their properties with the ones of {\em higher rank} symmetric spaces of noncompact type. We refer the reader to \cite{Mostow} and \cite{Parker} for a detailed discussion 
of negatively curved symmetric spaces. 

Recall that negatively curved symmetric spaces $X$ (also known as {\em rank 1 symmetric spaces of noncompact type}) come in four families: 
$$\H^n, {\mathbf C}\H^n, {\mathbf H} \H^n, {\mathbf O} \H^2,$$
i.e., real hyperbolic spaces, complex hyperbolic spaces, quaternionic hyperbolic spaces and the octonionic hyperbolic plane. We will normalize their Riemannian metrics so that the maximum of the sectional curvature is $-1$. 
The identity component of the isometry group of $X$ will be denoted $G$. 
 The basic fact of the geometry of negatively curved symmetric spaces is that two geodesic segments in $X$ 
are $G$-congruent if and only if they have the same length.  (This property fails in higher rank.) 

The {\em visual boundary} of a symmetric space $X$ consists of the equivalence classes of geodesic rays in $X$, equipped with a suitable topology. Here two rays are equivalent if and only if they are within finite distance from each other. The visual boundary of $X$ is denoted $S=\geo X$.  The elements of $\geo X$ are called  {\em ideal boundary points} of $X$.  A ray representing a point $\xi\in S$ is said to be {\em asymptotic} to $\xi$. Given any point $x\in X$, there is a unique ray emanating from $x$ and asymptotic to $\xi$; this ray is denoted $x\xi$. Thus, after we fix a basepoint $o\in X$, we can identify $S$ with the unit tangent sphere $U_oX$ in $T_oX$: Each ray $o\xi$ corresponds to its (unit) velocity vector at $o$.  This identification endows $S$ with a natural smooth structure 
and a Riemannian metric depending on $o$. 
The sphere $S$ is a homogeneous  $G$-space and the point stabilizers are the minimal parabolic subgroups $B< G$. 
 
 An important feature of negatively curved symmetric spaces (which also fails in higher rank) is:
 
 \begin{lemma}
For any two asymptotic rays $r_i: [0, \infty)\to X, i=1, 2$, there exist $t_1, t_2\in \R_+$ such that
the rays $r_i: [t_i, \infty)\to X$ are {\em strongly asymptotic:} 
 $$
 \lim_{t\to\infty} d(r_1(t_1+ t), r_2(t_2+t))=0.
 $$
 \end{lemma}

\medskip 
Attaching the visual boundary to $X$ provides the {\em visual compactification} of the symmetric space $X$:
$$
\ol X= X\sqcup S,
$$ 
where a sequence of points $x_n\in X$ converges to an ideal point $\xi\in S$ if and only if the sequence of geodesic segments $ox_n$ converges to the geodesic ray $o\xi$ representing $\xi$. 

Alternatively, one can describe this compactification as the {\em horofunction compactification} of $X$. This compactification is defined for a much larger class of metric spaces and we will use it later in section \ref{sec:finsler geometry} in the context of {\em Finsler metrics} on symmetric spaces. See Appendix \ref{sec:horoboundary}. We now return to negatively curved symmetric spaces. 

\medskip 
{\em Visibility property.} Any two distinct ideal boundary points $\xi, \hat\xi\in S$ are connected by a (unique)  geodesic 
line $l: \R\to X$: 
$$
\lim_{t\to\infty} l(t)=\xi, \quad \lim_{t\to-\infty} l(t)=\hat\xi. 
$$
This property again fails in higher rank. 

\medskip 
{\bf Isometries.} The isometries $g$ of $X$ are classified according to their (convex) displacement functions  
$$
d_g(x)= d(x, gx). 
$$

\begin{itemize}
\item Hyperbolic isometries: $\inf_X d_g > 0$. 
In this case, the infimum is attained on a $g$-invariant geodesic $a_g\subset X$, called the {\em axis} of $g$. The ideal endpoints of $a_g$ are fixed by $g$. 

\item Parabolic isometries: $\inf_X d_g=0$ and is not attained. Each parabolic isometry has a unique fixed ideal boundary point $\xi\in \geo X$,
and 
$$
\xi= \lim_{i\to\infty} x_i
$$
for every sequence $x_i\in X$ such that $\lim_{i\to\infty} d_g(x_i)=0$. Horospheres centered  at $\xi$ are invariant under $g$. 

\item Elliptic isometries: $\inf_X d_g=0$ and is attained, i.e.  $g$ fixes a point in $X$. 
\end{itemize}

In particular, there are no isometries $g$ for which $\inf_X d_g > 0$ and the infimum is not attained. 
This again fails in higher rank.

\medskip
There are two fundamental facts about symmetric spaces  of negative curvature which will 
guide our discussion: The {\em Morse Lemma} and the {\em Convergence Property}. 

\begin{theorem}
[Morse Lemma \cite{Mostow, Gromov_hypgps}] 
Quasigeodesics in $X$ are uniformly close to geodesics. More precisely, given constants $L,A$, there exists a number $D=D(L,A)$ such that  each $(L, A)$ quasigeodesic $q$ in $X$ is within Hausdorff distance $D$ from a geodesic. 
\end{theorem}

While the Morse Lemma fails in higher rank symmetric spaces $X$
(because it fails in euclidean plane), 
an important result that we prove in \cite{morse} is that {\em uniformly regular} quasigeodesics are uniformly close to {\em diamonds} in $X$, see Theorem \ref{thm:HRMorse} below.

\subsection{The rank 1 convergence property}
\label{sec:rank1conv}

Given two points $\al, \om\in S$ we define the {\em quasiconstant map}
$$
\al_\om:  \ol X - \{\om\}\to \{\al\}
$$
which is undefined at $\om$. A sequence in a locally compact topological space is {\em divergent} if it has no accumulation points in the space. We will use this in the context of sequences in Lie groups. We say that a divergent 
sequence $g_k\in G$ is {\em contracting} if it {\em converges to a quasiconstant map} 
$$
g_k\to \al_\om, 
$$
i.e.  
$$
g_k|_{\ol {X} - \{\om\}}\to \al$$
uniformly on compacts. The point $\al$ is the {\em limit point} (or the {\em attractor}) for the sequence $(g_k)$ and $\om$ is the {\em  exceptional point} (or the {\em repeller}) of the sequence.

\begin{rem}
If $(g_k)$ converges to $\al_\om$, then $(g_k^{-1})$ converges to $\om_\al$. 
\end{rem}

\begin{theorem}
[Convergence Property] 
\label{thm:conv}
Every divergent sequence $(g_k)$ in $G$ contains a contracting subsequence. 
\end{theorem}

While the naive generalization of the convergence property fails for the ideal boundaries of higher rank symmetric spaces and for flag manifolds, nevertheless, a certain version of the convergence property continues to hold, see section \ref{sec:conprop}.

The convergence at infinity of sequences in $X$ yields a notion of convergence at infinity for divergent 
sequences in $G$: 

\begin{definition}
A sequence $(g_k)$ in $G$ {\em converges} to a point $\al\in S$, $g_k\to \al$, if for some (equivalently, every) $x\in X$, 
$$
\lim_{k\to\infty} g_kx=\al. 
$$
\end{definition}

For instance, each contracting sequence converges to its attractor. 

The convergence property implies the following equivalent dynamical characterization of the convergence $g_k\to \al$
{i in terms of contraction. In particular, it yields a characterization in terms of the dynamics at infinity.}

\begin{lemma}\label{lem:flag-con-1}
For each $\al\in S$ and sequence $(g_k)$ in $G$ the following are equivalent:

1. $g_k\to \al$.

2. Every subsequence of $(g_k)$ contains a {contracting} subsequence which converges to $\al_\om$ for some $\om$. 

3. There exists a bounded sequence $(b_k)$ in $G$ such that the sequence $g_k b_k$ 
{is contracting and} converges to 
$\al_\om$ for some $\om$.

In addition, in parts 2 and 3, it suffices to verify the convergence to $\al_\om$  on $S - \{\om\}$. 
\end{lemma}

\subsection{Discrete subgroups}

\begin{definition}
A subgroup $\Ga < G$  is called {\em discrete} if it is a discrete subset of $G$. 
\end{definition}

\begin{definition} [Limit set] 
The {\em limit set} $\La(\Ga)$ of a discrete subgroup $\Ga < G$ is the  
accumulation set in $\geo X$ of one 
$\Gamma$-orbit $\Gamma x$ in $X$. 
\end{definition}

All orbits $\Ga x\subset X$ have the same accumulation set in $\geo X$.
This fact is an immediate application of the property that if $(x_i), (y_i)$ are two sequences in $X$ within bounded distance from each other and $x_i\to \xi\in S$, then $y_i\to \xi$ as well.

\begin{definition}
A discrete subgroup $\Ga< G$ is {\em elementary} if $\La(\Ga)$ consists of at most two points. 
\end{definition}

Every discrete subgroup $\Ga< G$  enjoys the convergence property, 
i.e.\ every divergent sequence $(\ga_k)$ in $\Ga$ contains a contracting subsequence,
compare Theorem~\ref{thm:conv}.
The convergence dynamics leads to a definition of the limit set in terms of the dynamics at infinity
and implies a dynamical decomposition of the $\Ga$-action into discontinuous and chaotic part:

\begin{lemma} For each discrete subgroup $\Ga< G$ we have:

1. $\La(\Ga)$ is the set of values $\al$ of attractors of contracting sequences of elements of $\Gamma$.  

2. $\La(\Ga)$ is the set of exceptional points of contracting sequences of elements of $\Gamma$. 

3. Unless $\Ga$ is elementary, its action on $\La(\Ga)$ is {\em minimal}: Every $\Ga$-orbit in $\La$ is dense.  

4. The domain $\Om(\Ga):=  S - \La(\Gamma)$ equals the {\em wandering set} of the action $\Ga\acts S$.

5. The action $\Ga\acts X\cup \Om(\Ga)$ is properly discontinuous. 
\end{lemma}

The subset $\Omega(\Gamma) = S - \La(\Gamma)$ is called the  {\em domain of discontinuity} of the 
subgroup $\Ga< G$. 

\subsection{Conical convergence}\label{sec:conical_convergence}

The notion of conical convergence plays a central role in describing geometric finiteness (both in rank 1 and in  higher rank), 
so we define it below in several different ways.

\begin{definition}
A sequence $x_k\in X$ converges to a point $\xi\in S$ {\em  conically} 
$$
x_k \stackrel{con}{\lra} \xi
$$
if $x_k\to \xi$ and for every ray $x\xi$ there is $R<\infty$ such that $x_k\in N_R(x\xi)$, the $R$-neighborhood of $x\xi$. A sequence $g_k\in G$ converges to a point $\xi\in S$ {\em conically}
$$
g_k \stackrel{con}{\lra} \xi, 
$$
if for some (equivalently, every) $x\in X$, the sequence $x_k=g_k x$ converges to $\xi$ conically. \end{definition}

The name {\em conical} in this definition 
comes from the fact that, e.g. in the upper half space model of $\H^n$, tubular neighborhoods of rays 
resembles  cones. 

\medskip 
As in the case of convergence $g_k\to \xi$  one can also characterize conical convergence 
in terms of the dynamics of $(g_k)$ on $S$. 

\begin{lem}\label{lem:eq-con}
Suppose that $g_k\to \xi\in S$, $g_k\in G$. Then $g_k \stackrel{con}{\lra} \xi$ 
if and only if either one of the two equivalent properties hold: 

1. For some, equivalently, every complete geodesic $l\subset X$ asymptotic to $\xi$, 
the sequence $g^{-1}_k l$ is relatively compact in the space of all geodesics in $X$. 

2.  For some, equivalently, every $\hat\xi\in S - \{\xi\}$ the sequence $g^{-1}_k(\xi, \hat\xi)$ is 
relatively compact in 
$$
(S\times S)^{opp}= S \times S - Diag(S\times S). 
$$
\end{lem} 

\subsection{The expansion property}\label{sec:expansion}

The notion of conical convergence $g_k \stackrel{con}{\lra} \xi$ is closely related to the concept of {\em expansivity}. We will use the expansivity concepts discussed in Appendix \ref{sec:expanding_actions}  
for actions on  $S$ equipped with an arbitrary Riemannian metric. The choice of the metric will not be important. 

\begin{prop}\label{prop:con-exp}
Suppose that $g_k\to \xi\in S$, $g_k\in G$. Then the conical convergence
$$
g_k \stackrel{con}{\lra} \xi 
$$
implies that the sequence $(g_k^{-1})$ has  diverging 
{ infinitesimal}
expansion at $\xi$. 
\end{prop}

\subsection{Conical limit points} 

We return to discussing discrete subgroups $\Ga< G$. 

\begin{definition}
[Conical limit points of $\Gamma$] A limit point $\xi$ of $\Ga$ is {\em conical} if there exists 
a sequence $\ga_k\in \Ga$ which converges to $\xi$ conically. The set of conical limit points of $\Ga$ is denoted 
$\La_c(\Ga)$. 
\end{definition}

In view of Lemma \ref{lem:eq-con},
one has the following characterization of conical limit points in terms of the dynamics at infinity:

\begin{lem}
Suppose that the limit set of $\Ga$ consists of at least two points. Then the following are equivalent:

1. $\xi\in \La(\Ga)$ is a conical limit point of $\Ga$. 

2. There exists a sequence $\ga_k\to\xi$ in $\Ga$ such that 
for some
$\hat\xi\in \La(\Ga) - \{\xi\}$ 
the sequence of pairs
$\ga_k^{-1}(\xi,\hat\xi)$ converges to a pair  of distinct points.  
\end{lem}

\medskip 
The situation when {\em all limit points are conical}
can also be characterized in terms of the action on the space $T\La$ of triples of distinct points:

\begin{theorem}
[See e.g. \cite{Bowditch_config}] \label{thm:traction}
Suppose that $\Ga$ is nonelementary. Then all limit points are conical iff the action of $\Ga$ on the triple space $T\La$ is cocompact. 
\end{theorem}

{ The triple space is an intrinsic replacement for the convex hull of the limit set in the symmetric space,
and the theorem provides one of the characterizations of {\em convex cocompactness} 
to be discussed in the next section, 
compare Theorem~\ref{thm:coco-GF2} below.}

\subsection{Geometrically finite groups}\label{sec:GFG}

The notion of geometric finiteness played a critical role in the development of the theory of Kleinian groups. It was 
originally introduced by Ahlfors, who defined geometric finiteness in terms of fundamental polyhedra. Subsequent equivalent definitions were established by Marden, Beardon and Maskit, Thurston, Sullivan and others. 

In this section we will give a list of equivalent definitions of convex-cocompactness in the rank 1 setting (equivalently, geometric finiteness without parabolic elements). In what follows, we will only consider discrete subgroups $\Ga < G$ of rank 1 Lie groups which contain no parabolic elements: These definitions require modifications if one allows parabolic elements, we refer the reader to \cite{Bowditch93, Bowditch_gf, Ratcliffe} for more details.

\subsubsection{Finitely sided fundamental domains}

\begin{definition}
[L.~Ahlfors] 
$\Ga<G$ is CC0 if for some point $o\in X$ not fixed by any nontrivial element of $\Ga$ the associated  
{\em Dirichlet fundamental domain} $D_o$ of $\Ga$, 
$$
D_o=\{x\in X: \forall \ga\in \Ga, d(x, o)\le d(x, \ga o)\},
$$
is {\em finite-sided}. The latter means that only finitely many ``half-spaces''  
$$
{\mathrm Bis}(\ga o, o)= \{x\in X: d(x, o)\ge d(x, \ga o)\}, \ga\in \Ga, 
$$
have nonempty intersection with $D_o$. 
\end{definition}

This definition was proposed by L.~Ahlfors in \cite{Ahlfors}; it was historically the first definition of geometric finiteness 
and the main one, until Thurston's work \cite{Thurston}.

\subsubsection{Convex cocompactness} 

\begin{definition}
[Convex cocompact subgroups]\label{def:cc1}
 $\Gamma<G$ is CC1 ({\em convex cocompact}) if there exists a nonempty $\Ga$-invariant closed convex subset $C\subset X$ such that $C/\Ga$ is compact. 
\end{definition}

This definition explains the terminology ``convex cocompact'' since it is stated in terms of cocompactness of the 
$\Gamma$-action on a certain convex subset of $X$.

There is a unique smallest nonempty $\Ga$-invariant closed convex subset 
if $|\La(\Ga)|\ge 2$,
namely the {\em convex hull} $C_\Ga$ of $\La(\Ga)$, which is the closed convex hull of the union of all 
geodesics connecting limit points of $\Ga$, see e.g.  \cite{Bowditch_gf}.\footnote{A convex subset $C\subset X$
as in Definition~\ref{def:cc1} contains $\Ga$-orbits.
Hence $\La(\Ga)\subseteq\geo C$, and therefore $C_\Ga\subseteq C$.}
Hence, to verify CC1, one  needs to test only $C_\Ga$:

\begin{lemma}
Assume that $|\La(\Ga)|\ge 2$. Then $\Ga$ is convex cocompact iff $C_\Ga/\Ga$ is compact.  
\end{lemma}

\medskip
Definitions CC0 and CC1 do not appear to be particularly useful in higher rank; below we present definitions which,
except for CC8,
do generalize to higher rank (after suitable modifications). 

\subsubsection{Beardon--Maskit condition: Dynamics on the limit set} 

The next definition is motivated by the work of Beardon and Maskit \cite{Beardon-Maskit}
who characterized the discrete subgroups of $PSL(2,\C)$ satisfying Ahlfors' CC0 condition 
in terms of their dynamics on the limit set.

\begin{definition}[A.~Beardon, B.~Maskit] 
$\Ga<G$ is CC2 if each limit point of $\Ga$ is conical. 
\end{definition}

Theorem \ref{thm:traction} can be reformulated as:  

\begin{thm}\label{thm:coco-GF2}
A nonelementary group $\Ga$ is CC2 iff $\Ga$ acts cocompactly on $T\La(\Ga)$. 
\end{thm}

\begin{rem}
 In the presence of parabolics one requires that each limit point is either conical or a ``bounded'' parabolic fixed point (A. Beardon, B. Maskit, B. Bowditch, see \cite{Beardon-Maskit}, \cite{Bowditch_gf}; cf. also \cite{Bishop}). 
\end{rem}

\medskip
Note that the condition CC2 {\em a priori} does not even imply finite generation of $\Ga$.

\subsubsection{Asymptotically embedded groups} 

Recall that each word hyperbolic group $\Ga$ has a {\em Gromov boundary} 
$\geo \Ga$, which is a metrizable compact on which $\Ga$ acts via homeomorphisms. (One constructs this boundary by looking at equivalence classes of geodesic rays in the Cayley graph of $\Ga$ or via horofunctions,
see \cite{CP}.)  

\begin{definition} [Asymptotically embedded] 
$\Ga< G$ is CC3 ({\em asymptotically embedded}) if it is Gromov-hyperbolic and 
$\geo \Ga$ is equivariantly homeomorphic to $\La(\Ga)$. 
\end{definition}

Equivalently:

\begin{definition} [Boundary embedded] 
$\Ga< G$ is {\em boundary embedded} if it is Gromov-hyperbolic and 
there exists an equivariant topological embedding $\be: \geo \Ga \to \geo X$.  
\end{definition}

The equivalence of CC3 and {\em boundary embedded} is easy to see using again the convergence property; 
it is also easy to see that $\be(\geo \Ga)=\La(\Ga)$.

\subsubsection{Coarse geometric definitions}

\medskip 
The next definition involves the {\em coarse geometry} of discrete subgroups: 

\begin{definition} 
$\Ga< G$ is CC4 if it is finitely generated and {\em undistorted} in $G$.
\end{definition}

Here $\Ga< G$ is {\em undistorted} if the word metric on $\Ga$ is comparable to the extrinsic metric coming from $G$. Equivalently, one (equivalently, each) orbit map $\Ga \to \Ga x\subset X$ is a QI (quasiisometric) embedding of $\Ga$ into $X$. 

\medskip
A minor variation on this definition (which will become major in higher rank) is:

\begin{definition}
A discrete subgroup $\Ga< G$ is {\em Morse}, or {\em satisfies the Morse property}, if $\Ga$ is 
word hyperbolic and each discrete geodesic in $\Ga$ maps (via the orbit map) to a 
discrete path in $X$ uniformly close to a geodesic.
\end{definition} 

Note that this definition does not a priori assume undistortion of $\Ga$ in $G$. 

\medskip 
The implication CC4 $\Rightarrow$ Morse follows immediately from the Morse Lemma. For the converse implication one observes that images of discrete geodesics in $\Ga$ under the orbit map are contained in uniform neighborhoods of geodesics in $X$ and have bounded backtracking.

\subsubsection{Quasiconvexity and coarse retractions} 

A subset $Y\subset X$ is called {\em quasiconvex} if there exists a constant $R$ such that for any pair of points $x, y\in Y$ the (unique) geodesic $xy$ between $x$ and $y$ is contained in $N_R(Y)$, 
the $R$-neighborhood of $Y$ in $X$. Each convex subset is, of course, also quasiconvex. While the opposite implication is false, it follows from the work of M.~Anderson that each quasiconvex subset $Y\subset X$ is within finite Hausdorff distance from its convex hull in $X$ (see \cite[Proposition 2.3.4]{Bowditch_gf}). 

\begin{definition}
[Quasiconvex subgroups] 
$\Ga < G$ satisfies CC5 if it  is {\em quasiconvex}, i.e., one (equivalently, every) orbit $\Ga x\subset X$ is a quasiconvex subset. 
\end{definition}

For each nonempty closed convex subset $C\subset X$ the nearest point projection $X\to C$ is $1$-Lipschitz, i.e., is distance non-increasing. Similarly, if $Y$ is a quasiconvex subset of a geodesic Gromov-hyperbolic space then there exists an $(L,A)$ coarse Lipschitz retraction $X\to Y$, which can be viewed as a coarsification of the nearest point projection. (A nearest point in $Y$ may not exist, instead, one projects $x\in X$ to $y\in Y$ such that for all $y'\in Y$, $d(x,y')\le d(x,y)+1$.) Here a map $f: X\to Y$ between metric spaces is $(L,A)$ coarse Lipschitz if
$$
d(f(x_1), f(x_2))\le L d(x_1, x_2) + A, \forall x_1, x_2\in X.  
$$

A subset $Z\subset X$ is called a 
{\em coarse retract} if there exists a coarse Lipschitz retraction $X\to Z$. 

\begin{definition}
[Coarse retract] 
\label{defn:retract}
A finitely generated subgroup $\Ga < G$ is a {\em coarse retract} if for one (equivalently, every) $x\in X$ there exists a coarse Lipschitz map $r: X\to \Ga$ such that the composition 
$$
\ga \mapsto \ga x \stackrel{r}{\mapsto} \ga'\in \Ga, 
$$
is within finite distance from the identity map. Here we equip $\Ga$ with a word metric. 
\end{definition}

{ 
\begin{rem}
This definition makes sense, of course, not only for negatively curved symmetric spaces but for all nonpositively curved symmetric spaces $X$, where $G$ is the identity component of the isometry group of $X$. 
\end{rem}
}

In view of the Morse Lemma and the coarse Lipschitz property of nearest point projections to quasiconvex subsets of $X$, one obtains: 

\begin{thm}\label{thm:CC4=CC5}
A finitely generated discrete subgroup $\Ga< G$ is undistorted iff it is quasiconvex iff it is a coarse retract. 
\end{thm}

\subsubsection{Expanding actions}

We refer the reader to Appendix \ref{sec:expanding_actions} for definitions of 
metric expansion and infinitesimal expansion.

\begin{definition}
[Expanding subgroups,  D.~Sullivan, \cite{Sullivan}] A discrete subgroup $\Ga< G$ is CC6 ({\em expanding})  
 if for each $\xi\in \La(\Ga)$ there exists  $\ga\in \Ga$ which is metrically expanding on $S$ at $\xi$. 
\end{definition}

Below are two variations on the expansion axiom:

\begin{theorem}\label{thm:CC2=CC6}
The following are equivalent:

1. $\Ga$ is {\em infinitesimally expanding} at $\La(\Ga)$:  For each $\xi\in \La(\Ga)$ there exists  $\ga\in \Ga$ 
which is infinitesimally expanding at $\xi\in S$.

2. $\Ga< G$ is CC6 (expanding). 

3. $\Ga$ is nonelementary and the action of $\Ga$ is metrically 
expanding on $\La(\Ga)$ (i.e., it suffices to check the expansion of distances only between limit points).

4. The group $\Ga$ is CC2. 
\end{theorem}
\proof It is clear that $1\Rightarrow 2 \Rightarrow 3$. The implication $3\Rightarrow 4$ is proven in Theorem \ref{thm:conical} in the Appendix \ref{app:congru}.  Lastly, the implication $4\Rightarrow 1$ follows from {\em extrinsic} conicality of the limit points of $\Ga$ (Lemma \ref{lem:eq-con}) and Proposition \ref{prop:con-exp}. \qed

\medskip 
The advantage of CC6 and its variations is that they make sense for general topological/smooth dynamical systems and, hence, are easy to extend to higher rank.

\subsubsection{Natural compactification of locally symmetric space} 

Our next definition is formulated in terms of existence of a natural compactification of the locally symmetric space  
$X/\Ga$: 

\begin{definition} 
[A.~Marden, \cite{Marden}] $\Ga$ is CC7 if the space $(X\cup \Om(\Ga))/\Ga$ is compact. 
\end{definition}

{ This definition first appeared in Marden's paper \cite{Marden} where he proved its equivalence to CC0 in the case 
of $X=\H^3$.}

\begin{figure}[tbh]
\centerline{\epsfxsize=4.5in \epsfbox{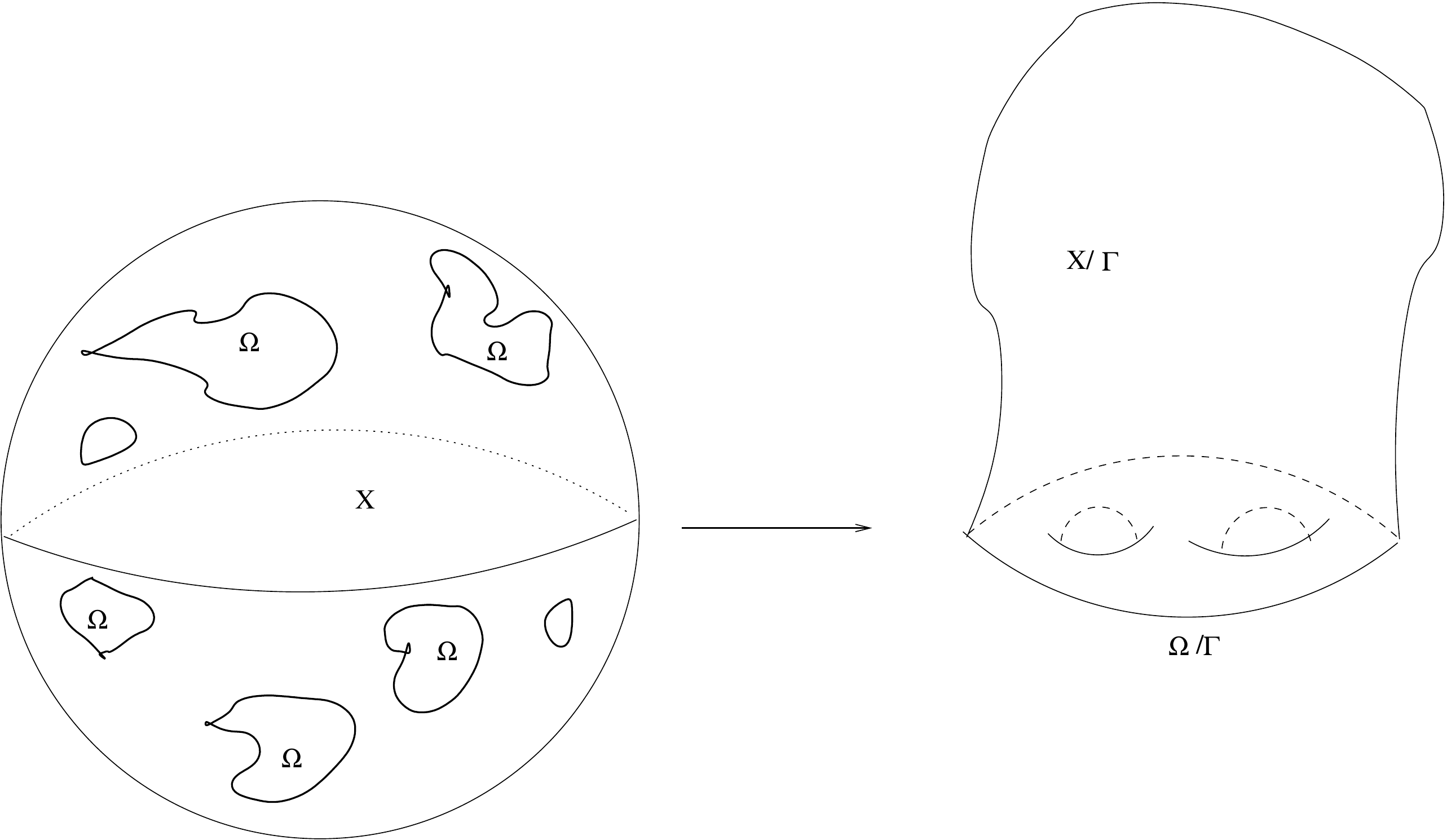}}
\caption{Quotient space of a geometrically finite group.}
\label{figure1.fig}
\end{figure}

\subsubsection{Finiteness of volume} 

The last definition states geometric finiteness in terms of the volume of the quotient space:

\begin{definition}
[W.~Thurston; B.~Bowditch] A discrete subgroup $\Ga< G$ is CC8 if either $|\La(\Ga)|\le 1$ or $|\La(\Ga)|\ge 2$ and:

1. The orders of the torsion elements of $\Ga$ are  bounded. 

2. For some (every) $\eps>0$ the quotient $N_\eps(C_\Ga)/\Ga$ has finite volume. 
\end{definition}

Here $C_\Ga$ is, as before, the closed convex hull of the limit set of $\Ga$ and $N_\eps$ is the $\epsilon$-neighborhood of $C_\Ga$ in $X$. 

\begin{rem}
{This definition is mostly due to W.~Thurston \cite{Thurston} who stated it for isometry groups of the hyperbolic 3-space 
without the extra conditions on torsion elements. The latter assumption was added by B.~Bowditch in the general setting.  The restriction on  orders of torsion elements is essential, unless $X$ is the (real) hyperbolic space of dimension $\le 3$ (E.~Hamilton, \cite{Hamilton}).} 
\end{rem}

\subsubsection{An equivalence theorem}

The following is a combination of work of many people:

\begin{theorem}\label{thm:main1}
For discrete isometry groups of rank 1 symmetric spaces (without parabolic elements), 
all the conditions CC1---CC8 are equivalent. 
\end{theorem}
\proof {
The equivalence of conditions CC1, CC2, CC7, CC8 is in Bowditch's paper \cite{Bowditch_gf}; note that Bowditch proved this result for discrete isometry groups of negatively pinched Riemannian manifolds, not just symmetric spaces of negative curvature. The equivalence of CC2 and CC6 is Theorem \ref{thm:CC2=CC6}. The equivalence of 
CC4 and CC5 is Theorem \ref{thm:CC4=CC5}. If $\Gamma$ is CC5 then the convex hull of $\Gamma x\subset X$ is Hausdorff-close to $\Gamma x$, hence, $\Gamma$ is CC1. If $\Gamma$ is CC1 then, taking $x\in C$ (as in the definition of CC1), and taking into account compactness of $C/\Gamma$, we conclude that  $\Gamma$ is CC5. 
Assume that $\Gamma$ is asymptotically embedded (CC3). Then $\Gamma$ is Gromov-hyperbolic and every $\xi\in \geo \Gamma$ is a conical limit point, see \cite{Tukia1994}. Hence, $\Gamma$ is CC2. Assume that $\Gamma$ is convex-cocompact (CC1) and acts cocompactly on the closed convex subset $C=C_\Ga\subset X$, the convex hull of the limit set of $\Gamma$. Then $C$ is a Gromov-hyperbolic  geodesic metric space quasiisometric to $\Gamma$. Hence  $\Gamma$ is Gromov-hyperbolic; the ideal boundary of $\Gamma$ is naturally homeomorphic to the ideal boundary of $C$, i.e. the limit set of 
$\Gamma$. Hence, $\Gamma$ is asymptotically embedded. \qed 
}

{ 
\begin{rem}
The equivalence of CC0 and CC1 in the case of the real hyperbolic spaces is proven in \cite{Bowditch93} and \cite[Theorem 12.4.5]{Ratcliffe}. Their proofs rely upon convexity of Dirichlet domains. While Dirichlet domains 
for general negatively curved symmetric spaces are not convex, they are {\em quasiconvex} which can be used to extend 
the arguments of \cite{Bowditch93} and \cite[Theorem 12.4.5]{Ratcliffe} to this more general setting. \end{rem}

}

\subsection{Consequences of geometric finiteness}\label{sec:corGF}

The first consequence of geometric finiteness is immediate:

\begin{corollary}
[C1] For a convex cocompact subgroup $\Ga< G$, the quotient $\Om(\Ga)/\Ga$ is compact. 
\end{corollary}

The next theorem, known as the {\em structural stability property} was first proven by D.~Sullivan  
\cite{Sullivan} using methods of symbolic dynamics and later by C. Yue \cite{Yue} using Anosov flows. 

\begin{thm}
[C2] Convex cocompactness implies structural stability: If $\Ga< G$ is convex cocompact then any homomorphism $\rho: \Ga\to G$ close to $id: \Ga\embed G$ is injective and $\rho(\Ga)< G$ is a 
convex cocompact subgroup which is topologically conjugate to $\Ga$ on the limit set: There exists a $\rho$-equivariant homeomorphism
$$
h_\rho: \La(\Ga)\to \La(\rho(\Ga)). 
$$
\end{thm}

Moreover:

\begin{thm}
[C3] In the context of C2: 

a. There exists a $\rho$-equivariant topological conjugation $f_\rho: \ol {X}\to \ol {X}$, 
which is smooth away from the limit set. 

b. If a sequence of representations $\rho_i$ converges to the identity representations, then the maps $f_{\rho_i}$ can be chosen so that 
$$
\lim_{i\to\infty} f_{\rho_i} =id.$$
Here convergence is uniform on $\ol {X}$ and $C^\infty$-uniform on compacts 
in the complement to the limit set.  
\end{thm}

This stronger stability theorem is a result of combined efforts of many people, see  \cite{Bowditch-stab, Izeki}.\footnote{Bowditch and Izeki only consider the case of the real-hyperbolic space but the proofs go through for other negatively curved symmetric spaces as well.} 
 
\begin{thm}
[C4] Convex cocompactness is {\em semidecidable}. 
\end{thm}

Recall that an algorithmic problem is {\em semidecidable} if there is an algorithm which answers 
YES in finite time if and only if the answer is positive (and runs forever if the answer is negative). 
Since we are dealing with computations over the reals, one has to specify the computability model: Here and below we are using the BSS (Blum-Shub-Smale), also known as the Real RAM, computability model. See \cite{BCSS} for the details. 

\medskip 
There are two ways to interpret the semidecidability of convex cocompactness.

\begin{thm}
\cite{morse}. Suppose that $\Ga$ is a word hyperbolic group defined in terms of a finite presentation. Then  there is an algorithm which, given a representation $\rho: \Ga\to G$ (defined in terms of the images of the generators) will terminate with the positive answer if and only if $\rho$ has finite kernel and the image $\rho(\Ga)$ is convex cocompact.   
\end{thm}

The first written proof of this theorem seems to be in \cite{morse} (in the context of Morse actions of hyperbolic groups on higher rank symmetric spaces), although some special cases of this theorem might have been known earlier.

One drawback of the above semidecidability theorem is that we are required to know in advance which hyperbolic group is being represented. The following theorem is limited (for various reasons) to  hyperbolic 3-space, but does not require a priori knowledge of the algebraic structure of $\Ga$;
the algorithm appears to be first discussed (and implemented) by R.~Riley, see \cite{Riley}; see also {a} paper by J.~Gilman \cite{Gilman2} and {one by} J.~Manning \cite{Manning}.

\begin{thm}
Geometric finiteness\footnote{Here we allow parabolic elements.} is semidecidable for subgroups of 
the isometry group of  hyperbolic 3-space $G=Isom(\H^3)$. 
\end{thm}
\proof The proof is in the form of a {\em ``Poincar\'e'' algorithm} for constructing a finite sided Dirichlet domain for discrete subgroups of $G$. 

\begin{figure}[tbh]
\includegraphics[width=90mm]{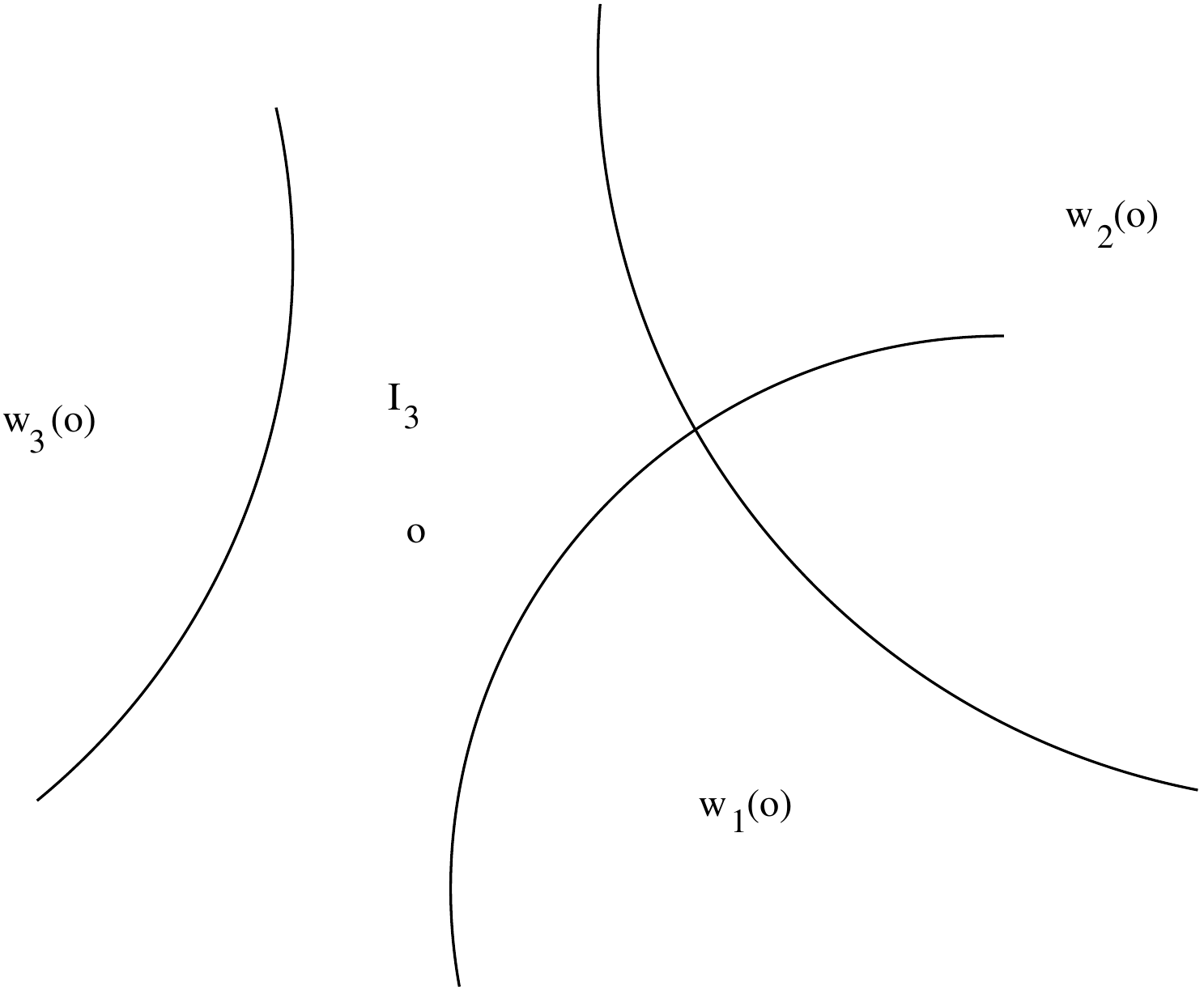}
\caption{Constructing a Dirichlet fundamental domain.}
\label{figure2.fig}
\end{figure}

The input for the (semi)algorithm is a tuple $(g_1,...,g_n)$ of elements of $G$. It attempts to construct a finite sided Dirichlet fundamental domain of the group $\Ga$ generated by $g_1,...,g_n$ by computing, inductively, intersections $I_k$ in $\H^3$ of half-spaces bounded by bisectors of pairs $o, w_i(o)$, where the $w_i$ are reduced words in $g_1^{\pm 1},...,g_n^{\pm 1}$, 
$$
I_k= \bigcap_{i=1}^k {\mathrm Bis}(o, w_i o),
$$
where
$$
{\mathrm Bis}(o, w o)= \{x\in \H^3: d(o, x)\le d(x, w o)\}. 
$$
See Figure \ref{figure2.fig}. (There is a separate issue of making sure that $o\in \H^3$ is not fixed by a nontrivial element of $\Ga$, we will not  address this problem here.)  The sequence $(w_i)$ is chosen to exhaust the free group on the generating set $g_1,...,g_n$. After constructing $I_k$ (by solving a system of linear inequalities in the Lorentzian space $\R^{3,1}$), the algorithm 
checks if the conditions of   Poincare's Fundamental domain theorem (see \cite{Maskit, Ratcliffe}) are satisfied by $I_k$. If they are satisfied for some $k$, then $\Ga$ admits a finite sided Dirichlet domain, namely, $I_k$. If $\Ga$ is geometrically finite, then this algorithm terminates (for any choice of  base point). If $\Ga$ is not geometrically finite, this algorithm will run forever. \qed

\medskip
Note that if $\Ga$ is geometrically finite, one can read off a finite presentation of $\Ga$ from 
the fundamental domain.

\begin{rem}
Note that for 2-generator subgroups of $PSL(2, \R)$ there are (more efficient) alternatives to Riley's algorithm, due to J.~Gilman and B.~Maskit, see \cite{Gilman1, Gilman-Maskit} and also \cite{Gilman2} for comparison of the computational complexities.
\end{rem}

\medskip
As an aside, we discuss the question of decidability of discreteness for finitely generated subgroups of connected 
Lie groups. For instance, one can ask if 
a representation $F_n\to G$ (where $G$ is a connected algebraic Lie group) has  discrete image. 
{
First of all, one has to eliminate certain classes of representations, otherwise the  discreteness problem 
is undecidable already for $n=1$ and $G=U(1)$, cf. \cite{Kapovich2015}. 

\begin{definition}
Let ${\mathcal E}\subset Hom(F_n, G)$ denote the subset consisting of representations $\rho$ such that 
$\rho(F_n)$ contains a nontrivial normal nilpotent subgroup. 
\end{definition}

For instance, in the case $G=PSL(2,\C)$, a representation $\rho$ belongs to ${\mathcal E}$ if and only if the group 
$\Ga= \rho(F_n)$ either has a fixed point in $\C P^1$ or preserves a 2-point subset of $\C P^1$. Algebraically, (for subgroups of $PSL(2,\C)$) this is equivalent to the condition that $\Ga$ is solvable. 

Secondly, one has to specify the {\em computability model}; as before we use the BSS (Blum--Schub--Smale) computability model. Restricting to representations in $Hom(F_n, G) \setminus {\mathcal E}$, one obtains the following 
folklore theorem, see Gilman's papers \cite{Gilman1, Gilman2} in the case $G=PSL(2,\C)$:  

\begin{thm}
For a connected algebraic Lie group $G$, it is semidecidable whether a representation 
$\rho\in Hom(F_n, G) \setminus {\mathcal E}$ is  nondiscrete. 
\end{thm}
\proof The key to the proof is a theorem of Zassenhaus, see e.g. \cite{Raghunathan,  Kapovich00}, where 
we regard $G$ as a real algebraic subgroup of $GL(k, \C)$ for some $k$, which we equip with the standard operator norm: 

\begin{theorem}
[H. Zassenhaus] There exists a (computable) number $\epsilon$ such that the neighborhood $U= G\cap B(1, \epsilon)$ (called a {\em Zassenhaus neighborhood} of $1$ in $G$) satisfies the following property: Whenever $\Gamma< G$ is 
a  subgroup, the subgroup $\Gamma_U< \Gamma$ generated by $\Gamma\cap U$ is either nondiscrete or nilpotent.  
\end{theorem} 

Suppose that $\Ga< G$ is nondiscrete, $\overline{\Ga}^0$ is the identity component of $1\in \overline{\Gamma}$, the 
closure of $\Gamma$ in $G$ (with respect to the standard matrix topology). Then 
$\overline{\Ga}^0$ is a normal subgroup of $\overline{\Gamma}$ of positive dimension. Therefore, the intersection 
$\overline{\Ga}^0\cap U$ is nondiscrete and, hence, the subgroup $\Gamma_U$ is nondiscrete as well.  

There are two cases which may occur:

1. The subgroup $N=\overline{\Ga}^0$  is nilpotent. Since $N$ is a Lie subgroup of $G$, there exists a  neighborhood $V\subset U$ of $1\in G$ such that $\Gamma \cap V$ is contained in $N$. In particular, $\Gamma$ contains a nontrivial normal nilpotent subgroup, namely $\Gamma\cap N$. This cannot happen if $\Ga=\rho(F_n)$, 
 $\rho\in Hom(F_n, G)\setminus {\mathcal E}$. 

2. The subgroup  $N=\overline{\Ga}^0$  is not nilpotent. Note that by Lie--Kolchin theorem, 
every connected nilpotent subgroup of $GL(k, \C)$ is conjugate to the group of upper triangular matrices. 
In particular, a connected Lie subgroup of $G$ is nilpotent if and only if it is (at most) $k-1$-step nilpotent. 
Thus, in our case, there exist elements $g_1,...,g_{k}\in N\cap U$ 
such that the $k$-fold iterated commutator
$$
[... [[g_1, g_2], g_3],...,g_k] 
$$
is not equal to $1\in G$. By continuity of the $k$-fold commutator map, there exist 
$\ga_1,..., \ga_{k}\in \Ga\cap U$ such that 
$$
[... [[\ga_1, \ga_2], \ga_3],..., \ga_k] \ne 1.  
$$

We now  describe our (semi)algorithm: We enumerate $k$-tuples of elements 
$(x_1,..., x_k)\in F_n\times ... \times F_n$ and, given $\rho\in Hom(F_n, G)\setminus {\mathcal E}$, 
look for the tuples such that 
$$
\ga_i=\rho(x_i), i=1,...,k
$$ 
satisfy $\ga_i\in U,  i=1,...,k$ and 
$$
[... [[\ga_1, \ga_2], \ga_3],..., \ga_k] \ne 1.  
$$
If $\Ga=\rho(F_n)$ is nondiscrete then we eventually find such a tuple thereby verifying nondiscreteness of $\Gamma$. \qed 
}

\medskip
In the case when $G=PSL(2,\R)$, a finitely generated subgroup of $G$ is discrete if and only if it is geometrically finite. Therefore, one can use Riley's algorithm in combination with nondiscreteness algorithm to determine if an $n$-generated nonsolvable  subgroup of $PSL(2,\R)$ is discrete. Hence, discreteness is decidable in $PSL(2,\R)$. On the other hand: 

\begin{thm}
\cite{Kapovich2015}.  Being discrete is undecidable  for nonsolvable 2-generated subgroups in $PSL(2,\C)$.  
 \end{thm}

\bigskip 
\section{Geometry of symmetric spaces of noncompact type}\label{sec:2}

\subsection{Basic geometry}\label{sec:Basicgeometry}

We refer to \cite{Eberlein, BGS} and \cite{Helgason} for a detailed treatment of symmetric spaces. 
From now on, $X$ is a symmetric space of noncompact type: It is a nonpositively curved symmetric space without euclidean factor, $G=Isom_o(X)$ is the identity component of the full isometry group of $X$. We will use the notation $xy$ for (oriented) geodesic segments in $X$ connecting $x$ to $y$. Recall that each symmetric space admits a {\em Cartan involution} $s_x$ about every point $x\in X$; such $s_x$ fixes $x$ and acts as $-id$ on the tangent space $T_xX$.

Then $X\cong G/K$, where $K<G$ is a maximal compact subgroup (the stabilizer in $G$ of a basepoint $o$ in $X$); $G$ is a {\em semisimple} real Lie group; the Riemannian metric on $X$ is essentially uniquely determined by $G$ (up to rescaling for each simple factor of $G$). An important example is 
$$
G=PSL(n,\R), K=PSO(n);
$$
the symmetric space $X=G/K$ can in this case be identified with the projectivized space of positive definite bilinear forms on $\R^n$.

\begin{rem}
For our examples we will frequently use $SL(n)$ instead of $PSL(n)$. The difference is that the group $SL(n)$ acts on the associated symmetric space with finite kernel. 
\end{rem}

A symmetric space $X$ is {\em reducible} if it metrically splits as a product $X_1\times X_2$. 
Each symmetric space $X$ of noncompact type admits a canonical (up to permutation of factors) product decomposition
$$
X= X_1\times \ldots \times X_n
$$
into irreducible symmetric spaces. 

\medskip 
{\bf Classification of isometries.} For each $g\in G$, as in rank 1, we consider its convex displacement function
$$
d_g(x)=d(x, gx). 
$$

\begin{definition}
The isometries of $X$ are classified as follows: 

1. An isometry $g$ of $X$ is {\em axial} or {\em hyperbolic} if $\inf_{x\in X}d_g>0$ and the infimum is realized. 
In this case, there exists a $g$-invariant geodesic in $X$, an {\em axis of $g$}, 
along which $g$ translates.\footnote{In general, this axis is not unique, but all axes of $g$ are parallel.} 
The union of axes is the minimum set of the convex function $d_g$. 

2. $g$ is {\em mixed} if $\inf_{x\in X}d_g>0$ but the infimum is not realized. 

3. $g$ is {\em parabolic} if $\inf_{x\in X}d_g=0$ but the infimum is not realized. 

4. $g$ is {\em elliptic} if $d_g=0$ and the infimum is realized. Equivalently, $g$ has a fixed point in $X$. 
\end{definition}

\medskip
An axial isometry $g$ is a {\em transvection} if it preserves parallel vector fields along one and hence any  
axis. Equivalently, $g$ is a product of two different Cartan involutions. 
A parabolic isometry $g$ is {\em  unipotent} 
if the closure of its conjugacy class contains the neutral element,
i.e.\ if there is a sequence $h_k\in G$ such that  
$$
\lim_{k\to\infty} h_k g h_k^{-1}=e\in G. 
$$

\medskip

\begin{definition} 
A  {\em flat} in  $X$ is a (totally geodesic) isometrically embedded euclidean subspace in $X$. A {\em maximal flat} is a flat in $X$ which is not properly contained in a larger flat. 
\end{definition}

\begin{ff}
All maximal flats in $X$ are $G$-congruent. 
\end{ff}

\begin{defn}
 $r=\rank(X)$ is the dimension of a maximal flat. 
\end{defn}

Note that $r=\rank_{\R}(G)$. 

\begin{ff}
[Cartan decomposition] 
$$G=K A_+ K ,$$ 
where $A < G$ is a {\em Cartan subgroup} (equivalently, a maximal abelian group of transvections, equivalently, a maximal $\R$-split torus), and $A_+\subset A$ is a certain sharp closed convex cone  with tip at $e$ (a subsemigroup). 
\end{ff}

More precisely,
the unique maximal flat $F\subset X$ preserved by $A$ contains the fixed point $o\in X$ of $K$,
$o\in F$. 
The cone $V=A_+o\subset F$ is a {\em euclidean Weyl chamber} 
with tip at $o$.
The Cartan decomposition corresponds to the fact that every $K$-orbit in $X$ intersects $V$ in precisely one point. 

\begin{example}
\label{ex:cartandeco}
For $G=SL(n,\R)$ and $K=SO(n)$, 
the Cartan subgroup $A< G$ can be chosen as the subgroup of diagonal matrices with positive entries,
and $A_+\subset A$ as the subset of diagonal matrices $a=diag(a_1,...,a_n)$ 
with decreasing  diagonal entries: 
$$a_1\ge a_2 \ge  ... \ge a_n>0$$
 The Cartan decomposition in this case is also known as the {\em singular value decomposition} of a matrix: $g= uav$ with $u, v\in SO(n)$ and $a\in A_+$. The diagonal entries $a_1,...,a_n$ of the matrix $a$ are known as the {\em singular values} of the matrix $g$. 
\end{example}

The $G$-stabilizer $G_F$ of a maximal flat $F\subset X$ acts on $F$ (in general unfaithfully); the image of the restriction homomorphism 
$$
G_F\to \Isom(F) 
$$
is a semidirect product
$$
W_{aff}= \R^r \rtimes W,
$$
where $\R^r\cong A$ is the full group of translations of $F$ and $W$ is a certain finite reflection group of isometries of $F$, called {\em the Weyl group} of $X$ (and of $G$). In view of the $G$-congruence of maximal flats, the action 
$W\acts F$ is independent of the choices. 

\begin{rem}
The subgroup $\R^r$ lifts to $G_F$ as the group of transvections (in $X$) along $F$. 
In contrast, the subgroup $W$ does in general not lift  to $G_F$. 
\end{rem}

\begin{figure}[tbh]
\includegraphics[width=90mm]{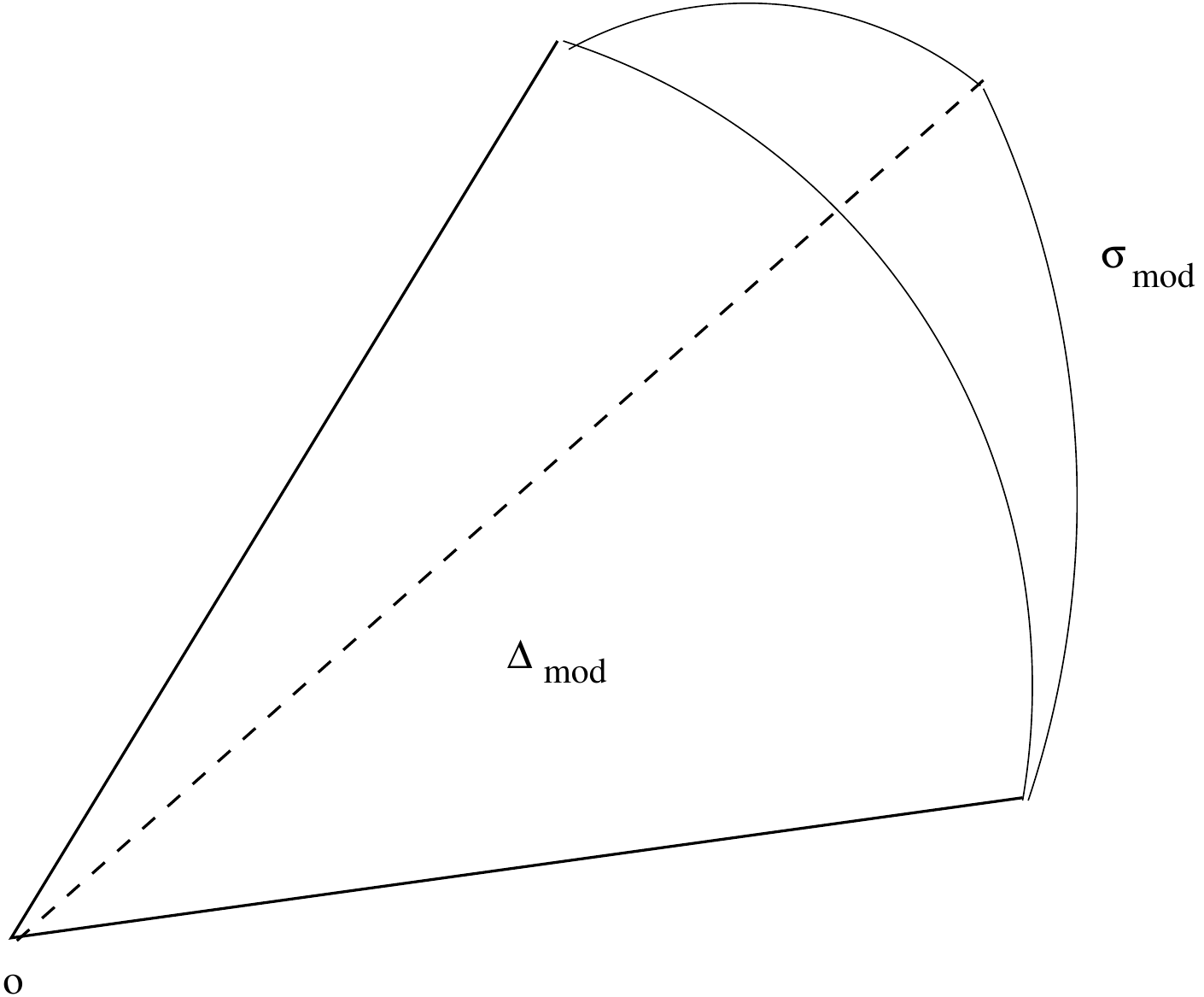}
\caption{Spherical and euclidean Weyl chambers.}
\label{figure3.fig}
\end{figure}

We pick  a maximal flat $F\subset X$ through the base point $o$ and regard it as the {\em model} flat $F_{mod}$.  
We will assume that $W$ fixes $o$ and  denote by $\Delta=\Delta_+=\Delta_{mod}\subset F_{mod}$ 
a certain fundamental domain of $W$, 
the {\em  model euclidean Weyl chamber} 
(see Figure \ref{figure3.fig}.) 
It is a complete  cone 
over a spherical simplex $\sigma_{mod}\subset \geo F_{mod}$, the  {\em model spherical Weyl chamber}. 
The tip of the cone $\Delta$ is the {origin} $o\in F_{mod}$. 
The cone $A_+\subset A$ is the subsemigroup of transvections preserving  the flat $F_{mod}$ and 
mapping $\De$ into itself, 
i.e.\ acting on $F_{mod}$ via translations 
$$
x\mapsto x+v, \quad v\in \Delta. 
$$
A {\em model Weyl sector} $V_{\taumod}$ is a face of $\Delta_{mod}$ 
which is the complete cone over a face $\taumod$ of $\simod$.
Euclidean Weyl chambers and Weyl sectors in $X$ are isometric copies of the model Weyl cone and the model Weyl sectors under $G$-congruences. 

We will frequently identify $\simod$ with a spherical simplex in the unit sphere in $F_{mod}$  centered at $o$, the intersection of this unit sphere with $\Delta$. 

\medskip 
{\em The opposition involution} $\iota: \Delta\to \Delta$ (also known as {\em Chevalley} or {\em duality} involution) is defined as the composition 
$$
\iota= w_0\circ (-id),
$$ 
where $w_0$ is the longest element of $W\acts F_{mod}$, the one sending the positive chamber in the model flat to the opposite one, and $-id$ is the antipodal map of $F_{mod}$ fixing  $o$. 

For each pointed maximal flat $(F,x)$ in $X$ 
there are finitely many euclidean Weyl chambers $V\subset F$ with tip $x$,
and they tessellate $F$. 

\begin{thm}
The following are equivalent:

1. The symmetric space $X$ is irreducible.

2. The action  $W\acts \R^r$ is irreducible.

3.  $G$ is a simple Lie group. 
\end{thm}

\medskip 
In the irreducible case, the Weyl groups $W$ are classified into $A$, $B=C$, $D$ (classical types) and $G_2, F_4, E_6, E_7, E_8$ (exceptional types). For instance, $SL_n$ has type $A_{n-1}$, $W\cong S_n$, the permutation group on $n$ symbols. 
The group $Sp_n$ has type $C_n$ and its Weyl group is isomorphic to the semidirect product 
$\Z_2^n \rtimes S_n$ where  $S_n$ acts on $\Z_2^n$ by permuting its basis elements.

\medskip 
{\em Walls} in $F_{mod}$ are the fixed hyperplanes of reflections in $W_{aff}$. Walls in $X$ are the images of walls in $F_{mod}$ under elements of $G$. 

\medskip 
{\em Regular} (geodesic) segments in $X$ are the segments  not contained in any wall. {\em Singular} segments are the segments contained in walls. 
Equivalently: A geodesic segment $xy$ is regular iff it is contained in a unique maximal flat. 

\medskip 
Each oriented segment $xy$ in $X$ defines a vector $v$ in $\Delta$, 
$$
v=d_\Delta(x,y),
$$
the {\em $\Delta$-valued distance from $x$ to $y$}.\footnote{The map $\mu: xy\mapsto v$ is also known as the {\em Cartan projection}, while the map $g\mapsto d_\Delta(x, gx)$ is sometimes called the {\em Lyapunov projection}.}
Namely, since $G$ acts transitively on pointed maximal flats in $X$, we can map $F$ to the model flat 
$F_{mod}$ and $x$ to the point $o\in F_{mod}$ via some $g\in G$. 
Now, project the oriented segment $g(xy)$ to the vector $v$ in $\Delta$ using the action of $W$. 

The vector $v$ 
is the {\em complete 
$G$-congruence invariant} of the pair $(x,y)$: Given two pairs $(x,y), (x',y')$, there exists $g\in G$ sending $(x,y)$ to $(x',y')$ iff 
$d_\Delta(x,y)= d_\Delta(x',y')$.

In the case of rank 1 spaces, $\Delta\cong[0,\infty)$ and $d_\Delta$ is the usual distance function. 

We refer to \cite{KLM} for the description of a complete set of {\em generalized triangle inequalities} for the chamber-valued distance function. The simplest of these inequalities has the form:
$$
d_\Delta(x,z)\le_{\Delta^*} d_\Delta(x,y) + d_\Delta(y,z),
$$
where $\Delta^*$ is the cone dual to $\Delta$, also known as the {\em root cone}:
$$
\Delta^*= \{u:  \< u, x\> \ge 0 \;\forall x\in \Delta \}. 
$$
 We also refer the reader to  \cite{Parreau} for  discussion of ``nonpositive curvature'' properties of  $d_\Delta$.

\begin{remark}\label{rem:ineq}
1. Here, given a convex cone $C$ with tip $0$ in a vector space $V$, we define the partial order $\le_C$ on $V$ by:
$$
u\le_C v \iff v-u\in C. 
$$

2. In general, $d_\Delta$ is not symmetric, but it satisfies the identity   
$$
d_\Delta(y,x) = \iota d_\Delta(x,y). 
$$
\end{remark}

\begin{rem} The theory of regular/singular segments has a relative analogue, relative to a face $\taumod$ of $\sigma_{mod}$; 
we will not cover the relative version in this paper. However, the relativization is important for the notion of $\tau_{mod}$-Morse maps and group actions, which correspond to $P$-Anosov subgroups 
in the sense of \cite{Labourie,GW} for parabolic subgroups $P<G$.
The discrete subgroups theory described in this survey is the one of $B$-Anosov subgroups, where $B < G$ is a minimal parabolic subgroup.
We refer the reader to \cite{morse} for the definition of $\taumod$-regularity. 
\end{rem}

\begin{example}
Consider the case of the symmetric space associated with the group $G=PGL(n,\R)$, 
i.e. $X$ consists of positive definite 
$n\times n$ matrices with unit determinant. Assume that $o\in X$ corresponds to the identity matrix. 
Then, up to scaling,
$$
d_\Delta(o,y)= \frac{1}{2} \bigl(\log(\la_1), \log(\la_2), \ldots, \log(\la_n)\bigr), $$
where $\la_1\ge \la_2\ge ....\ge \la_n$ are the eigenvalues of the matrix $y$ counted with multiplicity. 
The segment $oy$ is regular if and only if $\la_i> \la_{i+1}$ for all $i=1,...,n-1$. 
\end{example}

\subsection{Finsler geometry}\label{sec:finsler geometry}

Each symmetric space comes with a nonpositively curved Riemannian metric and the corresponding Riemannian distance function. Nevertheless, it turns out that many asymptotic aspects of $X$ (and of its quotients, locally symmetric spaces) are better captured by 
suitable $G$-invariant {\em polyhedral Finsler metrics} on $X$.  

Pick a regular vector $\bar\theta\in \sigma_{mod}$ (where we regard $\simod$ as a simplex in the unit sphere in $F_{mod}$), and define the  linear functional $\varphi$ on $F_{mod}$ dual to the vector $\bar\theta$. 
For simplicity, we assume $\bar\theta$ to be $\iota$-invariant. 
(See \cite{bordif} for the general treatment.)

\begin{rem}
There are several natural choices of the vectors $\bar\theta$ and, thus, of the dual 
linear functionals $\varphi$ and of the Finsler metrics defined below. For instance, 
one can take $\varphi$ to be the sum of all positive roots $\al\in R$ (positive with respect to the chamber 
$\Delta$). This linear functional will be {\em regular}, i.e.\ given by the inner product with a regular vector $\bar\theta$ 
in $\Delta$, and moreover $\iota$-invariant. While the metric $d_{\bar\theta}$ depends on the choice of $\bar\theta$, the compactification $\ol{X}^{Fins}$ of $X$ is independent of $\bar\theta$, see \cite{bordif}.  
For concreteness, the reader can assume that $\varphi$ is the sum of positive roots. 
\end{rem}

Given $\varphi$ (equivalently, $\bar\theta$), we define in \cite{bordif} 
a Finsler distance function $d_{\bar\theta}$ on $X$ as follows. 
First, we define a polyhedral Finsler norm on the vector space 
$F_{mod}=\R^r$ by 
$$
||v||_{\bar\theta}:= \varphi(d_\Delta(0, v)).  
$$
The unit ball $B_{mod}$ for this norm is the intersection of half-spaces
$$
\{x: (w^* \varphi)(x)\le 1\}, \quad w\in W. 
$$
Since this norm is $W$-invariant, it extends to a $G$-invariant Finsler metric on $X$ by defining the norm $||v||$ for a vector $v\in T_x X$ using the formula 
$$
||v||= ||dg(v)||_{\bar\theta}
$$
where $g: x\mapsto o\in X$, $dg(v)\in T_oF_{mod}$, $g\in G$. This norm on tangent spaces is 
a Finsler metric on the entire symmetric space $X$, 
and one has the Finsler distance function 
$$
d_{\bar\theta}(x,y):= \inf \int_{0}^{1} ||c'(t)|| dt,
$$
where the infimum is taken over all smooth paths $c: [0,1]\to X$, $c(0)=x, c(1)=y$. 
This distance function is also given by the explicit formula
$$
d_{\bar\theta}(x,y)= \varphi(d_\Delta(x,y)), \quad x, y\in X 
$$
which is the definition that we are going to use.
Due to our assumption that $\bar\theta$ is $\iota$-invariant,
the distance $d_{\bar\theta}$ is symmetric
and hence a metric in the usual sense.
We will refer to any such distance $d_{\bar\theta}$ as a {\em regular polyhedral Finsler metric} on $X$.

\begin{figure}[tbh]
\includegraphics[width=90mm]{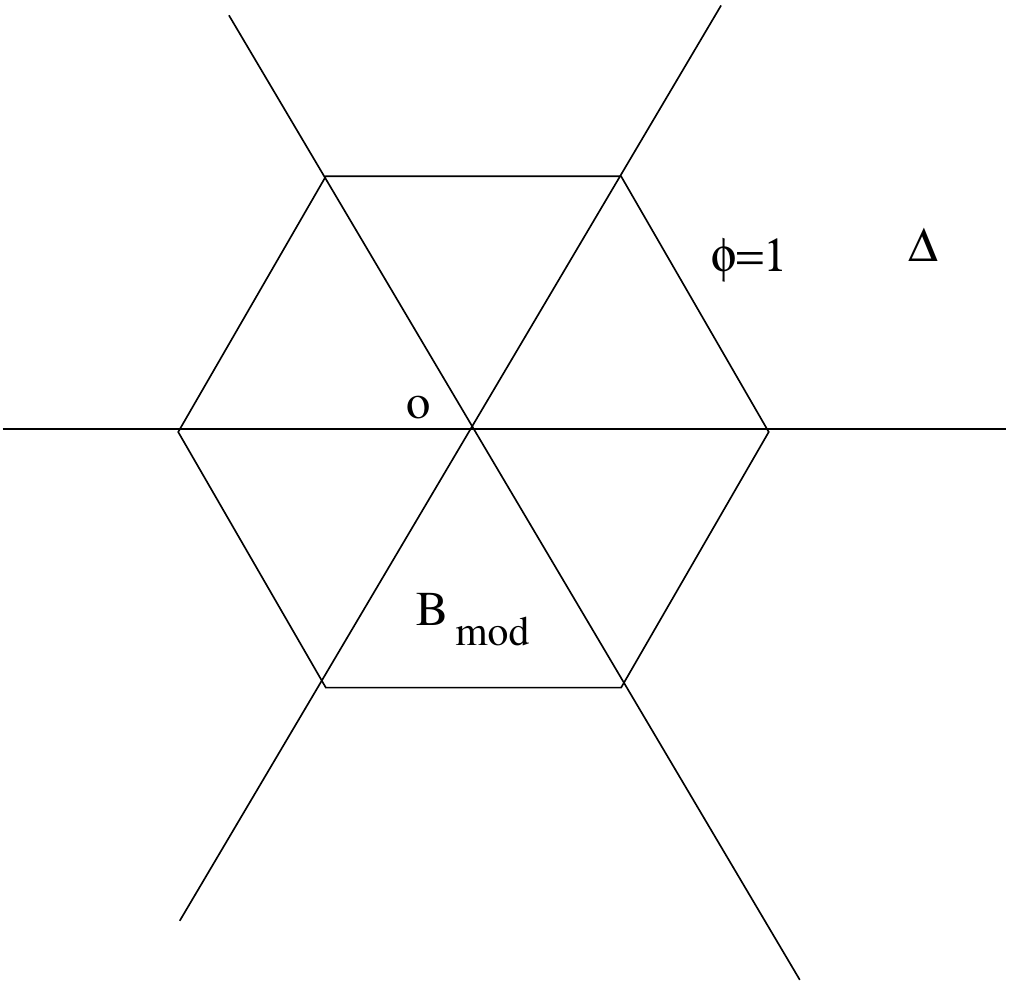}
\caption{Polyhedral Finsler norm.}
\label{figure4.fig}
\end{figure}

\medskip
{Regular polyhedral} Finsler metrics on $X$ are used in \cite{bordif} to construct a {\em Finsler compactification} $\ol{X}^{Fins}$ by adding to $X$ {\em Finsler horofunctions} in the manner similar to compactifying $X$ by adding to it Riemannian Busemann functions.  We will discuss this in more detail in \S \ref{Finsler-com}.

\subsection{Boundary theory} \label{sec:Boundary theory}

As in the rank 1 case, the {\em visual boundary} $\geo X$ of a symmetric space $X$ is defined as the set of asymptotic equivalence classes of geodesic rays in $X$: Two rays are asymptotic iff they are within finite distance from each other. There are two useful $G$-invariant topologies on $\geo X$: The first one is the visual topology, the topology of  a (the) unit tangent sphere in 
$X$.  This identification is achieved by choosing a reference point $x\in X$ and considering the set $Ray_x$ of geodesic rays $\rho$ emanating from $x$: Each geodesic ray in $X$ is equivalent to one and only one ray of such form. 
The set $Ray_x$ is identified with the unit tangent sphere $U_xX\subset T_xX$ by sending each ray $\rho$ to its velocity vector at $x$.

However, $\geo X$ also carries the structure of a (spherical) simplicial complex, defined via ideal boundary simplices of maximal flats in $X$. For each maximal flat $F$, the visual boundary $\geo F$  is identified with the unit sphere in $F$, and hence the $W$-action defines a Coxeter simplicial complex on $\geo F$.  

\begin{ff}
For any two maximal flats $F, F'$ the intersection $\geo F\cap \geo F'$ is a (convex) subcomplex of both
$\geo F$ and $\geo F'$.
\end{ff}

This proves that the tilings of the visual boundaries of the maximal flats are compatible.
The topology of this simplicial complex is called {\em Tits topology}.
It is induced by the {\em Tits metric}, which restricts to the angular metric on the visual boundary spheres of 
maximal flats. 
The simplicial complex is a {\em Tits building}, the {\em Tits boundary} of $X$, denoted $\tits X$. 
Its dimension equals $\rank(X)-1$. 
The identity map
$$
\tits X\to \geo X
$$
is a continuous (bijection), but never a homeomorphism,
i.e. the Tits topology is strictly finer than the visual topology.

{\em Apartments} in $\tits X$ are visual boundaries of maximal flats. 
Facets (i.e.\ top-dimensional simplices) of the apartments 
are called {\em chambers}.

Given a point $x\in X$ and a chamber $\si$ in $\tits X$, 
we let $V(x,\si)$ denote the {\em euclidean Weyl chamber} in $X$, which is the union of geodesic rays $x\xi$, $\xi\in \si$. Similarly, for a 
face $\tau$ of the simplicial complex $\tits X$,  we let $V(x,\tau)$ denote the {\em Weyl sector} equal to the union of  rays $x\xi$, $\xi\in \tau$. A point $\xi\in \geo X$ is {\em regular} if it belongs to the interior of a chamber $\si\subset \tits X$; equivalently, for some (every) $x\in X$ the geodesic ray $x\xi$ is regular.

\begin{ff} 
Any two ideal points (equivalently, chambers) belong to a common apartment. 
\end{ff}

Every $G$-orbit in $\geo X$ intersects every chamber exactly once,
and we have the {\em type map} 
$$\theta: \geo X\to \geo X/G\cong\sigma_{mod} .$$
For a maximal flat $F\subset X$,
the $G$-orbits in $\geo X$ intersect $\geo F$ in Weyl orbits,
and the restriction $\theta|_{\geo F}:\geo F\to\simod$ divides out the action of the Weyl group (of $F$ resp.\ $\geo F$).

\medskip 
{\bf Example:} (a) Rank 1 case: $\tits X$ is a discrete space.

(b) $SL(n, \R)$ case: $\tits X$ is 
the incidence complex of $\R P^{n-1}$.  Chambers are complete flags:
$$
V_\bullet=(V_1\subset \ldots\subset V_{n-1}\subset \R^n),
$$
where $\dim(V_i)=i$;  other faces are partial flags. 
The incidence relation: a partial flag $V_\bullet'$ is a face of a full flag $V_\bullet$ 
iff the full flag is a refinement of the partial flag.
For instance, if $n=3$, then full flags are pairs
$$
V_\bullet=(V_1\subset V_2),
$$
and partial flags are lines $V_1'$ or planes $V_2'$; they yield the vertices of the incidence graph. Then $V_1'$ is a vertex of 
$V_\bullet$ iff $V_1'=V_1$; $V_2'$ is a vertex of $V_\bullet$ iff $V_2'=V_2$. Thus, two vertices $V_1, V_2$ are connected by an edge iff $V_1\subset V_2$ (the line is contained in the plane). 

\begin{remark} $\rank(X)\geq2$ iff $\tits X$ is connected. 
\end{remark}

The {\em Furstenberg boundary}  $\DF X$ of $X$ is the space of {\em chambers} in $\tits X$.
The $G$-action on $\DF X$ is transitive
and the stabilizers in $G$ of the chambers are the {\em minimal parabolic subgroups} $B<G$.
Hence 
$$\DF X\cong G/B.$$
The topology on $\DF X$ induced by the visual topology coincides with its manifold topology 
as a homogeneous space. 
From the smooth viewpoint, 
$\DF X$ is a compact smooth homogeneous $G$-space,
and from the algebraic geometry viewpoint a homogeneous $G$-space with an underlying projective variety. 

For instance,
in the case $G=SL(n)$, 
the Furstenberg boundary is the full flag manifold,
and a minimal parabolic subgroup $B$ is given by the upper-triangular matrices, which is the stabilizer of the  full flag\footnote{Here 
and in what follows, $\< S\>$ denotes the linear span of a subset $S$ of a vector space.}
$$
\<e_1\>\subset \<e_1, e_2\>\subset \ldots \<e_1,\ldots, e_{n-1}\>. 
$$   
More generally,
for a face $\taumod\subseteq\simod$, we define the 
{\em generalized partial flag manifold} $\Flagt$ as the space of simplices $\tau\subset\tits X$ of type 
$\theta(\tau)=\taumod$.
The $G$-action on $\Flagt$ is again transitive.
The stabilizers of the simplices $\tau$ are the {\em parabolic subgroups} $P_{\tau}<G$ of type $\taumod$.
They form a conjugacy class and, denoting by $P_{\taumod}$ a representative,
we can write
$$\Flagt\cong G/P_{\taumod}.$$
Note that $\Flags=\DF X$.

For a simplex $\tau\in\Flagt$, we define its {\em star} $$\st(\tau)\subset\DF X$$ 
as the set of chambers of the Tits building $\tits X$ containing $\tau$:
\begin{equation}\label{eq:star}
\st(\tau):= \{\si\in \DF X: \tau\subset \si\}. 
\end{equation}

\begin{definition}
Ideal boundary points $\xi_\pm\in\geo X$ are {\em antipodal} if they are connected by a geodesic in $X$. 
Two chambers $\si_\pm$ are antipodal if they contain antipodal regular points. Equivalently, they are swapped by a Cartan involution in $X$. 
\end{definition}

{\bf Notation.} $\si^{opp}\subset\DF X$ denotes the set of chambers antipodal to $\si$. 

\begin{remark} 1. $\si^{opp}$  is an open subset of $\DF X$, called {\em open (maximal) Schubert cell of $\si$}.  

2. Antipodal implies distinct but not vice versa!

3. The complement $\DF X - \si^{opp}$ is a union of proper Schubert cycles in the projective variety $\DF X\cong G/B$, 
and hence a proper algebraic subvariety.
\end{remark}

\begin{example}
1. In {the} $SL(n)$ case, two full flags $V_\bullet, W_\bullet$ 
are antipodal iff they are {\em transversal}: $V_i$ is transversal to $W_{n-i}$ for each $i$. 

2. In the rank 1 case, antipodal is equivalent to distinct. 
The Tits boundary of a Gromov-hyperbolic space is a zero-dimensional building. 
\end{example}

\subsection{Quantified regularity}

Fix an $\iota$-invariant nonempty compact convex 
subset $\Theta\subset \simod^o$, where $\si_{mod}^o=\inte(\simod)$ is the interior of $\simod$. 
Define $V_\Theta\subset \Delta$, the $\Theta$-cone, as the cone with tip at the origin $o$ over the subset $\Theta$,
$$
V_\Theta= \R_{\ge 0} \cdot \Theta. 
$$
We define {\em $\Theta$-regular segments} in $X$ as segments whose $\Delta$-length  is in $V_\Theta$. 

More generally, given $x\in X$ and a euclidean Weyl chamber $V(x,\si)$, we define the {\em $\Theta$-cone}
$$
V_\Theta(x,\si)= \{y\in V(x,\si): d_\Delta(x, y)\in V_\Theta\}. 
$$

\begin{figure}[tbh]
\includegraphics[width=70mm]{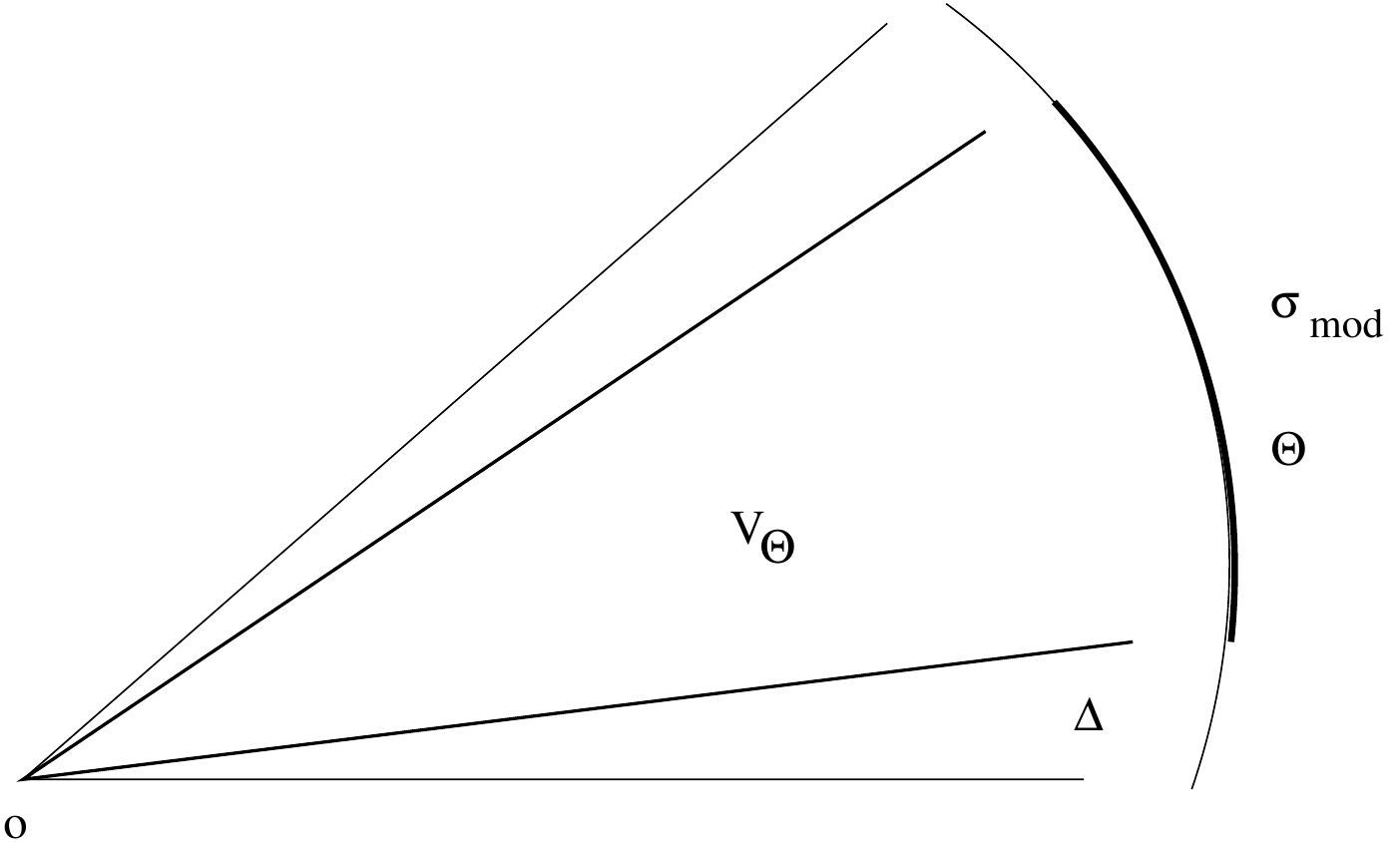}
\caption{Cone $V_\Theta=V_\Theta(0,\simod)$.}
\label{figure5.fig}
\end{figure}

\begin{remark}
Due to the $\iota$-invariance of $\Theta$, 
the notion of $\Theta$-regularity is independent of the orientation of the segments. 
\end{remark}

\medskip
For a negatively curved symmetric space $X$, a sequence $x_i\in X$ is divergent if and only if the sequence of distances $d(o,x_i)$ from a basepoint $o\in X$  diverges. Things become more complicated in higher rank symmetric spaces, since the ``right'' notion of a distance in $X$ is not a number but a vector in $\Delta$. This opens several possibilities for diverging to infinity and leads to (several)  notions of (asymptotic) regularity for sequences. In this survey we restrict to the simplest ones.

\begin{definition} 
[Regular sequence in $X$]
\label{def:regseq}
A sequence $(x_i)$ in $X$ is 

1. {\em regular} if the sequence of vectors $v_i=d_\Delta(o, x_i)$ diverges away from  the boundary of $\Delta$. 

2. {\em $\Theta$-regular} if it diverges to infinity and the sequence $(v_i)$ accumulates at $\Theta$. 
 
3. {\em uniformly regular} if it is $\Theta$-regular for some $\Theta$.  Equivalently, 
the accumulation set of the sequence $(v_i)$  is contained in $\simod^o$. Equivalently, there exists $\Theta'$ (a compact convex subset of $\simod^\circ$) 
such that for all but finitely many values of $i$ the vector $v_i$ belongs to $V_{\Theta'}$.  
\end{definition}

Analogously, we can define regularity for sequences of isometries $g_i$ of $X$: 
\begin{definition}
[Regular sequence in $G$]
\label{def:regseqgp}
A sequence $(g_i)$  in $G$ is {\em regular} (resp. {\em uniformly regular}, resp. {\em $\Theta$-regular}) 
if for some (equivalently, every) $x\in X$ the orbit sequence $x_i=g_i x$ has this property. 
\end{definition}

Thus, a divergent sequence in $X$ is uniformly regular iff all its subsequential limits in $\geo X$ are regular points. We will see later  how to characterize regular sequences $(g_i)$ in $G$ in terms of their action 
on the flag manifold $G/B$.

\medskip
\begin{remark} 1. Our notion of regularity for sequences is different from the notion 
introduced by Kaimanovich in \cite{Kaimanovich}, where a sequence in $X$ is called {\em regular}  
if it diverges at most sublinearly from a geodesic ray.

2. (Uniform) regularity of a sequence in $X$ is independent of the choice of  base point.

3. If $(x_i)$ and $(y_i)$ are sequences in $X$ within uniformly bounded distance from each other, 
$\sup_i d(x_i, y_i)<\infty$, then $(x_i)$ is (uniformly) regular iff $(y_i)$ is.  
\end{remark}

\begin{example}
\label{ex:regseq}
Suppose that $G=SL(n, \R)$ or $SL(n,\C)$. 
Then for each $g\in G$ we have its vector of singular values 
$$
a(g)=(a_1(g)\ge \ldots\ge a_n(g))
$$
where the $a_j$'s are the diagonal entries of the diagonal matrix $a$ in the Cartan resp.\ singular value decomposition,
cf.\ Example~\ref{ex:cartandeco}.
A sequence $(g_i)$ in $G$ is regular iff
$$
\lim_{i\to\infty} \frac{a_l(g_i)}{a_{l+1}(g_i)}=\infty
\quad \hbox{ for }l=1,...,n-1. 
$$
\end{example}

\begin{rem}
The singular values of a matrix depend on the choice of a euclidean/hermitian scalar product on $\R^n$ or $\C^n$ (this amounts to choosing a base point in the symmetric space of $G$), but the regularity of a sequence is independent of this scalar product.  
\end{rem}

In line with the notion of regular sequences in $X$ (which are maps $\N\to X$), one defines regular maps from other metric spaces into $X$. The most relevant for us is the following notion:

\begin{definition}
[Regular quasiisometric embedding]
\label{def:BTregular}
An $(L,A)$-quasiisometric embedding $f: Y\to X$ from a metric space $Y$ to a symmetric space $X$ 
is {\em $(\Theta,B)$-regular} if for all $y_1, y_2\in Y$ satisfying $d(y_1,y_2)\ge B$, 
the segment $f(y_1) f(y_2)$ is $\Theta$-regular in $X$. A map $f: Y\to X$ is a {\em uniformly regular quasiisometric embedding} if it is a $(\Theta,B)$-regular quasiisometric embedding for some $B$ and $\Theta$. 
\end{definition}

The most important cases when we will be using this definition are when $Y$ is a finitely generated group (equipped with a word metric) or a (possibly infinite) geodesic segment.  We will discuss regular quasiisometric embeddings and regular quasigeodesics in more details in the next section. 

\begin{figure}[tbh]
\includegraphics[width=90mm]{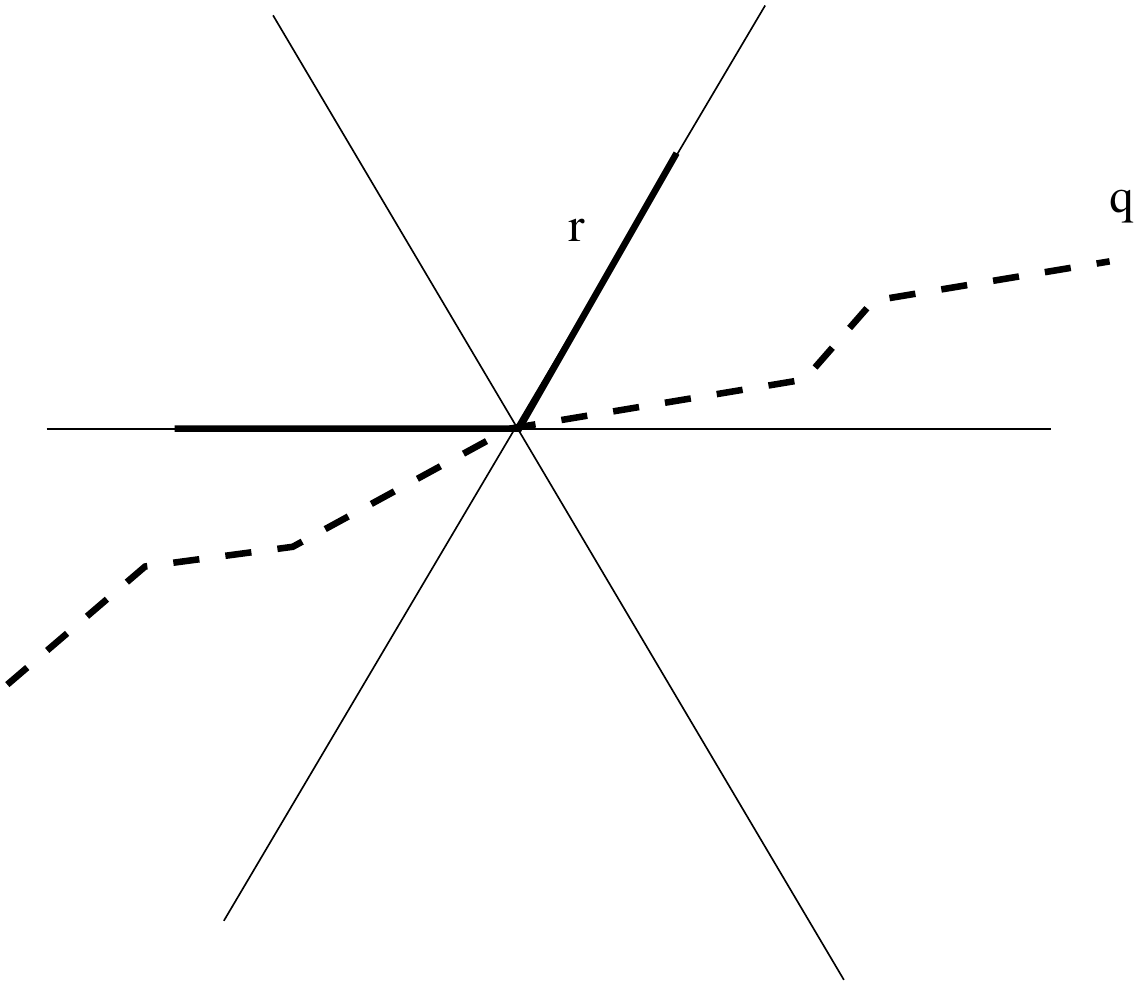}
\caption{The path $q$ in this figure is a Finsler geodesic, it is also a regular quasigeodesic in the model flat. The path $r$ is a Finsler geodesic but is not a regular quasigeodesic.}
\label{figure9.fig}
\end{figure}

\begin{example}[Regular and non-regular quasigeodesics]
Consider the case of quasiisometric embeddings into $\R^2=F_{mod}$, the model maximal flat of $SL(3)$.  We assume that the $x$-axis in $\R^2$ is a wall. Then, a piecewise linear function $f: \R\to \R$, yields a Finsler geodesic $q: x\mapsto (x, f(x))$, which is also a uniformly regular quasigeodesic in $F_{mod}$, provided that the slopes of linear segments in the graph of $f$ lie in the interval $[\eps, \sqrt{3} - \eps]$ for some  $\eps>0$. In contrast, the graph of the function $g(x)=|x|$ is not a regular quasigeodesic.
The reason is that for each $x>0$ the segment connecting the points $(-x, x), (x, x)$ in the graph of $g$ is horizontal and, hence, singular.  The graph $r$ of the function 
$$
h(x)= \begin{cases}
0 & \hbox{if~~}  x<0\\
\sqrt{3}x & \hbox{if~~}  x\ge 0
\end{cases}
$$
is a Finsler geodesic which is not a regular quasigeodesic. 
\end{example}

\medskip 
One of the geometric tools for studying regular quasiisometric embeddings are {\em diamonds} which we will define now. 
Diamonds 
can be regarded as the ``right'' generalization of geodesic segments 
when dealing with the metric $d_\Delta$  and with regular polyhedral Finsler metrics on $X$. 

\begin{definition}[Diamonds \cite{morse,anolec}]  
\label{def:diamo}
For a regular segment $xy$,
the {\em diamond} $\diamo_{x,y}\subset X$ 
is the intersection of the two 
euclidean Weyl chambers with tips $x$ and $y$ containing $xy$. 
\end{definition}

Diamonds are contained in maximal flats. 

\begin{figure}[tbh]
\includegraphics[width=90mm]{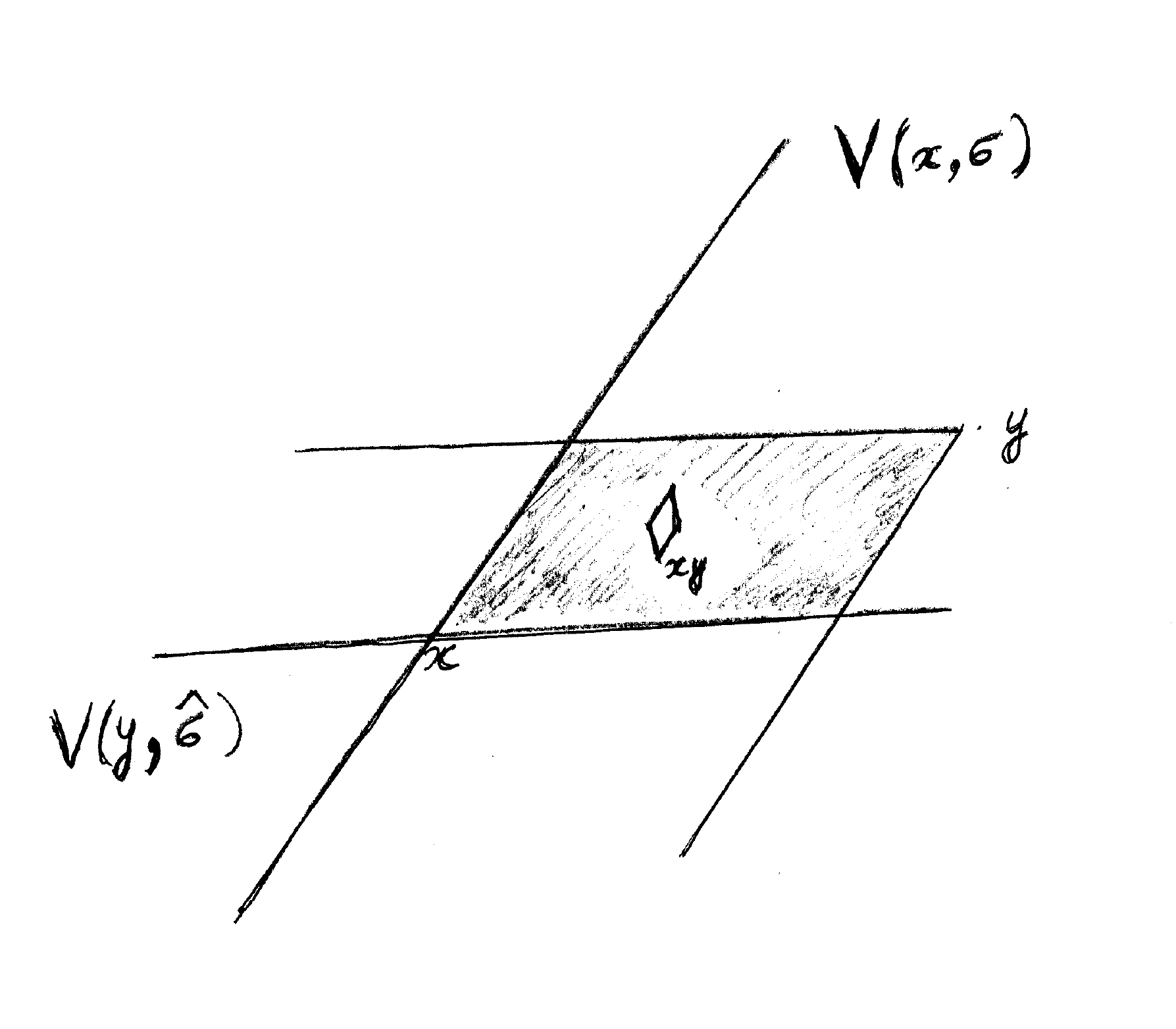}
\caption{A diamond in the model flat.}
\label{figure6.fig}
\end{figure}

\begin{example}
1. If $X$ is the product of rank 1 spaces
$$
X=X_1\times \ldots\times X_n
$$
then each diamond is the product of geodesic segments $s_i\subset X_i$. 

2. If $\rank(X)=2$, then each diamond is a parallelogram, its faces are contained in  walls.   
\end{example} 

\begin{figure}[tbh]
\includegraphics[width=90mm]{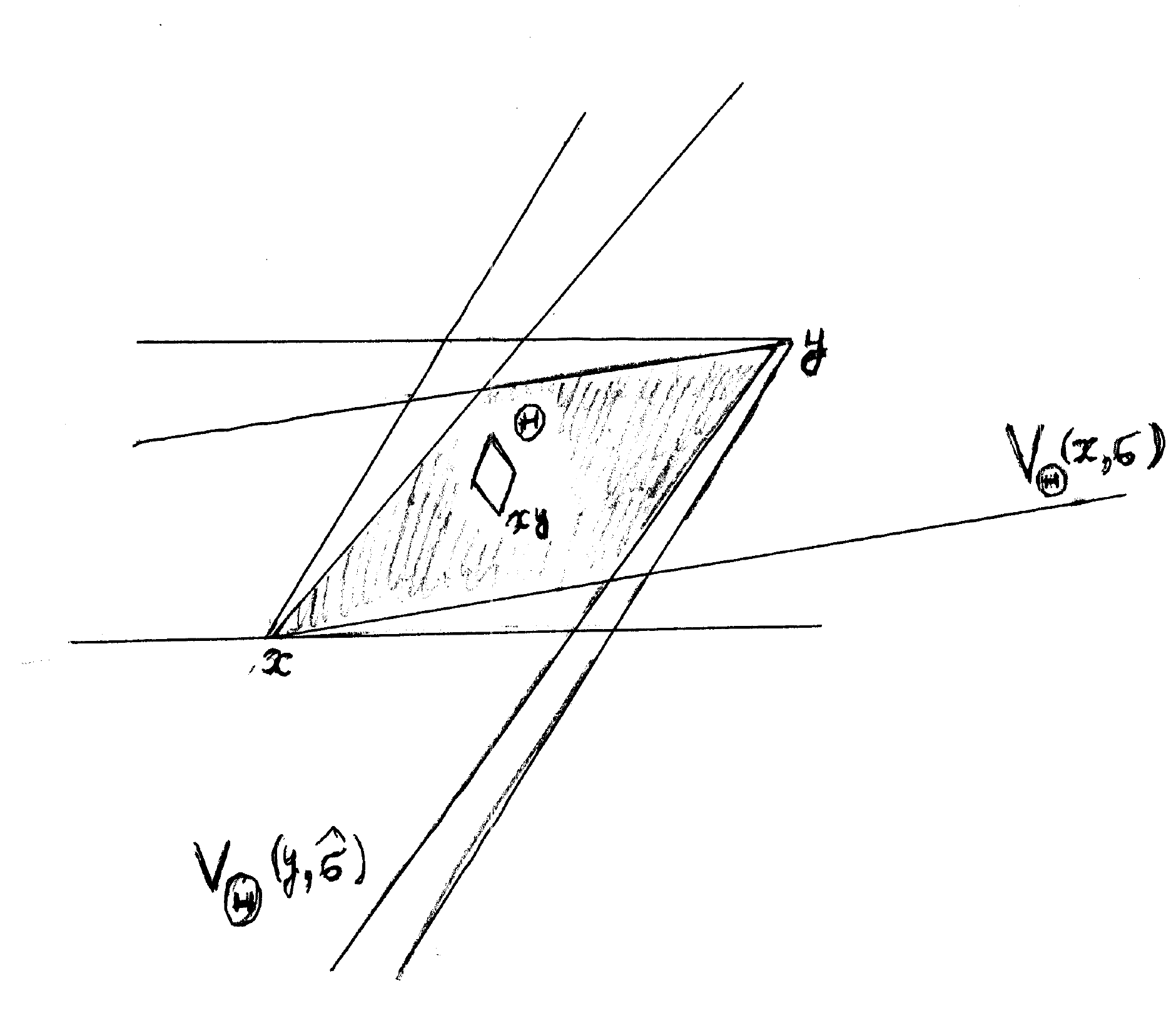}
\caption{$\Theta$-diamond in the model flat.}
\label{figure7.fig}
\end{figure}

Similarly, for $\Theta$-regular segments $xy$, 
one defines the {\em $\Theta$-diamond} $\diamo^{\Theta}_{x, y}\subset \diamo_{x,y}$ by
$$
\diamo^{\Theta}_{x,y}= V_\Theta(x,\si)\cap V_\Theta(y,\hat\si), 
$$  
where $\si, \hat\si$ are the (antipodal) chambers
such that $xy$ is contained in both $\Theta$-cones $V_\Theta(x,\si)$ and $V_\Theta(y,\hat\si)$.
As before, $xy\subset \diamo^{\Theta}_{x,y}$. See Figure \ref{figure7.fig}.

\begin{prop}
[Finsler description of diamonds \cite{bordif}]  
$\diamo_{x,y}$ is the union of all Finsler geodesics\footnote{With respect to a fixed regular polyhedral Finsler metric on $X$.} 
connecting $x$ to $y$. 
\end{prop}

It is quite clear that the diamond is filled out by Finsler geodesics connecting its tips.
The less part is to show that all Finsler geodesics connecting its tips are contained in the diamond.

Diamonds are enlargements of the Riemannian geodesic segments connecting their tips and,
in view of the proposition, 
may be regarded as their natural Finsler replacements,
reflecting the nonuniqueness of Finsler geodesics for polyhedral Finsler metrics.

\begin{figure}[tbh]
\includegraphics[width=90mm]{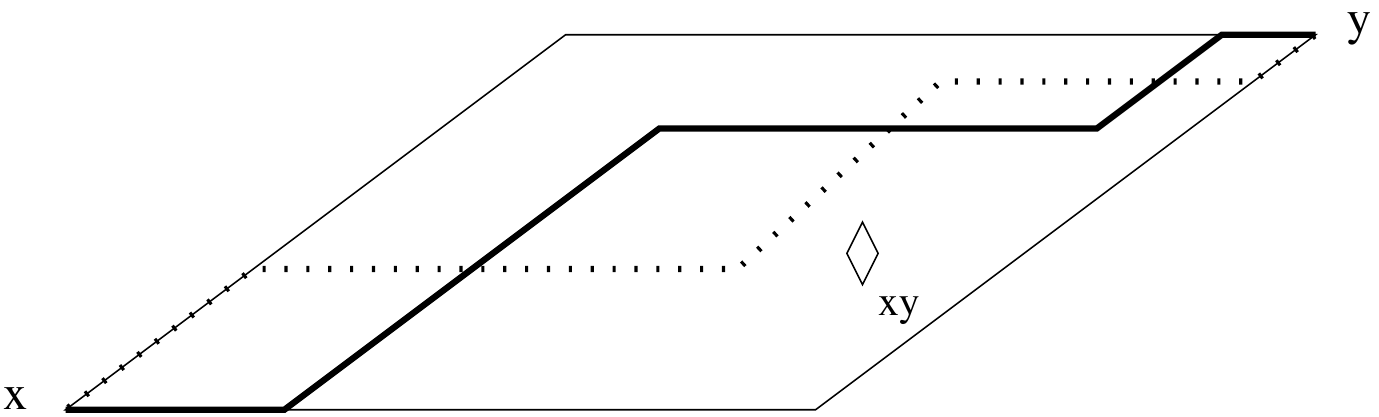}
\caption{Diamond as the union of Finsler geodesics.}
\label{figure8.fig}
\end{figure}

\medskip
\subsection{Morse quasigeodesics and the Higher Rank Morse Lemma}

\medskip
We will now discuss a higher rank version of the Morse Lemma for quasigeodesics in rank 1 spaces.

The Morse Lemma in rank 1 states that quasigeodesic segments are uniformly Hausdorff close to geodesic segments.
This is no longer true in higher rank, because it already fails in euclidean plane:

\begin{example}
[Failure of the naive version of the Morse Lemma] Take an $L$-Lipschitz function $f: \R\to \R$. Then $x\mapsto (x, f(x))$ is a quasigeodesic in $\R^2$, which, in general, is not close to any geodesic. For instance, take $f(x)=|x|$.  Further examples can be obtained by using suitable maps $r\mapsto (r,\theta(r))$ in polar coordinates.
\end{example}

We next define a class of quasigeodesics in $X$ 
which satisfy a higher rank version of the conclusion of the rank 1 Morse Lemma,
where geodesic segments are replaced with ``diamonds'',
see the previous section and Definition~\ref{def:diamo}.
That is,
we require the quasigeodesics and their subsegments to be uniformly close to $\Theta$-diamonds with tips at their endpoints:

\begin{definition}
[Morse quasigeodesics and maps \cite{morse}]  
\label{def:mrs}
Let $\Theta\subset \simod^o$ be nonempty $\iota$-in\-variant compact convex, 
and let $B, L, A, D,S>0$.

1. A map $q: I\to X$ from an interval $I$ 
is a {\em $(\Theta, B,L,A,D)$-Morse quasigeodesic} 
if it is an $(L,A)$-quasigeodesic, 
and if the image $q([s,t])$ is 
for any subinterval $[s,t]\subset I$ of length $t -s > B$ 
contained in the $D$-neighborhood of 
the $\Theta$-diamond $\diamo^{\Theta}_{x,y}$ with tips $x=q(s), y=q(t)$. 

2. A map $q: I\to X$ is a $(\Theta, B,L,A,D,S)$-{\em local Morse quasigeodesic} 
if its restrictions to subintervals of length $\leq S$ are $(\Theta, B,L,A,D)$-Morse quasigeodesics.

3. A map $f: Y\to X$ from a metric space $Y$ is {\em Morse} if 
it sends uniform quasigeodesics in $Y$ to uniform Morse quasigeodesics in $X$. 
 \end{definition}

Here, we call families of (Morse) quasigeodesics 
with fixed quasiisometry, respectively, Morse constants {\em uniform}.

Note that 
Morse quasigeodesics are {\em uniformly regular} in the sense of the previous section,
cf.\ Definition~\ref{def:BTregular}.\footnote{More precisely, 
$(\Theta, B,L,A,D)$-Morse quasigeodesics are $(\Theta',B')$-regular 
with $\Theta'\subset \simod^o$ and $B'>0$ depending only on 
$X$ and  the Morse data $(\Theta, B,L,A,D)$. }
Similarly, each Morse map is a uniformly regular quasiisometric embedding, provided that $Y$ is a {\em quasigeodesic metric space}, i.e.\ any two  points in $Y$ can be connected by a uniform quasigeodesic in $Y$. 

It is nontrivial that, vice versa, uniform regularity already forces quasigeodesics to be Morse:

\begin{theorem}[Higher Rank Morse Lemma \cite{mlem}] 
\label{thm:HRMorse}
Uniformly regular uniform quasigeodesics in $X$ are uniformly Morse. 
\end{theorem}

In other words, 
uniformly regular quasigeodesics in $X$ are Morse,
with Morse data depending on the quasiisometry constants and the uniform regularity data (and on $X$).

\begin{remark}\label{rem:addend}
1. The Morse Lemma holds as well for quasirays (diamonds are replaced with euclidean Weyl chambers)
and quasilines (diamonds are replaced with maximal flats\footnote{Or 
with unions of opposite Weyl cones in these maximal flats.}).

2. We also proved a version of our theorems with a weaker regularity assumption 
(relative to a face $\tau_{mod}$ of $\si_{mod}$). 
In this setting, diamonds are replaced by certain convex subsets of parallel sets, 
namely, intersections of opposite ``Weyl cones.''  
\end{remark}

For maps, one obtains accordingly:

\begin{corollary} 
Uniformly regular uniform quasiisometric embeddings $Y\to X$ from geodesic metric spaces $Y$ are uniformly Morse.
\end{corollary}

The closeness to diamonds in the Morse condition can be nicely reformulated 
in Finsler terms:

\begin{prop}
Uniform Morse quasigeodesics are uniformly Hausdorff close to Finsler geodesics.\footnote{Clearly, 
these Finsler geodesics are then uniformly regular.}
\end{prop}
The Morse Lemma for quasigeodesics then becomes:
\begin{corollary}
\label{cor:mlemfins}
Uniformly regular uniform quasigeodesics in $X$ are uniformly Hausdorff close to Finsler geodesics. 
\end{corollary}

There is a basic restriction on the coarse geometry of domains of uniformly regular quasiisometric embeddings
into symmetric spaces:

\begin{theorem}
[Hyperbolicity of regular subsets \cite{mlem}] 
If $Y$ is a geodesic metric space which admits a  uniformly regular quasiisometric embedding to $X$, 
then $Y$ is Gromov-hyperbolic. 
\end{theorem}

\begin{rem}
The Morse Lemma for quasirays implies that 
uniformly regular quasirays converge at infinity to Weyl chambers in a suitable sense. 
For uniformly regular quasiisometric embeddings from Gromov hyperbolic geodesic metric spaces, 
this leads to the existence of natural boundary maps.
We will make this precise at the end of next section,
see Theorem~\ref{thm:bdmp},
after introducing the notion of flag convergence of regular sequences in $X$ to chambers at infinity.
\end{rem}

It is a fundamental property of Morse quasigeodesics that they satisfy 
the following {\em local-to-global principle}: 

\begin{theorem}
[Local-to-global principle for Morse quasigeodesics \cite{morse}] 
\label{thm:L2G}
If a coarse Lipschitz path in $X$
is locally a uniform Morse quasigeodesic 
on a sufficiently large scale compared to the Morse data,
then it is globally a Morse quasigeodesic (for different Morse data).

More precisely:  
For Morse data $(\Theta, B,L,A,D)$ 
and another convex compact subset $\Theta'\subset\simod^o$ with $\Theta\subset\inte(\Theta')$ 
there exist constants $S,B',L', A', D'>0$ 
(depending also on $X$) 
such that every $(\Theta,B,L,A,D,S)$-local Morse quasigeodesic 
is a $(\Theta',B',L',A',D')$-Morse quasigeodesic.
\end{theorem}

This theorem parallels the local-to-global principle for quasigeodesics in Gromov-hyperbolic spaces, 
see e.g.\ \cite[Thm.\ 1.4 in ch.\ 3]{CDP}.  
It is derived from a basic local-to-global principle for straight paths
which we will explain later,
see Theorem~\ref{thm:straight-paths} below.

\subsection{Flag convergence} 
\label{sec:flagcv}

We introduce the following notion of convergence at infinity for regular sequences in $X$ to chambers in $\DF X$:
\begin{definition}[Flag convergence \cite{coco13, morse, anolec}]
\label{defn:flagconvergenceX}
A regular sequence $(x_n)$ in $X$ is said to 
{\em flag converge} to 
a chamber 
$\si\in \DF X$, 
if for some base point $o\in X$ and some sequence $(\si_n)$ in $\DF X$ with 
\begin{equation*}
\sup_nd\bigl(x_n, V(o,\si_n)\bigr) < \infty
\end{equation*}
it holds that $\si_n\to\si$ (in the manifold topology of $\DF X$).
\end{definition}
 
This convergence is independent of the choices of $o$ and $(\si_n)$, see \cite{coco13, anolec}. 

For {\em uniformly} regular sequences,
flag convergence can be described in terms of the visual compactification:
A uniformly regular sequence in $X$ flag converges to $\si\in \DF X$ 
iff its accumulation set in the visual compactification $\ol{X}$ is contained in $\si$. 

Flag convergence is induced by a natural topology on $X\sqcup \DF X$,
making it a {\em partial compactification} of $X$,
see \cite[\S 3.8]{mlem}. 
If $\rank(X)=1$, then $\DF X$ is the visual boundary of $X$ and the topology on 
$X\sqcup \DF X$ is the visual topology described in \S \ref{sec:rk1}, making $X\sqcup \DF X$ homeomorphic to a closed ball. The situation in higher rank is more complex, 
since then $\DF X$ is not even a subset of the visual boundary $\geo X$.\footnote{However,
$\DF X$ is a subset of the {\em Finsler} boundary, see \S~\ref{Finsler-com} below.}
The topology on $X\sqcup \DF X$ is obtained as follows.
Fix a basepoint $o\in X$
and define the {\em shadow} $Sh(B(y,R))$ of an {open} metric ball $B(y,R)\subset X$ in $X\sqcup \DF X$ as
$$
\{x\in X: ox \hbox{~~is regular and~~} \diamo_{ox}\cap B(y,R)\ne \emptyset\} \cup \{\si\in \DF X: V(o, \si)\cap B(y,R)\ne \emptyset\}. 
$$
Then a basis of the {\em shadow topology} on $X\sqcup \DF X$ at $\si\in \DF X$ consists of all sets $Sh(B(y,R))$ with $R>0$ and $y\in V(o, \si)$. We retain the metric topology on $X$. 

\begin{prop}
1. The shadow topology is independent of the basepoint $o\in X$. 

2. The shadow topology is 2nd countable and Hausdorff. 

3. The shadow topology restricts on $\DF X\cong G/B$ to the manifold topology. 

4. In rank 1, the shadow topology coincides with the visual topology.
\end{prop}

A regular sequence in $X$ flag converges to a chamber $\si\in \DF X$ 
iff it converges to $\si$ in the shadow topology, see \cite{mlem}.

\begin{figure}[tbh]
\includegraphics[width=90mm]{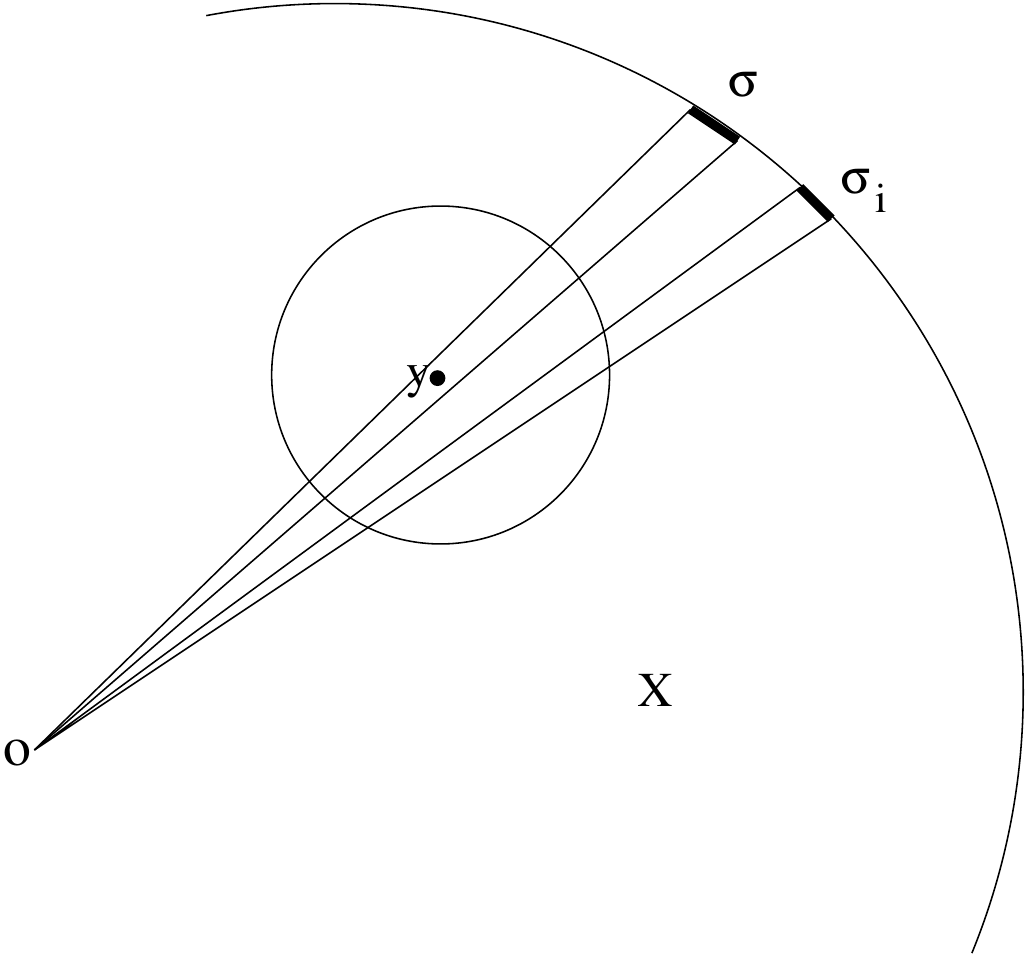}
\caption{Shadow topology and flag convergence: The sequence of chambers $\si_i$ flag converges to the chamber $\si$.} 
\label{figure18.fig}
\end{figure}

We extend the notion of flag convergence to sequences in $G$:

\begin{definition}[Flag convergence in $G$]
\label{defn:flagconvergenceG}
A regular sequence $(g_n)$ in $G$ {\em flag converges} to $\si\in \DF X$ 
if for some (equivalently, every) $x\in X$ the sequence $(g_nx)$ flag converges to $\si$.  
\end{definition}

Proposition \ref{prop:equiv-flag} in section \ref{sec:conprop} will provide equivalent conditions for flag convergence of sequences in $G$. 

\medskip
Now, with the flag topology at our disposal,
we can formulate the following {\em Addendum to the Higher Rank Morse Lemma}
regarding boundary maps for regular quasiisometric embeddings:

\begin{theorem} [Existence of boundary map \cite{mlem}] 
\label{thm:bdmp}
Each uniformly regular quasiisometric embedding $f: Y\to X$ 
from a $\delta$-hyperbolic geodesic metric space $Y$
continuously extends to a map $$Y\sqcup \geo Y\to X\sqcup \DF X.$$
The boundary extension $\geo f: \geo Y\to \DF X$ is antipodal,
i.e.\ maps distinct ideal points to opposite chambers.
\end{theorem}

\begin{remark} 
1. The appearance of $\DF X$ in the extension map comes from the fact that the restrictions $q=f\circ r$ of $f$ to geodesic rays $r=y\xi$ in $Y$ are Morse quasirays and hence are uniformly close to (subsets of) euclidean Weyl chambers $V(f(y),\sigma)\subset X$, see part 1 of Remark \ref{rem:addend}.  
The quasiray $q$ then accumulates at the chamber $\sigma=:\geo f(\xi)$. 

2. The antipodality of the boundary map follows from the Higher Rank Morse Lemma 
for quasilines,
see also part 1 of Remark \ref{rem:addend}.
\end{remark}

\subsection{Finsler compactifications}\label{Finsler-com}

We apply the horofunction compactification construction (see Appendix \ref{sec:horoboundary}) 
to the symmetric space $X$ and the Finsler distance 
function $d_{\bar\theta}$. The resulting compactification 
$$
\ol{X}^{\bar\theta}=\ol X^{Fins}=X\sqcup\geo^{Fins}X 
$$
is independent of $\bar\theta$ (as long as it is regular), in the sense that the identity map $X\to X$ extends to a homeomorphism of the compactifications
$$
\ol{X}^{\bar\theta_1} \to \ol{X}^{\bar\theta_2}
$$
for any two regular elements $\bar\theta_i\in \simod$.

\medskip 
In the case of the horofunction compactification of symmetric spaces (equipped with their standard Riemannian metrics), horofunctions were identified with asymptotic equivalence classes of geodesic rays in $X$. A similar identification can be done in the Finsler setting, but rays are replaced with {\em Weyl sectors} and the equivalence relation is a bit more complicated.  

Two Weyl sectors $V(x,\tau)$,  $V(x', \tau')$ are equivalent if and only if:

1. $\tau=\tau'$. 

2. For every $\eps>0$ there exists $y\in V(x,\tau), y'\in V(x',\tau')$ such that the Hausdorff distance between 
$V(y,\tau)$,  $V(y', \tau')$ is $<\eps$. 

\medskip
Note that if $\rank(X)=1$, Weyl sectors are geodesic rays in $X$ and two sectors are equivalent iff the rays are asymptotic. The connection of equivalence classes of sectors to Finsler horofunctions comes from the following theorem that allows an identification of element of $\geo^{Fins}X$ with equivalence classes $[V(x,\tau)]$ of Weyl sectors. In this theorem we let $d^{\bar\theta}_{x}$ be the function sending $y\in X$ to  $d_{\bar\theta}(x,y)$.  

\begin{thm}
[Weyl sector representation of points at infinity in Finsler compactification \cite{bordif}]  
1. Let $x_i\in V(p, \tau)$ be a sequence diverging away from the boundary faces of the sector $V(p,\tau)$. Then 
the sequence of Finsler distance functions $d^{\bar\theta}_{x_i} - d_{\bar\theta}(x_i, p)$ converges to a horofunction which will be denoted  
$b_{p,\tau}^{\bar\theta}$. The limit horofunction is independent of the sequence $(x_i)$. 

2. Every horofunction $b\in \geo^{Fins}X$ is equivalent (i.e. differ by a constant) to a horofunction 
$b_{p,\tau}^{\bar\theta}$. 

3. Two Finsler horofunctions $b_{p,\tau}^{\bar\theta}$, $b_{p',\tau'}^{\bar\theta}$ are equivalent if and only if the 
sectors $V(p, \tau)$, $V(p', \tau')$ are equivalent. 

4. The identification 
$$
\left[b_{p,\tau}^{\bar\theta}\right] \leftrightarrow [V(p, \tau)]$$
is $G$-equivariant, where $G$ acts on horofunctions by the precomposition. 
\end{thm}

This identification determines the following stratification of $\geo^{Fins}X$. The {\em small strata} are the sets
$$
S_{\tau}= \{[V(p, \tau)]: p\in X\},
$$
where $\tau$'s are simplices in $\tits X$. The {\em big strata} $S_{\taumod}$ are the unions
$$
S_{\taumod}= \bigcup_{\tau\in \theta^{-1}(\taumod)} S_\tau.  
$$
The group $G$ acts on each big stratum transitively. This $G$-invariant stratification extends to $\ol{X}^{Fins}$ by declaring the entire $X$ to be a single big stratum, $S_{\emptyset}$. The smallest big stratum is $S_{\simod}$, which is $G$-equivariantly homeomorphic to $\DF X$. This stratum is 
the unique closed $G$-orbit in $\ol X^{Fins}$; 
$$S_{\simod}\cong \DF X\cong G/B.$$
On the opposite extreme, the orbit $S_{\emptyset}=X$ is open and dense in $\ol{X}^{Fins}$.

The strata $S_{\taumod}$ for $\taumod\neq\emptyset$,
are {\em blow-ups} of the corresponding flag manifolds $\Flagt= G/P_{\taumod}$, 
where $P_{\taumod}$ are representatives of conjugacy classes of parabolic subgroups of $G$, 
parameterized by faces $\taumod$ of $\simod$. More precisely,
there are $G$-equivariant fibrations 
$$S_{\taumod}\lra\Flagt$$
with contractible fibers.
The fiber $S_{\tau}\subset S_{\taumod}$ over $\tau\in\Flagt$ 
can be interpreted geometrically 
as the space of {\em strong asymptote classes of Weyl sectors} $V(x,\tau)$ asymptotic to $\tau$,
cf.\ \cite[\S 3]{bordif}. 
In particular, it is a symmetric space of rank 
$$\dim\simod-\dim\taumod<\rank X=1+\dim\simod.$$
The topological boundary $\D S_{\tau}=\ol S_{\tau}-S_{\tau}$ of $S_\tau$ 
is the union of small strata,
namely of the $S_{\nu}$ for the simplices $\nu$ strictly ``refining" $\tau$
in the sense that $\nu\supsetneq\tau$.

\begin{thm}[\cite{bordif}] 
$\ol X^{Fins}$ is $K$-equivariantly homeomorphic to the closed unit ball in $X$ with respect to the  {\em dual Finsler metric} $d^*_{\bar\theta}$ on $X$. This compactification is $G$-equivariantly homeomorphic to the {\em maximal Satake compactification} $\ol{X}^S_{max}$ of the symmetric space $X$, see \cite{Borel-Ji} for the definition.
\end{thm}

\begin{cor}
$\ol X^{Fins}$ is a real-analytic manifold with corners on which $G$ acts real-analytically. 
\end{cor}

\begin{figure}[tbh]
\includegraphics[width=90mm]{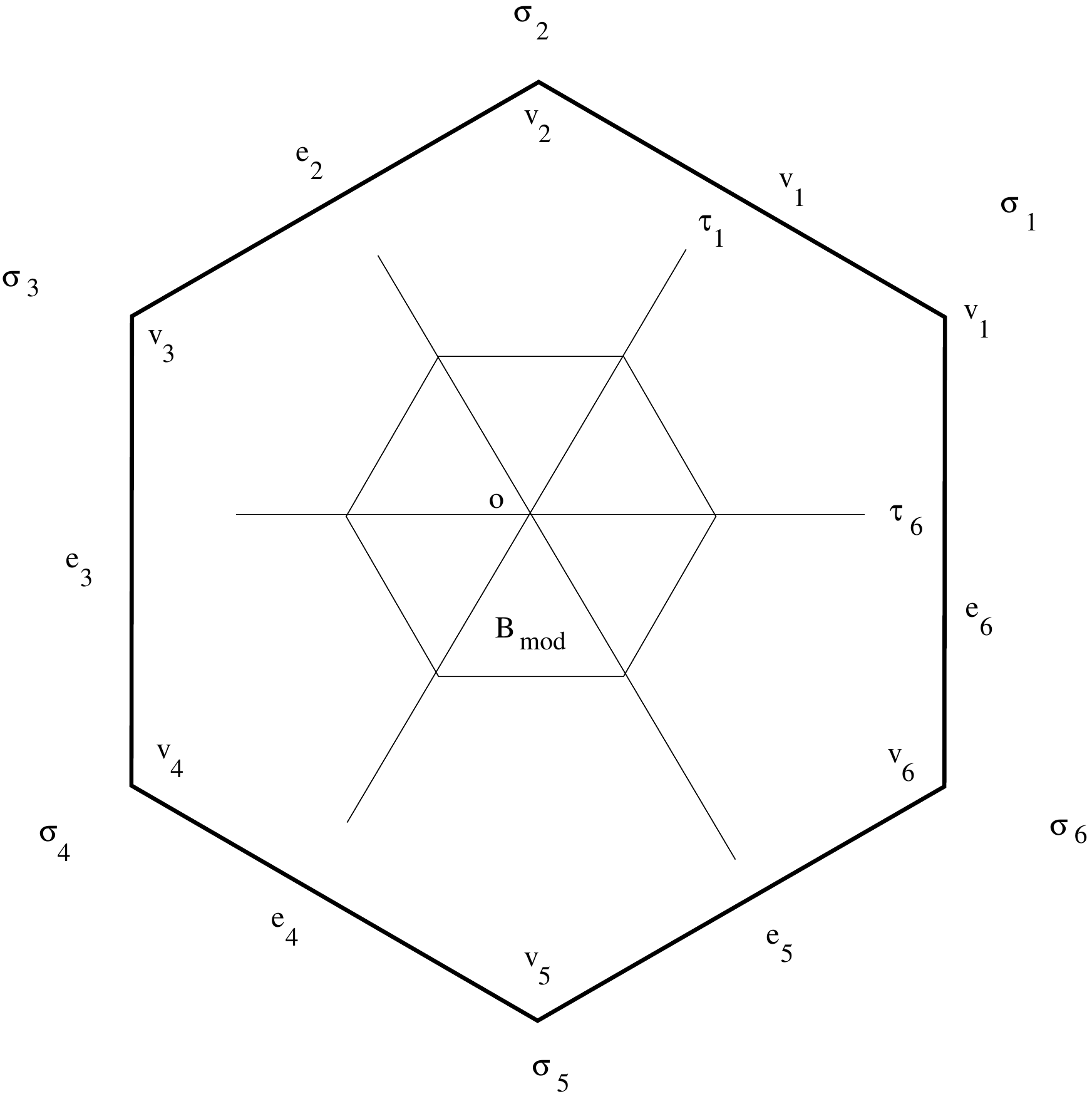}
\caption{Finsler compactification of the model flat.}
\label{figure10.fig}
\end{figure}

\begin{example} We now describe the (regular polyhedral) 
Finsler compactification of the model flat $F_{mod}$ for $SL(3,\R)/SO(3)$, in which case $\bar\theta$ is 
the midpoint of {the edge} $\simod$. 
Let $\si_1,\ldots, \si_6$ denote the spherical chambers of $F_{mod}$ listed in the cyclic order. Let $\zeta_i\in \si_i$ denote the midpoint of $\si_i$. Let $\tau_{i}$ denote the common vertex of $\si_{i}, \si_{i+1}$. 
Each chamber $\si_i$ determines a vertex $v_i$ of the Finsler compactification $\ol{F_{mod}}^{Fins}$ of $F_{mod}$; each $\tau_i$ determines an edge $e_i$ of $\ol{F_{mod}}^{Fins}$. In terms of Finsler horofunctions: Each vertex $v_i$ corresponds to the horofunction whose restriction to $F_{mod}$ is 
$$
b_{\zeta_i}(x)= - x\cdot \zeta_i. 
$$
Each edge $e_i$ corresponds to the 1-parameter (modulo additive constants) family of Finsler horofunctions
$$
\max \left( b_{\zeta_i} + s,  b_{\zeta_{i+1}} + t\right), \quad s, t\in \R. 
$$
Using the normalization $s+t=0$ ({we are} assuming that all Finsler horofunctions are normalized to vanish at the origin), we can write this family as
$$
b_{i,s}=\max \left( b_{\zeta_i} + s,  b_{\zeta_{i+1}} - s\right) - |s|, \quad s\in \R. 
$$
As $s\to+ \infty$, $b_{i,s}$ converges (uniformly on compacts in $F_{mod}$) to $b_{\zeta_i}$, while as $s\to-\infty$, $b_{i,s}$ converges to $b_{\zeta_{i+1}}$, representing the two vertices of the edge $e_i$. We, thus, obtain a description of the stratified space $\ol{F_{mod}}^{Fins}$ as a hexagon, dual to the unit ball $B_{mod}$ of the regular polyhedral Finsler norm on $F_{mod}$. 

Regarding the small strata of $\ol{X}^{Fins}$: They are points (corresponding to the spherical chambers, elements of $\DF X$) and open 2-dimensional disks, which have natural geometry of hyperbolic 2-planes, and $X$ itself. Note that there are two types of open 2-disks, corresponding to two types of vertices of the spherical building $\geo X$. Taking two opposite vertices $\tau, \hat\tau$ of $\geo X$ and the parallel set $P(\tau, \hat\tau)$ (the union of all geodesics asymptotic to $\tau, \hat\tau$) splits as $\H^2\times \R$. The Finsler compactification of this parallel set contains $\H^2\times \{\pm \infty\}$,  the open disk strata of $\ol{X}^{Fins}$ which have different type. See Figure \ref{figure10.fig}. 
\end{example}

We refer the reader to \cite{bordif}  for more details and to \cite{JS} for the description of compactifications of  finite dimensional vector spaces equipped with polyhedral norms. 

\begin{definition}
We say that a subset of $\geo^{Fins}X$ is {\em saturated} if it is a union of small strata. 
\end{definition}

It is worth noting that 
the stabilizers of  points in the Finsler compactification 
are {\em pairwise different} closed subgroups of $G$.
The stabilizers of the points at infinity in $S_{\tau}$ are contained in the parabolic subgroup $P_{\tau}$, where $P_\tau$ is the stabilizer in $G$ of the simplex $\tau$.

We conclude this section with the following theorem which provides a satisfying metric interpretation of the shadow topology: 

\begin{thm}[Prop. 5.41 in \cite{bordif}] 
The subspace topology on $X\sqcup \DF X$ induced from $\ol X^{Fins}$ is equivalent to the shadow topology on 
$X\sqcup \DF X$. 
\end{thm}

\subsection{The higher rank convergence property} 
\label{sec:conprop}

We consider the action of $G$ on the full flag manifold $G/B=\DF X$. The 
usual convergence property,
compare section~\ref{sec:rank1conv}, 
fails in this context: 
In higher rank, 
a divergent sequence $(g_k)$ 
never converges
to a constant map 
on the complement of a point in $\DF X$. 
However, as we noted earlier, in higher rank {\em distinct} should be replaced with {\em antipodal}. 

Given two chambers $\al, \om \in \DF X$ we define the {\em quasiprojective map}
$$
\al_\om: \om^{opp} \to \{\al\},
$$
left undefined on the set $\DF X - \om^{opp}$ consisting of chambers which are not antipodal to $\om$. The chamber $\al$ is called the {\em attractor} and $\om$ is called the {\em repeller}. We say that a sequence $(g_k)$ in $G$ {\em converges} to a quasiprojective map $\al_\om$ if $g_k$ converges to $\al$ uniformly on compacts in $\om^{opp}$. 

\begin{theorem}
[The higher rank convergence property \cite{coco13, morse,anolec}] 
\label{thm:regconv}
Each regular sequence $(g_k)$ in $G$ contains a subsequence $(g_{k_i})$ which converges to  the map $\al_\om$ for some $\al, \om\in \DF X$. Conversely, if a sequence $(g_k)$ has such a limit $\al_\om$, then it is regular. 
\end{theorem}

\begin{figure}[tbh]
\includegraphics[width=90mm]{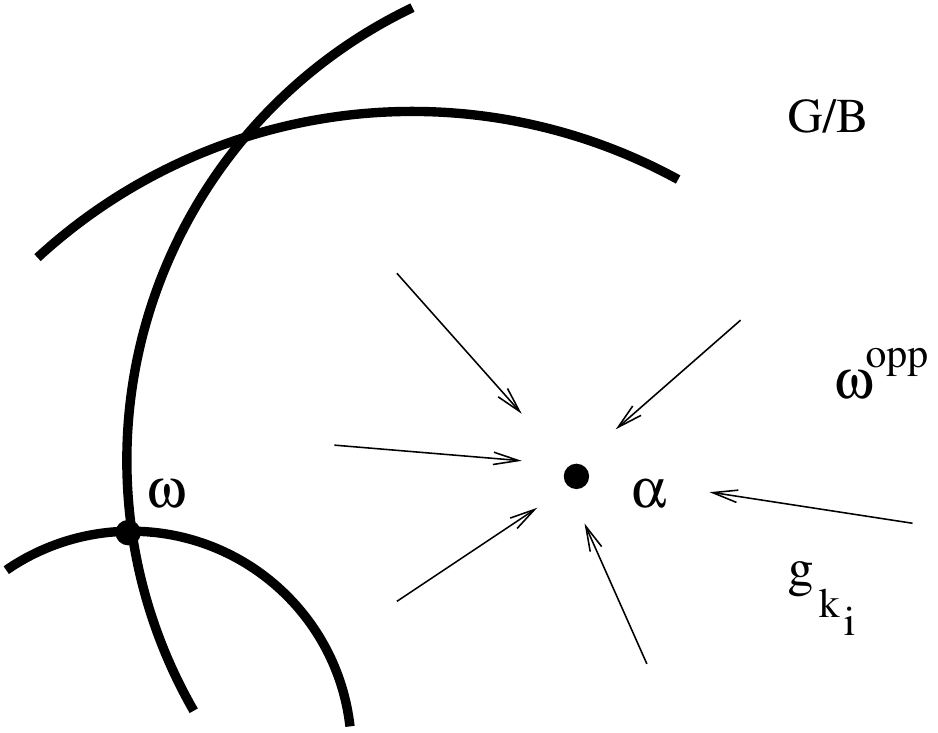}
\caption{The convergence property.}
\label{figure11.fig}
\end{figure}

\begin{remark} 
The complement $G/B -  \om^{opp}$ is the exceptional set for this convergence (where uniform limit fails locally). 
\end{remark}

This theorem gives a {\em dynamical characterization} of regular sequences in $G$:

\begin{cor}
\label{cor:regconveq}
A sequence $(g_k)$ in $G$ is regular iff every subsequence $(g_{k_i})$ contains a further subsequence 
$(g_{k_{i_j}})$ which converges to some quasiprojective map $\al_\om$. 
\end{cor}

As in the rank 1 case:
$$
g_k\to \al_\om \iff g_k^{-1}\to \om_\al. 
$$

\begin{rem}
1. More generally, one defines $\tau_{mod}$-regularity of a sequence relative to a face $\tau_{mod}\subseteq \sigma_{mod}$. Each $\tau_{mod}$ determines a (partial) flag manifold 
$$\Flagt=G/P_{\tau_{mod}}.$$
Then the convergence property (for arbitrary sequences $(g_k)$ in $G$) reads as:

Each sequence $(g_k)$ in $G$ contains a subsequence $(g_{k_i})$ which is either bounded in $G$ or is 
$\tau_{mod}$-regular for some face $\tau_{mod}$. 
The latter is equivalent to convergence (uniform on compacts) of $(g_{k_i})$ to a quasiprojective map 
$$
\al_{\om}: \om^{opp}\subset \Flagt  \to \al\in \Flagt.$$
Here $\om$ is a face of $\tits X$ of the type opposite to $\tau_{mod}$. 
We refer to \cite{morse,anolec} for details.  

2. An equivalent notion of convergence of sequences in $G$ had been introduced earlier by Benoist in \cite{Benoist},
see in particular part (5) of his Lemma 3.5.  
\end{rem}

\begin{example}
Consider the case $G=SL(n, \R)$ and a sequence of diagonal matrices $g_k= Diag(a_{1,k},\ldots, a_{n,k}) \in G$ 
with $a_{1,k}\geq\ldots\geq a_{n,k}>0$. 
Recall from Example~\ref{ex:regseq}
that regularity of the sequence $(g_k)$ amounts to the conditions 
$$
\lim_{k\to\infty} \frac{a_{i,k}}{a_{i+1,k}}=\infty \quad\hbox{ for } i=1,...,n-1. 
$$
The attractive flag for the sequence $(g_k)$ is 
$$
\al= \left(\<e_1\> \subset \< e_1, e_2\> \subset \ldots \subset \< e_1,...,e_{n-1}\>\right),
$$ 
and the repelling flag is 
$$
\om = \left(\<e_n\> \subset \< e_n, e_{n-1}\> \subset \ldots \subset \< e_n,...,e_{2}\>\right).
$$
\end{example}

{It is useful to reformulate 
the definition of flag convergence of a regular sequence $(g_n)$ to a chamber $\al\in \DF X$
in terms of the dynamics on the flag manifold $\DF X$:}

\begin{prop}
[Flag convergence criteria, \cite{coco15}] 
\label{prop:equiv-flag}

The following are equivalent for a regular sequence $(g_n)$ in $G$:

1.  The sequence $(g_n)$ flag converges to $\al\in \DF X$. 

2. Every subsequence in $(g_n)$ contains a further subsequence which converges to a quasiprojective map 
$\al_\om: \DF X\to \DF X$. 

3. There exists a bounded sequence $b_n\in G$ such that the sequence $(g_n b_n)$  converges to a quasiprojective map $\al_\om: \DF X\to \DF X$. 
\end{prop}

We equip the smooth compact manifold $G/B$ with an auxiliary Riemannian metric 
(not necessarily $K$-invariant). 
This allows us to define expansion properties for elements $g\in G$ at chambers $\si\in G/B$ in the same way it was done in the rank 1 situation, see section \ref{sec:expansion}. 

We can now introduce the stronger notion of conical convergence of regular sequences in $G$ to chambers in $\DF X$, which first appeared in \cite{Albuquerque}. 

\begin{definition} [Conical convergence in $X$ to chambers at infinity]
\label{def:coniconv}
Let $\si\in\DF X$ be a chamber. 
A regular sequence $(x_k)$ in $X$ converges to $\si$ {\em conically} if there exists a constant $R$ and a euclidean Weyl chamber $V=V(x,\sigma)\subset X$ such that the sequence $(x_k)$ is contained in the 
$R$-neighbor\-hood of $V$. A sequence $(g_k)$ in $G$ converges to $\si\in \DF X$ conically if for some (equivalently, every) 
$x\in X$ the orbit sequence $(g_kx)$ converges conically to $\si$.  
\end{definition} 

Thus, conical convergence implies flag convergence, but the converse is false, as in rank 1.  We conclude with alternative formulations of  the  conical convergence $g_k\to \si$. 

\begin{prop}[Conical convergence criteria, \cite{morse,anolec}] 
\label{prop:equiv_conicality}

a. Suppose that a regular sequence $g_k\in G$ flag converges to $\si\in \DF X$. Then the following are equivalent:

1. $(g_k)$ converges conically to $\si$. 

2. For some (equivalently, every) point $x\in X$ and maximal flat $F\subset X$ whose visual boundary contains $\si$, 
the sequence of maximal flats $g_k^{-1}(F)$ is precompact in the space of all maximal flats in $X$. 

3. For some (equivalently, every) chamber $\hat\si\in \si^{opp}$, the sequence $g_k^{-1}(\si, \hat\si)$ is precompact in the space of antipodal pairs of chambers in $\DF X$. 

b. Conical convergence $g_k\to \si$ 
implies that the sequence $(g_k^{-1})$ has diverging infinitesimal expansion\footnote{See 
\S \ref{sec:expanding_actions} for the definition.}
at $\si$. 
\end{prop}

\bigskip 
\section{Discrete subgroups: Geometric and dynamical conditions} \label{sec:3}

As before, $G$ denotes the identity component of the isometry group of a symmetric space $X$ of noncompact type,
and $\Ga<G$ denotes a discrete subgroup.

We first define certain classes of discrete subgroups $\Ga<G$ within which we will be working 
throughout most of this paper,
namely of discrete subgroups which exhibit {\em rank 1 behavior} relative to 
(conjugacy classes of) parabolic subgroups $P< G$. 
We then discuss and compare various geometric and dynamical conditions for such subgroups.
As we noted earlier, 
in this survey we describe the theory for simplicity only in the (regular) case relative to minimal parabolic subgroups $P=B$. 
In the general ($\taumod$-regular) case, almost all the results go through with suitable modifications;
for the details we refer the reader to either one of the papers \cite{coco15, morse, anolec}. 
Among the major  differences in the $\taumod$-regular case are that one has to replace limit sets in the full flag manifold $G/B$ with limit sets in partial flag manifolds $G/P$, Weyl cones over chambers are replaced with suitable Weyl cones over 
stars of $\taumod$-type simplices, the expansion property occurs in the 
partial flag manifolds, various notions of regularity have to be modified and the Bruhat order on the Weyl group is replaced with orders on its coset spaces.

\subsection{Regularity and limit sets}

\begin{definition}
The 
{\em visual limit set} 
$$\La(\Ga)\subset\geo X$$ of $\Ga$ is the set of accumulation points of a $\Ga$-orbit
$\Ga x\subset X$ 
in the visual compactification $\ol X=X\sqcup\geo X$.
The elements of $\La(\Ga)$ are the 
{\em visual limit points}. Similarly, the {\em Finsler limit set} 
$$\La^{Fins}_{x}(\Ga)\subset \geo^{Fins}X$$ 
of $\Ga$ is the accumulation set of the orbit
$\Ga x$ in the Finsler compactification $\ol{X}^{Fins}=X\sqcup\geo^{Fins}X$. 
\end{definition}

\begin{rem}
While the visual limit set is independent of the orbit $\Ga x$, the Finsler limit set depends on it.
\end{rem}

\smallskip
We define regularity of subgroups as an asymptotic geometric condition on their orbits: 

\begin{definition}
[Regular subgroups \cite{coco13}]
A discrete subgroup $\Ga< G$ is {\em regular} 
(resp.\ {\em uniformly regular})
if each  divergent sequence of elements in $\Ga$ 
has this property
(cf.\ Defs.~\ref{def:regseq} and~\ref{def:regseqgp}). 
\end{definition}

Regularity can be read off the location of limit sets:

\begin{remark} 
1. Uniform regularity of $\Ga$ is equivalent to the property that the visual limit set 
consists only of regular ideal boundary points.

2. Regularity of $\Ga$ is equivalent to the property that the Finsler limit set $\La^{Fins}_x(\Ga)$ 
of some (equivalently, any) orbit $\Ga x$ 
is contained in the stratum $\DF X\subset \geo^{Fins} X$. 
\end{remark}

\begin{rem}
The notion of regularity for subgroups $\Ga< Isom(X)$ equally 
makes  sense when $X$ is a euclidean building.
The definition 
remains the same 
since for euclidean buildings one also has a $\De$-valued ``distance''  function $d_\Delta: X\times X\to \Delta$, 
where $\Delta$ is the model euclidean Weyl chamber of $X$,
see Appendix~\ref{sec:modsp}.
Most of the results mentioned in this survey go through without much change in the case when $X$ is a locally 
compact euclidean building. 
\end{rem}

\begin{definition}
[Chamber limit set \cite{coco13}]
The {\em chamber limit set}
$$
\La_{ch}(\Ga)\subset \DF X$$
consists of the chambers $\si\in \DF X$ for which there exists a sequence $(\ga_k)$ in $\Ga$ flag converging to $\si$,
$\ga_k\to\si$ (see Definitions~\ref{defn:flagconvergenceX} and~\ref{defn:flagconvergenceG}). 
The subgroup 
$\Ga< G$ is ($\simod$-){\em elementary}
if $\La_{ch}(\Ga)$ consists of at most 
two points. 
\end{definition}

\begin{rem}
1. More generally, in \cite[section 6.4]{coco15} we define the notion of $\taumod$-limit sets $\Lat(\Ga)\subset\Flagt=G/P_{\taumod}$ for discrete subgroups $\Ga< G$. 
(One has $\La_{ch}=\Las$.)

2. Benoist introduced in \cite[\S 3.6]{Benoist} 
a notion of limit set $\La_{\Ga}$ for Zariski dense subgroups $\Ga$ of reductive algebraic groups over local fields which in the case of real semisimple Lie groups 
is equivalent to our concept of chamber limit set $\La_{ch}$.\footnote{Benoist's limit set $\La_{\Ga}$ is contained in a partial flag manifold $Y_{\Ga}$ which in the case of real Lie groups is the full flag manifold $G/B$, see the beginning of \S 3 of his paper. 
In this case, $\La_{\Ga}$ consists of the limit points of the sequences in $\Ga$ contracting on $G/B$, cf.\ his Definitions 3.5 and 3.6.} What we call the $\taumod$-limit set $\Lat$ for other face types $\taumod\subsetneq\simod$ is mentioned in his Remark 3.6(3), and his work implies that, in the Zariski dense case, 
$\Lat$ is the image of $\La_{ch}$ under the natural projection $\DF X=\Flags\to\Flagt$ of flag manifolds.
\end{rem}

\begin{example}\label{ex:product-case}
Consider $X=X_1\times X_2$, the product of two real hyperbolic spaces, 
$g=(g_1, g_2)$ an infinite order isometry of $X$, where $g_1, g_2$ are isometries of $X_1, X_2$. Then the cyclic subgroup $\Ga=\<g\>$ is regular if and only if neither $g_1$ nor $g_2$ is elliptic. The subgroup $\Ga$ is uniformly regular if and only if both $g_1, g_2$ are hyperbolic 
isometries of $X_1, X_2$ or both are parabolic isometries. 
A cyclic group generated by an element of {\em mixed type} is {\em not} uniformly regular. 
The Furstenberg boundary of $X$ is the product $\geo X_1\times \geo X_2$. 
If $\Ga$, as above, is regular and $\la_i^+,\la_i^-$ are the fixed points of $g_i$
in $\geo X_i$,\footnote{I.e.\ $\la_i^+$ and $\la_i^-$ are the attractive and repulsive fixed points if $g_i$ is hyperbolic, 
and $\la_i^+= \la_i^-$ is the unique fixed point if $g_i$ is parabolic.}
then $\La_{ch}(\Ga)= \{(\la_1^-, \la_2^-),  (\la_1^+, \la_2^+)\}$. 
In particular, in the mixed case 
if, say, $g_1$ is hyperbolic and $g_2$ is parabolic with the unique fixed point $\la_2^+=\la_2^-=:\la_2$, then $\La_{ch}(\Ga)= \{(\la_1^-, \la_2),  (\la_1^+, \la_2)\}$. 
Note that if $\Ga$ is uniformly regular then the limit set $\La_{ch}(\Ga)$ is antipodal, but it is not antipodal if $\Ga$ is merely regular. 
The limit chambers are conical limit points if $g$ is uniformly regular of type hyperbolic-hyperbolic,
and otherwise they are not.
\end{example}

The next proposition gives alternative descriptions of chamber limit sets for regular and
uniformly regular subgroups in terms of Finsler and visual compactifications:

\begin{prop}[{\cite{bordif,coco13}}]
1. If $\Ga< G$ is regular then $\La_{ch}(\Ga)\subset \DF X$ is the Finsler limit set of $\Ga$. 
(In this case, it is independent of the $\Ga$-orbit.)

2. If $\Ga< G$ is uniformly regular then $\La_{ch}(\Ga)\subset \DF X$ is the set of 
chambers $\si\in \DF X$ 
which contain visual limit points.
(These are then contained in their interiors.)
\end{prop}

Let us mention in this context the following structural result for the visual limit set:

\begin{thm}
[{\cite[part of Thm.\ 6.4]{Benoist}}] For every Zariski dense discrete
subgroup $\Ga< G$ there exists 
an $\iota$-invariant closed convex subset $\ell(\Ga)\subset \simod$
with nonempty interior, 
such that for each chamber 
$\si\in\DF X$ satisfying $\inte(\si)\cap \La(\Ga)\neq\emptyset$ it holds that $\theta(\si\cap \La(\Ga))=\ell(\Ga)$.  
\end{thm}

Thus, in the case of uniformly regular Zariski dense 
subgroups $\Ga< G$, the visual limit set $\La(\Ga)$ is a 
$\Ga$-equivariant product bundle over $\La_{ch}(\Ga)$ with fiber $\cong \ell(\Ga)$. 

\medskip
In general, verifying (uniform) regularity of a subgroup is not an easy task. 
See e.g.\ \cite[Thm 3.51]{anolec} 
and Theorem \ref{thm:Anosov} of this paper 
for results of this kind. 
For Zariski dense subgroups the verification of regularity becomes easier.
The next result provides a sufficient condition:
\begin{thm}
Let $\rho: \Ga\to G$ be a representation whose image is Zariski dense in $G$. 
Suppose that $Z$ is a compact metrizable space, $\Ga\acts Z$ is a discrete convergence group action (with finite kernel),
and $f: Z\to \DF X$ is a $\rho$-equivariant topological embedding. 
Then $\rho$ has finite kernel and $\rho(\Ga)$ is regular.  
\end{thm}
\proof
In view of the Zariski density of $\rho(\Ga)$, also $f(Z)$ is Zariski dense in $\DF X$. 
Consequently, 
the assumption that $\Ga$ acts on $Z$ with finite kernel implies that $\rho$ has finite kernel. 

We assume that $\rho(\Ga)$ is not regular. 
We will be using certain notions and a proposition from 
\cite[\S 9.1.2]{bordif}. Given a simplex $\tau\in\Flagt$, the subvariety 
$$\stF(\tau) \subset\DF X$$
is the set of chambers $\si$ containing $\tau$ as a face. Similarly, for $\tau_-\in \Flagit$, 
$$\CF(\tau_-)\subset\DF X$$
is the Zariski open and dense subset equal to the union
$$
\bigcup_{\tau\in C(\tau_-)} \stF(\tau). 
$$

Suppose that for some sequence $\ga_i\to\infty$ in $\Ga$, 
the sequence $g_i=\rho(\ga_i)\in G$ is not regular. Hence, it contains a subsequence contained in a tubular neighborhood of the boundary of $\Delta$. Then, after extraction,  since $\Delta$ has only finitely many faces, 
 the sequence $(g_i)$ is {\em $\taumod$-pure} for some proper face $\taumod$ of $\simod$. This means that there exists a constant $D$ such that for each  $i$ the vectors $v_i:=d_\Delta(o, g_i(o))\in \Delta$ 
 belong to the $D$-neighborhood of some proper face $V_{\taumod}$ of $\Delta$ (the Weyl sector over the face  
 $\taumod$). Therefore, according to \cite[Prop. 9.4] {bordif}, 
after further extraction, there exists a pair of simplices 
 $\tau_+ \in \Flagt, \tau_-\in \Flagit$ 
 such that the sequence $(g_i)$ converges on the Zariski open and dense subset 
$\CF(\tau_-)\subset \DF X$ to a nonconstant {\em algebraic} map $\phi: \CF(\tau_-) \to \stF(\tau_+)$. 
Since $\phi$ is algebraic, it cannot be constant on a Zariski dense subset. 
On the other hand, by the convergence property on $Z$, after extraction, $(g_i)$ 
converges to a constant map on $f(Z) -\{f(z)\}$
for some (exceptional) point $z\in Z$.
A contradiction. \qed

\subsection{Generalized convergence subgroups}

It is useful to reformulate the concepts of (chamber) limit set and regularity for discrete subgroups 
purely in terms of their dynamics on $\DF X$. 

\begin{definition}
[Convergence subgroups \cite{coco15}]
\label{def:simod-convergence}
A discrete subgroup $\Ga < G$ is a {\em $\simod$-convergence\footnote{We add the prefix $\simod$ 
in order to distinguish from the notion of abstract {\em convergence group} in topological dynamics,
cf.\ Appendix~\ref{app:congru}.} subgroup} 
if for every divergent sequence of elements $\ga_k\in \Ga$, 
every subsequence of $(\ga_k)$ contains a further subsequence 
which converges to a quasiprojective map $\al_\om: \DF X\to \DF X$,
cf.\ section~\ref{sec:conprop}.
\end{definition}

Corollary~\ref{cor:regconveq} yields:
\begin{thm}[{\cite{morse,anolec}}]
\label{thm:regconvsbgp}
A discrete subgroup $\Ga< G$ is regular iff it is a $\simod$-convergence subgroup.
\end{thm}
Furthermore, 
the chamber limit set $\La_{ch}(\Ga)$ 
is the set of chambers  
$\al\in \DF X$ for which there exists a sequence $\ga_k\in \Ga$ such that $\ga_k\to \al_\om$ for some $\om \in \DF X$.

\medskip
We note that in \cite{coco15} we formulate a more abstract notion of generalized convergence actions of groups on topological spaces in terms of {\em accumulation of sequences $g_k\in G$ at subsets},  which covers $\simod$-convergence groups and their limit sets above. This more abstract notion explains why {\em balanced thickenings} 
(see Definition \ref{defn:thickenings} and \eqref{eq:thickening}) 
of limit sets appear naturally in the study of regular subgroups of $G$. We also note that the convergence type behavior in the sense of accumulation 
has been studied earlier by Karlsson, Papasoglu and Swenson 
in the general context of nonpositive curvature, see \cite[Thm.\ 1]{Karlsson} and \cite[Thm.\ 4]{PapaSwen}. 

\subsection{Rank 1 discrete subgroups $\Ga<G$} 

Currently, it appears that regularity (or even uniform regularity) alone is not enough to fully capture ``rank 1'' behavior of discrete subgroups $\Ga< G$. We introduce one extra condition on the chamber limit set:

\begin{definition}
A discrete subgroup $\Ga< G$ is {\em  antipodal} (A), if its limit chambers are pairwise antipodal.  
\end{definition}

We now can define a class of discrete subgroups which exhibit {\em rank 1 behavior:} 

\begin{definition}
A discrete subgroup $\Ga< G$ is {\em regular antipodal} (RA) or a {\em rank 1 discrete subgroup} of $G$
if it is regular and antipodal. 
\end{definition}

The higher rank convergence property (Definition \ref{def:simod-convergence}) then implies:

\begin{corollary}
If $\Ga$ is RA,
then the action $\Ga\acts\La_{ch}(\Ga)$ 
is an abstract convergence action.\footnote{See Appendix 
\ref{app:congru}.}
\end{corollary}
\proof Suppose that $(\ga_k)$ is a sequence of distinct elements in $\Ga$. 
In view of the regularity of $\Ga$, after extraction, $\ga_k\to \al_\om$ uniformly on compacts in $\om^{opp}$
for some limit chambers $\al,\om\in\La_{ch}(\Ga)$. 
Due to antipodality, $\La_{ch}(\Ga)- \{\om\}\subset \om^{opp}$.
Therefore, $\ga_k$ converges to $\al$ uniformly on compacts in 
$\La_{ch}(\Ga)- \{\om\}$. \qed 

\medskip 
Thus, restricting to the chamber limit set of an RA subgroup
brings us back to the familiar rank 1 setting!

\subsection{RCA subgroups}\label{sec:conical_convergence2}

We now begin discussing various geometric and dynamical conditions for regular discrete subgroups.
The first one concerns the asymptotic geometry of orbits,
namely whether limit chambers can be reached along orbits in a ``straight'' way:

\begin{definition}[Conical]
A limit chamber $\si\in \La_{ch}(\Ga)$ is {\em conical} 
if there exists a sequence 
$\ga_k\in\Ga$ such that $\ga_k\to \si$ conically,
cf.\ Definition~\ref{def:coniconv}.

A discrete subgroup $\Ga < G$ is {\em conical} if all its limit chambers are conical.
\end{definition}

\begin{theorem}
[Extrinsic  conicality is equivalent to intrinsic conicality \cite{morse,anolec}] 
For nonelementary RA subgroups $\Ga<G$, 
conicality is equivalent to intrinsic conicality in terms of the action $\Ga\acts\La_{ch}(\Ga)$, 
as defined in Appendix \ref{app:congru}. 
\end{theorem}

\begin{corollary} 
A nonelementary RA subgroup $\Ga<G$ is conical iff the action $\Ga \acts T\La_{ch}$ 
on triples of distinct limit chambers
is cocompact. 
\end{corollary}

We now arrive to the first definition of geometric finiteness in higher rank, generalizing the Beardon-Maskit definition:

\begin{definition}[RCA subgroups \cite{coco13, coco15}] 
A discrete subgroup $\Ga< G$ is {\em RCA} if it is regular, conical and antipodal. 
\end{definition}

\begin{remark} 
An analogous definition and theory exist in the $\taumod$-regular case.
One replaces the $\Ga$-action on $\DF X$ with the action on the partial flag manifold $\Flagt$. 
\end{remark}

Note that, a priori, it is unclear even why RCA groups are finitely generated. 
However, as a consequence of Bowditch's theorem \cite{Bowditch_char}
about the dynamical characterization of word hyperbolic groups,
one obtains: 

\begin{corollary}
Each nonelementary RCA subgroup $\Ga$ is word hyperbolic 
and its Gromov boundary $\geo \Ga$ is equivariantly homeomorphic to $\La_{ch}(\Ga)$. 
\end{corollary}

For RCA groups regularity is equivalent to uniform regularity:

\begin{theorem}
[RCA implies uniform regularity \cite{morse,anolec}] 
If $\Ga< G$ is nonelementary RCA then it is uniformly regular. 
\end{theorem}

\subsection{Expansion at infinity}

The RCA condition discussed above is in terms 
of the asymptotics of the group action on the symmetric space $X$. 
Our next definition is in terms of the dynamics at infinity, more precisely, of the action 
on the Furstenberg boundary $\DF X$.

\begin{definition}[CEA subgroup \cite{coco15,anolec}] 
A discrete subgroup $\Ga < G$ is {\em CEA} (convergence, expanding, antipodal) if:

1. $\Ga< G$ is a $\simod$-convergence subgroup. 

2.  The action $\Ga\acts\DF X$ is expanding at $\La_{ch}(\Ga)$
in the sense of Appendix~\ref{sec:expanding_actions}. 

3. The chamber limit set $\La_{ch}(\Ga)$ is antipodal. 
\end{definition}

We recall that the convergence condition is equivalent to regularity, see Theorem~\ref{thm:regconvsbgp}.\footnote{We 
impose the convergence instead of the regularity condition 
in order to make the notion purely dynamical.}

For nonelementary subgroups $\Ga<G$,
the second condition is satisfied 
if the restricted action $\Ga\acts\La_{ch}(\Ga)$
is expanding,
compare Theorem~\ref{thm:conical} in Appendix~\ref{app:congru}.

Expansivity is useful for proving cocompactness of the $\Ga$-action on domains of discontinuity, see \cite{coco13}. 
We will come back to this later.

\subsection{Asymptotically embedded subgroups and Anosov representations}

The next condition we give is in terms of boundary maps into $\DF X$.
It requires the discrete subgroups $\Ga<G$ to be intrinsically word hyperbolic,
unlike our earlier conditions where hyperbolicity was a consequence.

We will consider boundary maps of the following kind:

\begin{definition}
A map into $\DF X$ is {\em antipodal} 
if it sends distinct points to antipodal chambers. 
\end{definition} 

\begin{remark} 
1. Antipodal maps are injective. 

2. An antipodal continuous map is the same thing as a {\em transversal map} in the sense of \cite{GW}. 
\end{remark}

\begin{definition}[Asymptotically embedded subgroups \cite{coco13, anolec}] 
A discrete subgroup $\Ga<G$ is {\em asymptotically embedded}  
if it is RA, intrinsically word hyperbolic, and if there exists a $\Ga$-equivariant homeomorphism 
$$\be: \geo \Ga \stackrel{\cong}{\lra} \La_{ch}(\Ga)\subset\DF X.$$
\end{definition}

Note that the boundary map $\beta$ is necessarily antipodal in this case. 
Furthermore,
any orbit map $\Ga\to X$
continuously extends by $\beta$ 
to a map 
$$\Ga\sqcup\geo\Ga\to X\sqcup \DF X$$
from the visual (Gromov) compactification of $\Ga$ into the partial compactification $X\sqcup \DF X$
equipped with the topology of flag convergence\footnote{Compare \S \ref{sec:flagcv}.},
see \cite[Prop 3.20]{anolec}. 

\medskip
Below we present two related notions,
{\em boundary embedded subgroups} and {\em Anosov subgroups}. 
Instead of requiring an identification of the Gromov boundary with the chamber limit set, we can at first require only the existence of an 
equivariant antipodal embedding into $\DF X$:

\begin{definition}[Boundary embedded subgroups \cite{coco13, anolec}] 
A discrete subgroup $\Ga< G$ is {\em boundary embedded} 
if it is intrinsically word hyperbolic and there exists an equivariant antipodal continuous map  
$$
\beta': \geo \Ga \to \DF X,
$$
called a {\em boundary embedding} for $\Ga$.
\end{definition}

\begin{remark} 1. {\em Boundary embedded} is the topological part of the definition of {\em Anosov} subgroups 
in \cite{Labourie,GW}.
The {\em dynamical} part is omitted in this definition. 

2. Boundary embeddings are in general {\em not unique}.
This is so by trivial reasons if $|\geo\Ga|=2$, 
but it also happens when $|\geo\Ga|\geq3$, see e.g.\ \cite[Example 6.20]{morse}.
\end{remark}

There is the following {\em dichotomy} for the relation of boundary embeddings with the chamber limit set:

\begin{thm}[Boundary embedding dichotomy {\cite[Thm.\ 3.11]{anolec}}]
Suppose that $\Ga< G$ is a boundary embedded regular subgroup. 
Then for each boundary embedding $\beta': \geo \Ga \to \DF X$ we have the following dichotomy:

1. Either $\beta'(\geo \Ga)= \La_{ch}(\Ga)$, or

2. $\beta'(\geo \Ga)\cap \La_{ch}(\Ga)=\emptyset$ and
$$
\La_{ch}(\Ga)\subset \bigcap_{\si\in \beta'(\geo \Ga)} (\DF X - \si^{opp}). 
$$
In the first case, $\Ga$ is {asymptotically} embedded,
while the second alternative implies that $\La_{ch}(\Ga)$ is contained in a proper subvariety of $\DF X$. 
\end{thm}

The second alternative cannot occur in the Zariski dense case. Therefore:

\begin{cor}
If $\Ga< G$ is regular, Zariski dense and $\beta': \geo \Ga \to \DF X$ is a boundary embedding,
then 
$\be'(\geo \Ga)=\La_{ch}(\Ga)$.  In particular, $\Ga< G$ is asymptotically embedded. 
\end{cor}

We note that the last part of this corollary was already proven in \cite{GW}. 
The next theorem shows that one does not need Zariski density in order to conclude that 
 $\Ga< G$ is asymptotically embedded. 

\begin{thm}[{\cite{morse} and \cite[Thm.\ 3.15]{anolec}}]
\label{thm:bae}
A regular subgroup $\Ga < G$ is boundary embedded 
iff it is asymptotically embedded. 
\end{thm}

\begin{remark} 
1. In general, there may exist several boundary embeddings for $\Ga$,
and only one of them yields the asymptotic embedding. 

2. Theorem \ref{thm:bae} is one of the few results in the theory 
which hold only in the regular case (as opposed to the $\taumod$-regular case). 
\end{remark}

Now we give our versions, see \cite{morse,anolec}, 
of the definition of Anosov subgroups, 
which were originally defined 
in \cite{Labourie,GW}
using expansion properties of geodesic flows.
Our definitions do not use geodesic flows of word hyperbolic groups 
but replace them by a simpler coarse geometric object,
namely by the space of discrete geodesics with respect to a word metric.
These were the first such definitions which are close to Anosov in spirit 
but do not use flows.  

We fix a Riemannian metric on $\DF X$.
Moreover, we will assume in the next definition that the word hyperbolic group $\Ga$ is equipped with a fixed word metric.

\begin{definition}
[Anosov subgroups, {\cite[\S 6.5]{morse}}]
1. A subgroup $\Ga< G$ is {\em  Anosov} if it is boundary embedded
with boundary embedding $\be'$ and, in addition, 
for each normalized\footnote{$r(0)=1\in \Ga$} discrete geodesic ray $r: k\mapsto \ga_k\in \Ga$ asymptotic to $\xi\in \geo \Ga$ 
the sequence $(\ga_k^{-1})$ is uniformly exponentially infinitesimally expanding at $\be'(\xi)\in\DF X$. 
More precisely, there are constants $C,A>0$ 
depending only on the subgroup  $\Ga< G$,
the word metric on $\Ga$ 
and the Riemannian metric on $\DF X$
such that 
$$
\eps(\ga_k^{-1}, \be'(\xi)) \ge A e^{Ck}
$$ 
for $k\geq0$. 
Here, $\eps$ is the infinitesimal expansion factor 
defined in Appendix~\ref{sec:expanding_actions}.

2. A subgroup $\Ga< G$ is {\em non-uniformly Anosov} 
if it is boundary embedded with boundary embedding $\be'$
and, in addition, 
for each discrete geodesic ray $r: k\mapsto \ga_k\in \Ga$ asymptotic to $\xi\in \geo \Ga$, the sequence 
$(\ga_k^{-1})$ contains a subsequence with diverging infinitesimal expansion at $\be'(\xi)\in\DF X$,
$$
\sup_{k\in\N}\eps(\ga_k^{-1}, \be'(\xi)) =\infty.
$$ 
\end{definition}

Note that due to the stability of quasigeodesics in word hyperbolic groups,
the definition is independent of the word metric on $\Ga$.

\begin{thm}[\cite{morse,anolec}]\label{thm:Anosov}
Suppose that $\Ga<G$ is intrinsically word hyperbolic 
and not virtually cyclic. 
Then: 

1. $\Ga< G$ is non-uniformly Anosov iff $\Ga< G$ is Anosov iff $\Ga< G$ is asymptotically embedded. 

2. If these conditions are satisfied, 
then the boundary maps for the Anosov and asymptotic embeddedness conditions
coincide.
\end{thm}

\begin{rem}
Note that the original Anosov condition 
in \cite{Labourie,GW}
involves the space ${\mathcal G}$ 
of (equivalence classes) of all {\em parameterized} geodesics in $\Ga$, equipped with a suitable topology. 
This space admits two commuting actions: 
a left action of $\Ga$ and a right action of $\R$ (shifting geodesics). Let 
$(\DF X\times\DF X)^{opp}$ denote the subset of $\DF X\times\DF X$ consisting of pairs of opposite chambers. We regard the product space  
$$
{\mathcal B}:= {\mathcal G} \times (\DF X\times\DF X)^{opp}
$$
as a trivial bundle over ${\mathcal G}$; then the boundary map $\be: \geo \Ga\to\DF X$ 
defines a section of this bundle which projects to a section $s_\be$ of the quotient bundle
$$
\Ga\backslash {\mathcal B} \to \Ga\backslash {\mathcal G}. 
$$  
The commuting actions of 
$\Ga$  and $\R$ lift to commuting actions on 
${\mathcal B}$, where $\R$ acts trivially on the second factor, while $\Ga$ acts on the second factor via the restriction of the natural product action  on $\DF X$.  
The original Anosov axiom amounts to an expansion/contraction condition (along $s_\beta$) 
for the {\em right} $\R$-action on $\Ga\backslash {\mathcal B}$. The basic dynamical duality principle suggests that this condition can be reinterpreted as expansion/contraction property for the {\em left} $\Ga$-action on ${\mathcal B}/\R$: 
This is what the our interpretation of the Anosov property amounts to (after a careful rewriting of the definitions involved),
see \cite[\S 6.5]{morse} for a detailed discussion.  
\end{rem}

\subsection{URU subgroups}

The next set of definitions is in terms of 
extrinsic {\em coarse geometric properties}. 

A finitely generated subgroup $\Ga<G$
is said to be {\em undistorted} 
if for some (every) point $x\in X$ 
the orbit map
$$
\ga\mapsto \ga x\in X
$$
is a quasiisometric embedding $\Ga\to X$, where $\Ga$ is equipped 
with a word metric. 
Equivalently,
the inclusion $\Ga\subset G$ is a quasiisometric embedding.

\begin{definition}[URU subgroups, {\cite{mlem,bordif}}]
A finitely generated discrete subgroup $\Ga < G$ is {\em URU} 
if it is uniformly regular and undistorted. 
\end{definition}

\begin{remark} 
There are regular undistorted subgroups which are not uniformly regular.
Take, for instance,
the cyclic groups in Example \ref{ex:product-case} 
where $g_1$ is hyperbolic and $g_2$ is parabolic. 
There are also finitely generated nonabelian free subgroups of this kind in $PSL(2,\R)\times PSL(2,\R)$; 
moreover, some of such subgroups are not even $P$-Anosov for any $P$, 
see \cite[Appendix A]{GGKW1}. 
Similar examples also exist among closed surface subgroups of  
$PSL(2,\C)\times PSL(2,\C)$ and closed surface subgroups of $\Isom(T\times T)$, where $T$ is a simplicial tree, see \cite{KL-undistorted}.  
\end{remark}

The next condition imposes a ``Morse'' property on the images of discrete geodesics in $\Ga$.
A priori stronger than URU,
it is equivalent to it because of the Higher Rank Morse Lemma
and can be viewed as describing the extrinsic coarse geometry of URU subgroups.

\begin{dfn}[Morse subgroups, {\cite{morse}}]
\label{def:morssbgp}
A finitely generated discrete subgroup $\Ga<G$ is {\em Morse}
if some (equivalently, every) orbit map $\Ga\to X$ is Morse, 
cf.\ Definition~\ref{def:mrs}.
\end{dfn}

The next condition is motivated by the Finsler geometric interpretation of the Morse Lemma, 
see Corollary~\ref{cor:mlemfins},
that uniformly regular quasigeodesics in $X$ are Finsler quasiconvex.
Here, 
a subset $A\subset X$ is called {\em Finsler quasiconvex}
if there exists a constant $R>0$
such that for any pair of points $x_1,x_2\in A$
there {\em exists} a Finsler geodesic\footnote{I.e. a geodesic with respect to a fixed regular polyhedral Finsler metric $d_{\bar\theta}$ on $X$.} 
from $x_1$ to $x_2$ contained in the $R$-neighborhood of $A$.

\begin{definition}[Finsler quasiconvex subgroups, {\cite{bordif}}] 
A subgroup $\Ga<G$ is {\em Finsler quasiconvex} 
if some (equivalently, every) $\Ga$-orbit $\Ga x\subset X$ is Finsler quasiconvex.
\end{definition}

This notion mimics the notion of quasiconvexity for discrete subgroups of rank 1 Lie groups 
discussed in section \ref{sec:GFG}. 
The key difference with rank 1 is that, when $\rank(X)\geq2$, 
Finsler geodesics connecting pairs of points in $X$ are no longer unique.

\subsection{Equivalence of conditions}

We can now put together a theorem which states 
the equivalence of various geometric and dynamical notions of geometric finiteness for 
discrete isometry groups of symmetric spaces exhibiting rank 1 behavior. 
This theorem is a combination of results 
of \cite{morse}  and \cite{mlem}. It will be augmented by two more equivalent notions in section \ref{sec:bordif} 
(Corollary~\ref{cor:S-coco}).

\begin{theorem}[Equivalence]
\label{thm:main}
For  discrete subgroups $\Ga< G$ the following conditions are equivalent
in the nonelementary\footnote{Here, ``nonelementary'' means $|\geo\Ga|\geq3$ in the Anosov conditions 5 and 6, 
which assume word hyperbolicity but no regularity, 
and means $|\La_{ch}(\Ga)|\geq3$ in all other cases.} 
case:

1. $\Ga< G$ is RCA. 

2. $\Ga < G$ is CEA. 

3. $\Ga< G$ is asymptotically embedded. 

4. $\Ga< G$ is boundary embedded.\footnote{This, unlike the other equivalences, 
is limited to  $\simod$-regular subgroups $\Ga <G$.}  

5. $\Ga< G$ is Anosov.  

6. $\Ga< G$ is non-uniformly Anosov.  
 
7. $\Ga < G$ is Morse. 

8. $\Ga < G$ is URU. 

9. $\Ga<G$ is uniformly regular and Finsler quasiconvex. 
\end{theorem}

The most difficult step in the proof of this theorem is from URU to Morse: 
It follows from the Morse Lemma for uniformly regular quasigeodesics 
and the companion results on hyperbolicity and boundary maps  in \cite{mlem}. 

\begin{rem}
The {\em nonelementary} assumption in the theorem most likely could be dropped. 
The reason for including it 
is that it is currently unknown if there are RCA (or CEA) subgroups $\Ga< G$ with $\La_{ch}(\Ga)$ consisting of a single point. 
\end{rem}

\begin{rem}[Relation with the paper \cite{GGKW1}]
A weaker form of the equivalence of the conditions {\em Anosov} and {\em URU} 
was established in \cite[Thm.\ 1.3]{GGKW1} 
after \cite{morse,mlem} had been available.
There are two major differences:
First,
the discussion in \cite{GGKW1} is restricted to word hyperbolic subgroups,
while URU only assumes finite generation.
Second,
the URU condition is replaced in \cite{GGKW1} with the (a priori) stronger ``CLI''  (coarse linear increase) condition,
see \cite[Thm 1.3(iv)]{GGKW1}.
The difference in character between the URU and CLI conditions is, roughly,
like the difference between asymptotic linear growth and quasiisometric embedding 
for Lipschitz maps $\N\to\R_+$;
the former implies the latter but not conversely.
More precisely,
to describe the difference between URU and CLI
for a finitely generated subgroup $\Ga<G$,
fix a word metric on $\Ga$ and a point $x\in X$. 
Then consider for all discrete geodesic rays $r:\N\to\Ga$ normalized by $r(0)=e$
their $\De$-distance projections 
\begin{equation*}
\bar r:=d_{\De}(x,rx):\N\to\De .
\end{equation*}
The subgroup $\Ga<G$ is URU iff 
the paths $\bar r$ are drifting away from $\D\De$ 
at a {\em uniform linear rate}, in a coarse sense.\footnote{In the general $\taumod$-regular case,
$\D\De$ is replaced with $\D_{\taumod}\De=V(0,\D_{\taumod}\simod)$,
the union of the walls of $\D\De$ not containing the sector $V(0,\taumod)$.}
This is equivalent to the {\em uniform linear growth} of $\al\circ \bar r$ for all {\em simple} roots $\al$.
On the other hand,
$\Ga<G$ is CLI 
iff the $\al\circ \bar r$ {\em uniformly increase}, again in a coarse sense, 
for all {\em positive} roots $\al$.\footnote{In the $\taumod$-regular case, 
one takes the simple, respectively, positive roots
which do not vanish on $V(0,\taumod)$,
equivalently, which are nonnegative on the symmetrized cone $\Wt\De\subset\Fmod$.
Here $\Wt<W$ denotes the stabilizer of $\taumod$.}
There is the following geometric interpretation of the CLI condition from the Morse viewpoint:
CLI is equivalent to the $\bar r$ being {\em uniform Morse quasigeodesic rays}.
For arbitrary Lipschitz paths $\N\to\De$,
the linear drift condition is strictly weaker than the Morse condition.
In particular, URU follows from CLI and, on the face of it, appears weaker.
However, 
it is not hard to see that the paths $\bar q:=d_{\De}(q(0),q):[0,T]\to\De$ 
coming from uniform Morse quasigeodesics $q:[0,T]\to X$ are themselves uniform Morse quasigeodesics in $\De$.
In particular, CLI is a consequence of Morse.
Thus, 
\cite[Thm.\ 1.3]{GGKW1} and in particular the implication {\em Anosov$\Rightarrow$CLI} 
follow from \cite{morse,mlem}
(that also the latter implication follows 
is stated as unclear in \cite[\S 1.3 of version 5]{GGKW1}).
On the other hand,
the implication {\em URU$\Rightarrow$Anosov}, 
which is based on our Higher Rank Morse Lemma \cite{mlem}, 
does not follow from \cite{GGKW1}.
\end{rem}

\subsection{Consequences} 

In this section, we briefly discuss some properties shared by the groups 
satisfying (one of) the conditions listed in Theorem~\ref{thm:main}.

In some cases, it is more natural to talk about representations rather than subgroups.
We call a representation $\Ga\to G$ of a word hyperbolic group $\Ga$ a
{\em Morse representation} if some (every) orbit map $\Ga\to X$ is Morse, 
see \cite{morse}.

\bigskip
{\bf 1. Local-to-global principle.} 
The local-to-global principle for Morse quasigeodesics (Theorem~\ref{thm:L2G}) 
implies one for Morse representations:

\begin{thm}
[Local-to-global principle for Morse representations \cite{morse}] 
\label{thm:ltgpmrp}
Suppose that $\rho: \Ga\to G$  is a representation of a word hyperbolic group $\Ga$
such that some orbit map $\Ga\to X$, restricted to a sufficiently large ball in the Cayley graph
(say, centered at $1$),
is {\em locally Morse} of sufficiently good quality. 
Then $\rho$ is a Morse representation.   

More precisely, given Morse data $(\Theta, B, L, A, D)$ and the scale $S$ determined by them via Theorem~\ref{thm:L2G}
(for some $\Theta'$),
if the orbit map 
$\Ga\to X$ 
sends discrete geodesic segments of length $\leq S$ in the Cayley graph 
(say, passing through $1$) 
to $(\Theta, B, L, A, D)$-Morse quasigeodesic segments in $X$,
then  $\rho$ is a Morse representation. 
\end{thm}

{\bf 2. Structural stability.}
Structural stability was first established for convex-cocompact subgroups of rank 1 Lie groups by Sullivan, 
see section \ref{sec:corGF}. The following theorem is a generelization of Sullivan's result: 

\begin{thm}[Structural stability of Morse representations]
\label{thm:strstb}
For a word hyperbolic group $\Ga$,
the space of Morse representations $\rho:\Ga\embed G$ 
is an open subset of $\Hom(\Ga,G)$.
On this subset, 
the boundary map 
$\beta_{\rho}: \geo \Ga \to  \La_{ch}(\rho(\Ga))$ 
depends continuously on the Morse representation $\rho$.
\end{thm}

Here, the representation space $\Hom(\Ga, G)$ 
is equipped with the topology of pointwise convergence,
i.e. $\lim_{n\to\infty} \rho_n=\rho$ iff for every element $\ga$ of $\Ga$  
$$
\lim_{n\to\infty} \rho_n(\ga)=\rho(\ga). 
$$

\begin{remark} Structural stability of Anosov representations has first been proven in \cite{Labourie} 
for fundamental groups of closed negatively curved manifolds, 
and in \cite{GW} for general word hyperbolic groups. 
Our proof in \cite{morse} 
derives structural stability,
from the Morse viewpoint, 
as a direct consequence of the local-to-global principle for Morse maps (see Theorem \ref{thm:L2G}).
\end{remark}

\medskip
{\bf 3. Semidecidability.}\footnote{See section \ref{sec:corGF} for the definition of semidecidability.} 
It is {\em semidecidable} if a representation $\rho: \Ga\to G$ 
of a word hyperbolic group $\Ga$ is Morse.

\medskip 
The proof of semidecidability  given in \cite{morse} is also based on the local-to-global principle for Morse maps: The algorithm explores finite subsets $F$ of the Cayley graph of $\Ga$ and ranges of Morse data $(\Theta, B,L,A,D)$ to determine if an orbit map $\Ga\to X$ is $(\Theta, B,L,A,D)$-Morse on $F$. 

\medskip
{\bf 4. Cocompactness.} Each Morse subgroup $\Ga< G$ acts properly discontinuously and cocompactly on various domains associated with the action $\Ga\acts X$. These domains are contained in the flag manifold $\DF X= G/B$ and in the Finsler compactification of $X$. We refer to sections \ref{sec:dd} and \ref{sec:bordif} for the precise statements. 

\subsection{Examples: Morse-Schottky subgroups}

Let $\Ga$ be a free group on $k$ generators, denoted $\al_1, \al_2$,..., $\al_k$. 
Realize each $\al_i$ as a {\em regular} hyperbolic isometry $g_i$ of $X$,
i.e.\ $g_i$ preserves a regular geodesic line and translates along it.
Assume furthermore that the isometries $g_i$ are in {\em general position} with respect to each other, 
in the sense that the subset of $\DF X$ consisting of the $2k$ attractive and repulsive chambers 
$\si_1^\pm,...,\si_k^\pm$ of the isometries $g_1,...,g_k$ is antipodal. 

The following theorem was proven in \cite[Thm. 7.40]{morse} for $k=2$ (2-generated free groups), but the same proof goes through for arbitrary $k\in\N$. 

\begin{thm}\label{thm:MS-actions}
There exists $N_0$ such that for all $N\ge N_0$
$$\rho: \al_i\mapsto g_i^N, \quad i=1,...,k,$$ 
defines a (faithful) Morse representation $\rho:\Ga\to G$. 
\end{thm}

\begin{rem}
(i) Regarding earlier work on the construction of free subgroups of Lie groups,
note that Tits, when proving the Tits alternative using his ping-pong argument,
only shows the {\em injectivity} of certain representations of free groups, 
although his proof clearly also implies the {\em discreteness} of their images. 
Benoist \cite{Benoist} improved on Tits' result and obtained control on the {\em asymptotic geometry}.
In particular, 
he produced discrete free subgroups which, in our terminology, are {\em uniformly regular}.
Our construction \cite[Theorem 7.40]{morse} 
is the first to control the {\em coarse geometry}.
We prove that the resulting free subgroups are {\em Morse},
which amounts to describing the extrinsic coarse geometry of their orbits
(see Definitions~\ref{def:mrs}(i) and~\ref{def:morssbgp}).
In particular, they are {\em undistorted}.

Whereas the arguments of Tits and Benoist use the dynamics at infinity,
our approach is different.
We work inside the symmetric space and 
build up representations of free groups
using a version of the local-to-global principle for Morse representations,
see Theorem~\ref{thm:straight-paths} below. 

(ii) In \cite[Theorem~7.40]{morse} we prove a more general version of Theorem~\ref{thm:MS-actions}
which allows for $\taumod$-regular generators.
\end{rem}

Our proof of Theorem~\ref{thm:MS-actions}
is based on the notion of {\em straight paths} in symmetric spaces.\footnote{In fact,
the notion of straight paths and Theorem~\ref{thm:straight-paths} below 
are the key technical tools for proving our local-to-global principles 
Theorems~\ref{thm:L2G} and ~\ref{thm:ltgpmrp} for Morse quasigeodesics and Morse subgroups.}
This concept is a higher rank analogue of a piecewise geodesic path in a rank 1 symmetric space,  where each edge is sufficiently long and the vertex angles are close to $\pi$;
such paths are known to be uniformly quasigeodesic. In higher rank, the angle condition has to be suitably modified,
in order to make sure that the path bends ``transversally to the flat directions''. 
Below is the precise definition. 

Let 
$$
x_0 x_1 \ldots x_n$$ 
be a piecewise geodesic path in $X$ with vertices $x_i$. 
We call such a path {\em $s$-spaced} if 
$$
d(x_{i-1}, x_i)\ge s
$$
for all $i$.
In order to define {\em straightness},
we consider the chain of diamonds 
$$
D_i= \diamo_{x_{i-1} x_{i}}
$$
associated to our path
and define {\em $\bar\zeta$-angles} 
$$
\angle^{\bar\zeta}(D_{i-1}, D_i)
$$
between consecutive diamonds as follows.
As an auxiliary datum, we fix an $\iota$-invariant regular type $\bar\zeta\in \inte(\simod)$.
For every regular segment $xy$, respectively, for the associated diamond 
$$
\diamo_{x y}= V(x,\si)\cap V(y, \hat\si)
$$
we define the tangent vector $v_{xy} \in T_x V(x,\si)$  as the unique unit vector of type $\bar\zeta$,
i.e. such that the geodesic ray from 
$x$ in the direction $v_{xy}$ is asymptotic to the point $\zeta\in \si$ of type $\bar\zeta$, 
$\theta(\zeta)=\bar\zeta$. 
Then define the $\bar\zeta$-angle 
between two diamonds $\diamo_{xy}, \diamo_{xz}$ at $x$
as the Riemannian angle
$$
\angle^{\bar\zeta}(\diamo_{xy}, \diamo_{xz}) :=\angle(v_{xy}, v_{xz}).
$$

\begin{definition}
Let $\eps>0$ and $\Theta \subset \inte(\simod)$ be compact convex.
The piecewise geodesic path $x_0 x_1 \ldots x_n$ is called $(\Theta,\eps)$-{\em straight} if for all $i$ the segments
$x_{i-1} x_i$ are $\Theta$-regular 
and
$$
\angle^{\bar\zeta}(D_{i-1}, D_i)\ge \pi- \eps. 
$$
\end{definition}

\begin{figure}[tbh]
\includegraphics[width=90mm]{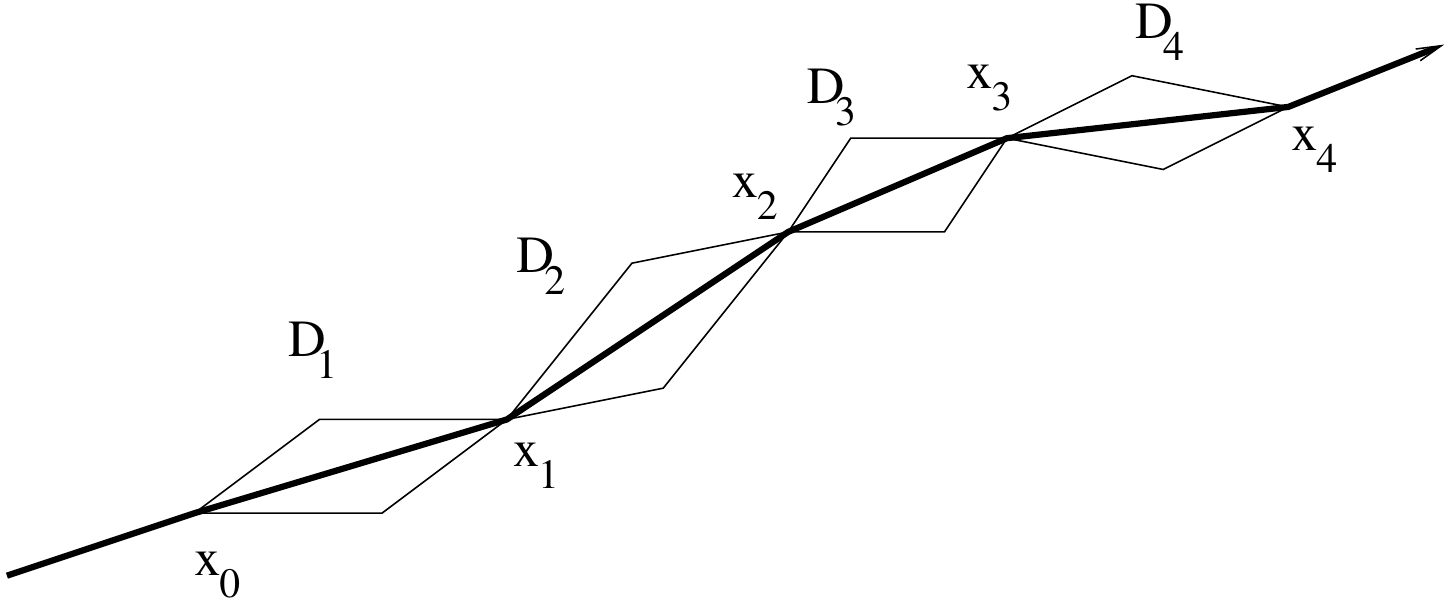}
\caption{A  string of diamonds.}
\label{figure12.fig}
\end{figure}

Now we can formulate:
\begin{thm}
[Local-to-global principle for straight paths \cite{morse}]
\label{thm:straight-paths}
Each sufficiently spaced and sufficiently straight 
piecewise geodesic path in $X$
is a uniform Morse quasigeodesic. 

More precisely, 
given 
compact convex subsets $\Theta, \Theta'\subset \inte(\simod)$ with 
$\Theta\subset \inte(\Theta')$, 
there exist numbers $\eps, s, B, L, A, D>0$,  
depending also on $\bar\zeta$,
such that each $s$-spaced and $(\Theta,\eps)$-straight piecewise geodesic path 
is $(\Theta',B,L,A,D)$-Morse. 
\end{thm}

Now we return to the setup of Theorem~\ref{thm:MS-actions}
and apply this local-to-global result to construct Morse representations of free groups.

We let $T$ denote the Cayley tree of $\Ga$ associated with the generating set $\{\al_1,...,\al_k\}$;
its vertex set is identified with $\Ga$.
For a point $x\in X$,
we extend the orbit map $o_x:\Ga\to\Ga x\subset X, \ga\mapsto\rho(\ga)x$,
to a piecewise geodesic map $f_x: T\to X$
of the Cayley tree 
by sending its edges to geodesic segments in $X$

We would be done if we could arrange $f_x$ to be straight in the sense that 
it maps lines in $T$ to $s$-spaced and $(\Theta,\eps)$-straight paths in $X$ for good data $s,\Theta$ and $\eps$,
because these image paths would be uniformly Morse by Theorem~\ref{thm:straight-paths},
which means that the representation $\rho$ would be Morse.
However, this is impossible to arrange if $k\geq2$, 
due to the ``lack of space'' in the unit tangent spheres:
The $2k$ image edges connecting the orbit point $x$
to the adjacent orbit points $\rho(\beta)x$, $\beta\in\{\al_1^{\pm1},\ldots,\al_k^{\pm1}\}$,
cannot have pairwise $\bar\zeta$-angles close to $\pi$,
equivalently, 
the $2k$ directions $v_{x\rho(\beta)x}$ at $x$ cannot have pairwise Riemannian angles close to $\pi$.\footnote{In 
euclidean buildings, 
it is easy to construct straight piecewise geodesic trees.}

We circumvent this difficulty by looking at {\em midpoint paths}:
Let $l\subset T$ be a line passing through the sequence of consecutive vertices
$\ldots,\ga_{-1},\ga_0,\ga_1,\ldots$.
Its $f_x$-image is the biinfinite piecewise geodesic path
$$ \ldots x_{-1}x_0x_1\ldots$$
with vertices at the orbit points $x_i=\rho(\ga_i)x$.
Let $m_i$ denote the midpoint of the segment $x_{i-1}x_i$
and consider the {\em midpoint path}
$$ \ldots m_{-1}m_0m_1\ldots$$
We are again done if we can show that these midpoint paths for all lines $l\subset T$ are uniformly
well spaced and straight.

This approach works 
and it is how our proof of Theorem \ref{thm:MS-actions} proceeds:
The point $x\in X$ can be chosen arbitrarily. 
We show that for a suitable compact convex subset $\Theta\subset \inte(\simod)$ and arbitrary $\eps,s>0$ 
the midpoint paths for all lines $l\subset T$ 
are $s$-spaced and $(\Theta,\eps)$-straight, provided that $N$ is sufficiently large. 

The $s$-spacedness for large $N$ easily follows from our genericity assumption 
that the chambers $\si^\pm_1,...,\si^\pm_k$ are pairwise antipodal. 
Due to $\Ga$-equivariance, the $(\Theta,\eps)$-straightness condition 
can be verified locally by looking at special short midpoint paths: 
For every triple of generators 
$$
\al, \be, \ga\in \{\al_1^{\pm 1},...,\al_k^{\pm 1}\}
$$
with $\al\ne \be, \be\ga\ne 1$ consider the quadruple 
$$
(\ga_0,\ga_1,\ga_2,\ga_3):=(\al, 1,\be, \be\ga)
$$ 
of elements in $\Ga$.
Then the $\ga_0\ga_1\ga_2\ga_3$ are geodesic paths in $T$,
and it suffices to check $(\Theta,\eps)$-straightness for the associated midpoint paths $m_0m_1m_2$. 
The latter is deduced from smallness of $\bar\zeta$-angles, 
$$
\angle^{\bar\zeta}(\diamo_{m_1m_0}, \diamo_{m_1x}) =\angle(v_{m_1m_0}, v_{m_1x}) <\frac{\eps}{2},
$$
and sufficient spacing (to ensure the regularity).
The smallness of angles and the regularity are verified 
by a direct geometric argument using the regularity of the elements $g_i$ and their general position. 
We refer the reader to \cite[sect. 7.6]{morse} for the details.  

\begin{figure}[tbh]
\includegraphics[width=90mm]{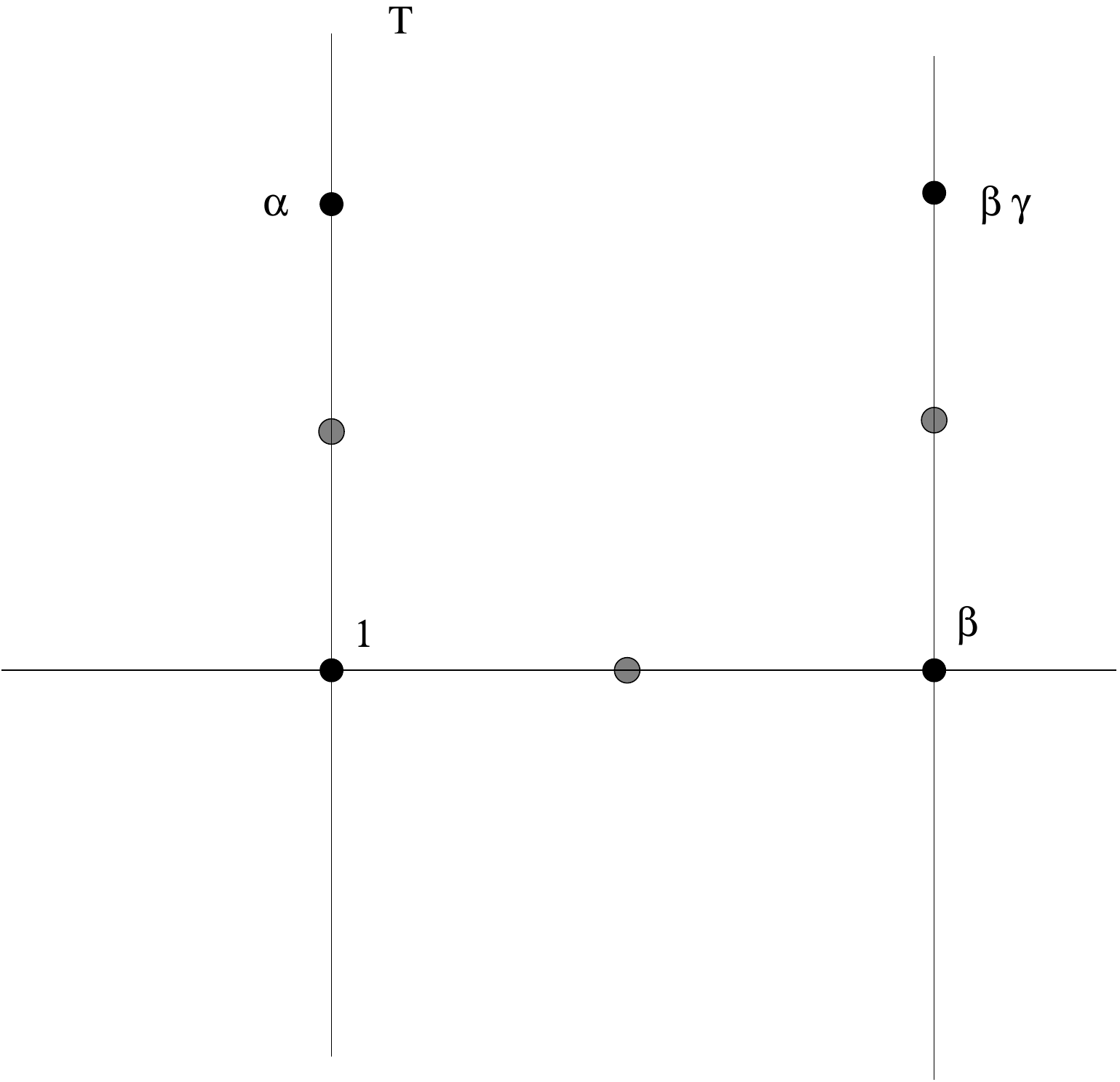}
\caption{Special short midpoint paths.}
\label{figure13a.fig}
\end{figure}

\medskip 
\subsection{Further examples} 

Other examples of Morse, equivalently, Anosov subgroups are provided by Hitchin representations, see \cite{Labourie}, 
which are the origin of the notion of Anosov representations. Similarly, one obtains the {\em complex version} of these examples: Start with the (unique up to isomorphism) irreducible representation
$$
\rho_n: SL(2,\C)\to G=SL(n,\C). 
$$
Then for each discrete subgroup $\Ga< SL(2,\C)$ its image $\rho_n(\Ga)=\Ga_n< G$ is an RA subgroup. 
If, in addition, $\Ga$ is convex-cocompact, then $\Ga_n$ is RCA, equivalently, Anosov. 
Due to structural stability, 
any representation $\rho$ sufficiently close to $\rho_n$, is also Anosov.
In the case when $\H^3/\Ga$ is noncompact one obtains many ``interesting'' deformations of $\rho_n$, cf. \cite[Thm. 8.44]{Kapovich00} and \cite{HP}. 
Note, that, unlike in the case of Hitchin representations, the connected component of $\rho_n: \Ga\to SL(n,\C)$ also contains   representations which are not Anosov 
(and some which are not discrete and faithful), since this is already the case for $SL(2,\C)$-representations. 

\medskip 
Weakening the regularity condition to $\tau_{mod}$-regularity and, accordingly, Anosov actions to $\tau_{mod}$-Anosov actions, one obtains more classes, e.g. 
groups of projective transformations acting properly discontinuously cocompactly on bounded strictly  convex solids in the affine space $\R^n$. Such groups are $\tau_{mod}$-Anosov subgroups of $PGL(n+1,\R)$, where $\tau_{mod}$ is the edge of $\sigma_{mod}$ corresponding to the partial flag  {\em line $\subset$ hyperplane}, see \cite[Prop. 6.1]{GW}.

\bigskip 
\section{Discrete subgroups: Domains of proper discontinuity}\label{sec:4}

Note that, so far, we were only looking at limit sets and ignoring domains of discontinuity. The first successful, but limited, treatment of domains of discontinuity was given in \cite{GW}: 
It was proven there that each Anosov subgroup $\Ga< G$ admits a (possibly empty!) 
domain of {proper} discontinuity in a certain bundle $G/AN$ over the full flag manifold $\DF X\cong G/B$. 
These domains were obtained by using a certain embedding of $G$ into a larger Lie group. We will now describe a more comprehensive 
{\em and intrinsic}
treatment of domains of discontinuity, following \cite{coco13,coco15}.
There are two key points in this treatment:

1. Domains of {\em proper} discontinuity are not unique and therefore {\em not canonical} (unlike in the rank 1 case).
There are several natural choices which depend on a certain auxiliary combinatorial datum. 

2. Mumford's GIT (Geometric Invariant Theory) in algebraic geometry serves as a guiding principle.

\subsection{Digression: Mumford's GIT}
 
We begin with the basic {\em topological dynamics} framework of GIT, to be found not in Mumford's book \cite{Mumford}, but e.g.\ in Newstead's lectures \cite{Newstead} and  Dolgachev's book  \cite{Dolgachev}.  

\medskip 
{\bf GIT Mantra:} Let $H$ be a topological group (say, a discrete group or an algebraic group), $Z$ a compact Hausdorff space and $H\times Z\to Z$ a continuous action of $H$ on $Z$. We would like to form a quotient $Z//H$ which is again compact and Hausdorff. In order to do so, we have to partition $Z$ ($H$-invariantly) into {\em semistable} and {\em unstable} points:
$$
Z= Z_{sst} \sqcup Z_{u}
$$
so that $Z_{u}$ is closed. Note that $Z_{sst}$ could be empty. This partition is further refined as follows:

1. $Z_u$ is filtered as an increasing union of closed subsets
$$
Z_0\subset Z_1\subset \ldots \subset Z_{u},
$$
where $Z_0$ is the set of {\em maximally unstable points}. 

2. $Z_{sst}$ contains an open subset $Z_{st}$ of {\em stable} points,
on which the $H$-action is {\em proper}. 
(There is also a subset of {\em nice semistable points}, but we will ignore this.)  

The set of maximally unstable points is, typically, canonical and depends only on the action $H\acts Z$, while the rest of the unstable filtration (including the choice of $Z_{u}$ itself) depends on an auxiliary datum.
In the algebro-geometric context, 
this datum consists of a {\em positive algebraic line bundle} $L\to Z$ (defined up to its tensor power\footnote{In 
the sense that $L$ and $L^{\otimes n}$, $n>0$, lead to the same sets of stable/semistable points.}),
while in our geometric group theory context, it will be a {\em thickening}, see section~\ref{sec:thick} below.

One can think of this partition as: {\em good} (stable), {\em bad} (unstable) and {\em ugly} (semistable but not stable). 

The construction works best when the partition it is {\em neat}\footnote{Note that this is our terminology, it appears that algebraic geometers do not have one.} 
in the sense that {\em stable= semistable}, i.e., ugly = $\emptyset$. 

In order to form the GIT quotient $Z//H$ do the following:

a. Remove the {\em bad} ($Z_{u}$). 

b. Keep the {\em good} ($Z_{st}$) and take the usual topological quotient $Z_{st}/H$ (with the quotient topology), it will be Hausdorff. The set $Z_{st}$ will be {\bf a} domain of 
{proper} discontinuity (in the framework of discrete groups), or {\em domain of properness} in general.  

In the {\em neat} case, you are done. If not:  

c. Deal with the {\em ugly}: For the semistable points use the {\em extended orbit equivalence relation}:
$$
z\sim z'\iff \overline{H z}\cap \overline{H z'}\ne \emptyset
$$
where the closure is taken in $Z_{sst}$. For stable points this  amounts to the usual orbit equivalence: 
$$
H z= H z'. 
$$
Equip the quotient with the quotient topology. Now, if the stars are aligned in your favor, then the resulting quotient is both compact and Hausdorff.

\begin{remark} One can (and should!) vary the auxiliary datum and watch how the quotient space transforms. In the context of symplectic geometry ({\em symplectic reduction}), one sees the variation of  the symplectic structure; generically, one has the {\em neat case} and some degeneration occurs when {\em semistable, non-stable points} appear. 
This is called {\em wall-crossing}.
\end{remark}

What we managed to do in \cite{coco13} is to adapt this mantra to the Morse group actions on various flag manifolds $\Ga\acts Z=G/P$. Note that:

1. \cite{coco13}  dealt only with the regular case, but \cite{coco15} covers the general case of $\taumod$-regular Morse subgroups. 

2. \cite{coco13,coco15} succeeded only in the {\em neat} case: We do not have a theory dealing with the {\em ugly}. 

\medskip
\begin{bexample}
\label{bex:newst}
(cf. Newstead's Example 1.1): Consider the action of $\Ga=\<\gamma\>\cong\Z$ on the real projective plane, which, in an affine patch, is given by:
$$
\ga(x,y)= (\la x, \la^{-1} y), \quad\la > 1. 
$$
The domain of discontinuity is the projective plane minus the three fixed points 
$[1:0:0]$, $[0:1:0]$ and $[0:0:1]$,
which are the only points with infinite stabilizer.
However, the action on this domain is not proper and the quotient is non-Hausdorff.

The maximally unstable set consists only of the two points 
$$[1:0:0], [0:1:0].$$ 
The projective plane minus the $x$- and $y$-axes belongs to the stable part $Z_{st}$ 
(for any choice of a ``line bundle'', or a ``thickening'', in our terminology). 
There, the action is proper.
In order to obtain a larger domain of {\em proper} discontinuity,
one must add a suitable {\em part} of the coordinate axes minus the origin $[0:0:1]$.

Now, we describe three different such enlargements and the corresponding quotients (two of which will be homeomorphic),
resulting from suitable choices of auxiliary data: 
Left, Right (both {\em neat}) and the Center ({\em non-neat}). 

{\bf Left:} Make the entire $x$-axis unstable and include the $y$-axis (minus the origin) into the set of stable points. 
The action of $\Ga$ will be properly discontinuous, cocompact with the quotient $Z_{st}/\Ga\cong T^2$. 

{\bf Right:} Do the same, but make the $y$-axis unstable and add the $x$-axis (minus the origin) to the stable set. 
The quotient is again $T^2$.  
 
Both  left and right partitions of the projective plane will result from what we call {\em balanced thickenings}; 
these will be introduced in Definition~\ref{defn:thickenings} and equation~\eqref{eq:thickening}.  
 
{\bf Center:} Declare the coordinate axes (including the origin) to be {\em semistable but not stable}. 
Then the action on $Z_{sst}$ is not proper (of course!), but the GIT quotient $Z//\Ga$ is compact and Hausdorff: 
It results from $T^2$ by collapsing the union of two parallel essential simple loops to a point. 
\end{bexample}

Note that in Mumford's setting, the degree of unstability is determined by a real (actually, rational) number, the {\em slope}, the value of the Hilbert-Mumford {\em numerical function}. The unstable filtration is given by the standard order on the real numbers: Positive slope means unstable, the more positive the slope is, the more unstable the point is. 
In the KLP setting, the real numbers are replaced with the {\em Weyl group $W$} and its (partial) Bruhat order: 
The smaller the value $w\in W$, the more unstable a point is.
The value $w$ will measure ``how far'' an element of $G/B$ is from the chamber limit set $\La_{ch}(\Ga)$. What corresponds to ``zero'' in this order is not at all clear (and depends on the  choice of a thickening).
Loosely speaking, we have to ``cut $W$ in half'' in order to define an analogue of zero.    

\begin{rem}
This description captures only the $\Ga$-action on $G/B$; in full generality, we also need Bruhat orders on quotients of $W$ by its parabolic subgroups, but we will ignore this here.  
\end{rem}

\subsection{Relative position, Bruhat order and thickenings}
\label{sec:thick}

Given two chambers $\si, \si'\in G/B$ we define their $W$-valued {\em distance} 
or the {\em position of $\si$ relative to $\si'$} 
$$
\delta(\si, \si') = w\in W
$$
as follows. 
Let $\kappa: \amod\to a\subset \tits X$ be a chart\footnote{See Appendix \ref{sec:modsp}.}  
whose image is an apartment $a$ in $\tits X$ 
containing $\si, \si'$ and such that $\kappa$ sends the  model chamber $\simod\subset \amod$ to the chamber $\si'$. 
Then $w\in W=Aut(\amod)$ is defined as the unique element sending $\simod$ to $\kappa^{-1}(\si)$. 
Since transition maps 
between different charts are restrictions of elements of $W$, it follows that 
$w$ is independent of the choice of $\kappa$.\footnote{Note that, by the convexity of apartments,
$a$ must contain the convex hull of $\si'\cup\si$.
Since $\kappa^{-1}$ is determined on the chamber $\si'$,
it follows that it is determined on this convex hull, and in particular also on $\si$.}

The relative position is $G$-invariant: $\delta(g\si, g\si')= \delta(\si, \si')$ for all $g\in G$. 
In general, it is nonsymmetric: 
$$\delta(\si', \si)= \delta(\si, \si')^{-1} .$$
We also define the {\em complementary position}
$$\hbox{c-}\delta(\si,\si'):= w_0 \delta(\si,\si'),$$
where $w_0$ is the longest element of $W$, see \S \ref{sec:Basicgeometry}. In other words, 
if $a\subset \tits X$ is an apartment containing $\si, \si'$ and 
 $\widehat{\si'}\subset a$ is the chamber opposite to $\si'$ then 
$$
\hbox{c-}\delta(\si,\si') = \delta(\si, \widehat{\si'}). 
$$
Since we wish to use $\delta$ as a ``distance'' on $G/B$, we need a (partial) order on $W$ which allows us to compare distances. This order is the {\em Bruhat order}, which we discuss next.  
We first introduce the Bruhat order combinatorially.
(A very detailed discussion of the Bruhat order with many examples can be found in \cite[ch.\ 2]{BB}.)
Afterwards we discuss a geometric way of defining it as the {\em folding order},
which is how we use it in our papers. 

\medskip
{\bf Bruhat order.}  We fix a standard generating system $S$ for $W$
(its elements are called {\em simple reflections} 
and they are the reflections in the faces of the positive fundamental chamber $\simod$),
which defines the {\em word length} $\ell$ on $W$. 
A partial order on $W$, called the {\em (strong) Bruhat order}, is induced by the following convention (and transitivity):

If $v=ur$, where $r$ is a reflection in $W$ (a conjugate of one of the generators) and $\ell(u)< \ell(v)$, then $u<v$. 
In this case one writes 
$$
u \stackrel{r}{\longrightarrow} v. 
$$
In particular, $1$ is the smallest element of $W$ and $w_0$ is the largest. 

Equivalently, this order is given by the condition that $u\le v$ iff a subword (consisting of not necessarily consecutive letters) of a reduced word for $v$ equals a reduced word for $u$.  

We note that left multiplication with $w_0$ {\em reverses} the Bruhat order. 

\begin{example}
Consider $W=S_n$, the permutation group on $n$ letters. As usual, we identify permutations $\pi$ with the 
strings $\pi(1)\ldots \pi(n)$, where we put or not put commas between the adjacent symbols $\pi(i), \pi(i+1)$ when convenient. 
The (standard) simple reflections in $W$ are the transpositions $s_1=(1,2), s_2=(2,3),..., s_{n-1}= (n-1,n)$. 
In the examples below we will always equip $S_n$ with this generating set. The reflections in $W$ are the 
transpositions $(i,j), i < j$. For $r=(i,j)$, the notation $\pi \stackrel{r}{\longrightarrow} \pi'$ means that one moves from $\pi$ to $\pi'$ by transposing $\pi(i), \pi(j)$ in the string $\pi(1)\ldots \pi(n)$, where $\pi(i)<\pi(j)$. In the poset diagrams of $S_3$ and $S_4$ below we connect nodes  
$u$ and $v$ whenever   $u \stackrel{r}{\longrightarrow} v$. 
\end{example}

\begin{figure}[tbh]
\includegraphics[width=90mm]{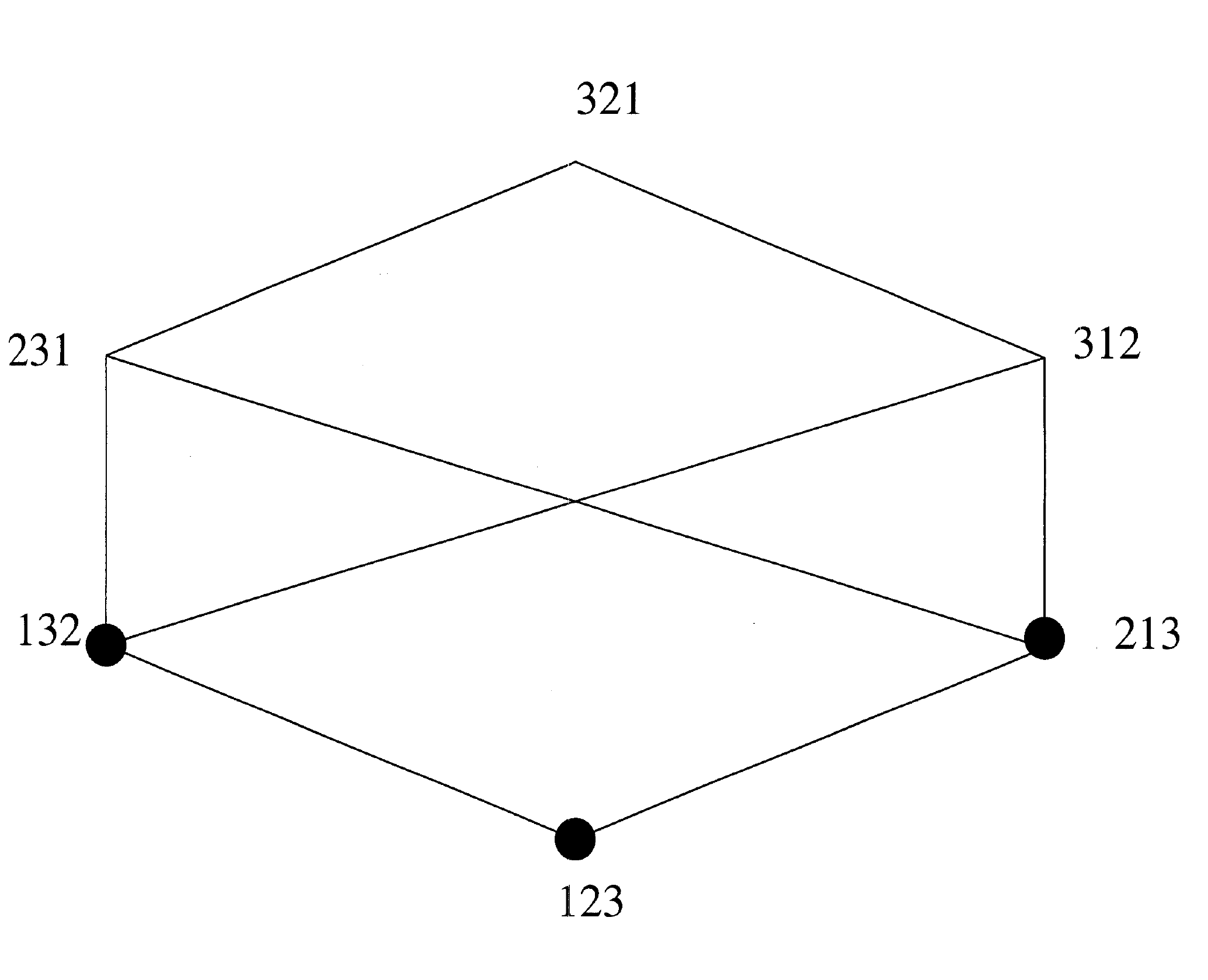}
\caption{The poset diagram of the Bruhat order for $W=S_3$. The larger permutations are higher in the figure, $w_1\ge w_2$ iff the corresponding nodes of the poset diagram are connected by a descending edge path. The circled nodes of the diagram constitute the unique balanced thickening.}
\label{figure13.fig}
\end{figure}

\begin{figure}[tbh]
\includegraphics[width=90mm]{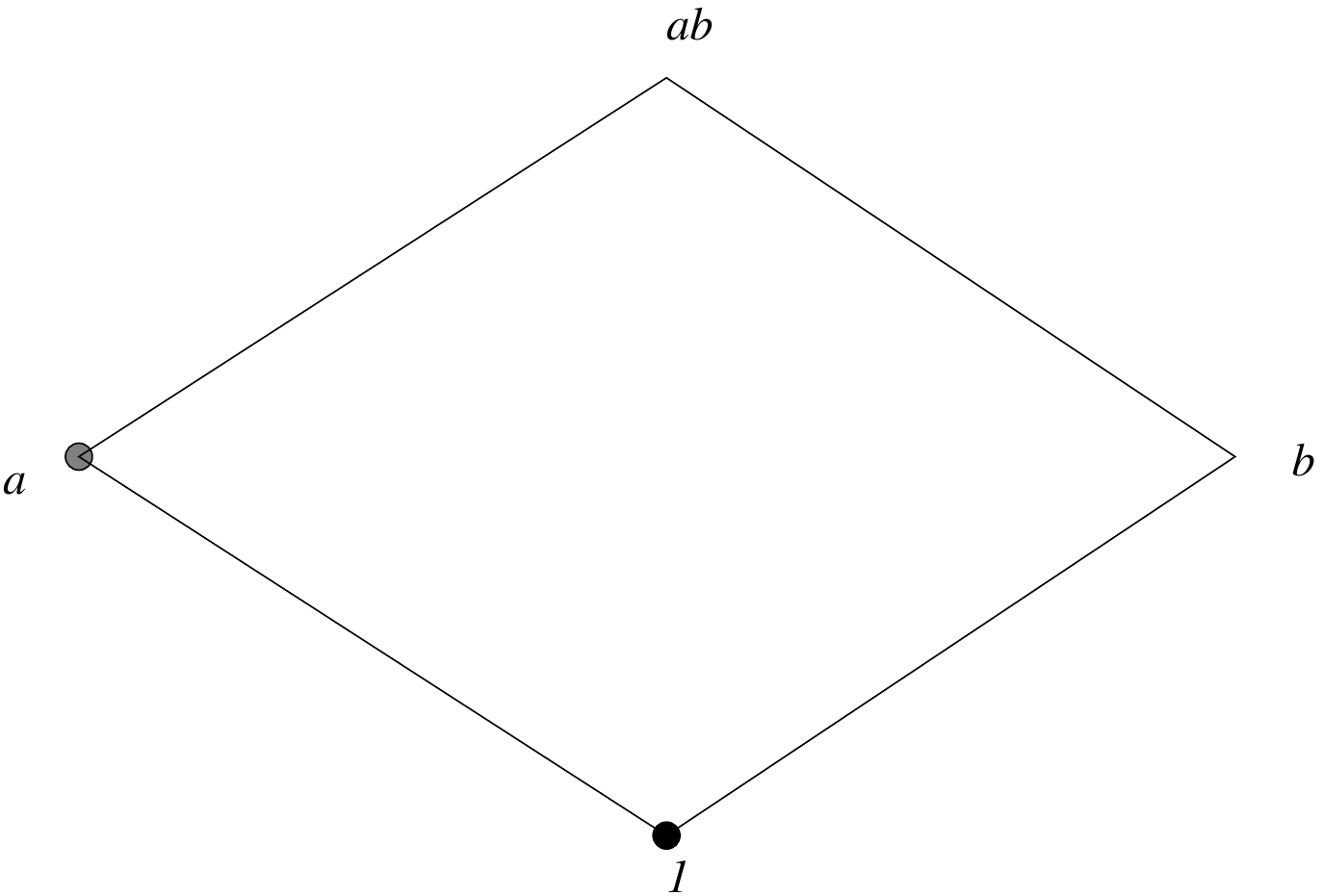}
\caption{The poset diagram of the Bruhat order for $W=\Z_2\times \Z_2$ with the generators $a, b$ and $w_0=ab$. 
The larger permutations are higher in the figure. The circled nodes of the diagram constitute one of the two balanced thickenings.}
\label{A11.fig}
\end{figure}

\begin{figure}[tbh]
\includegraphics[width=90mm]{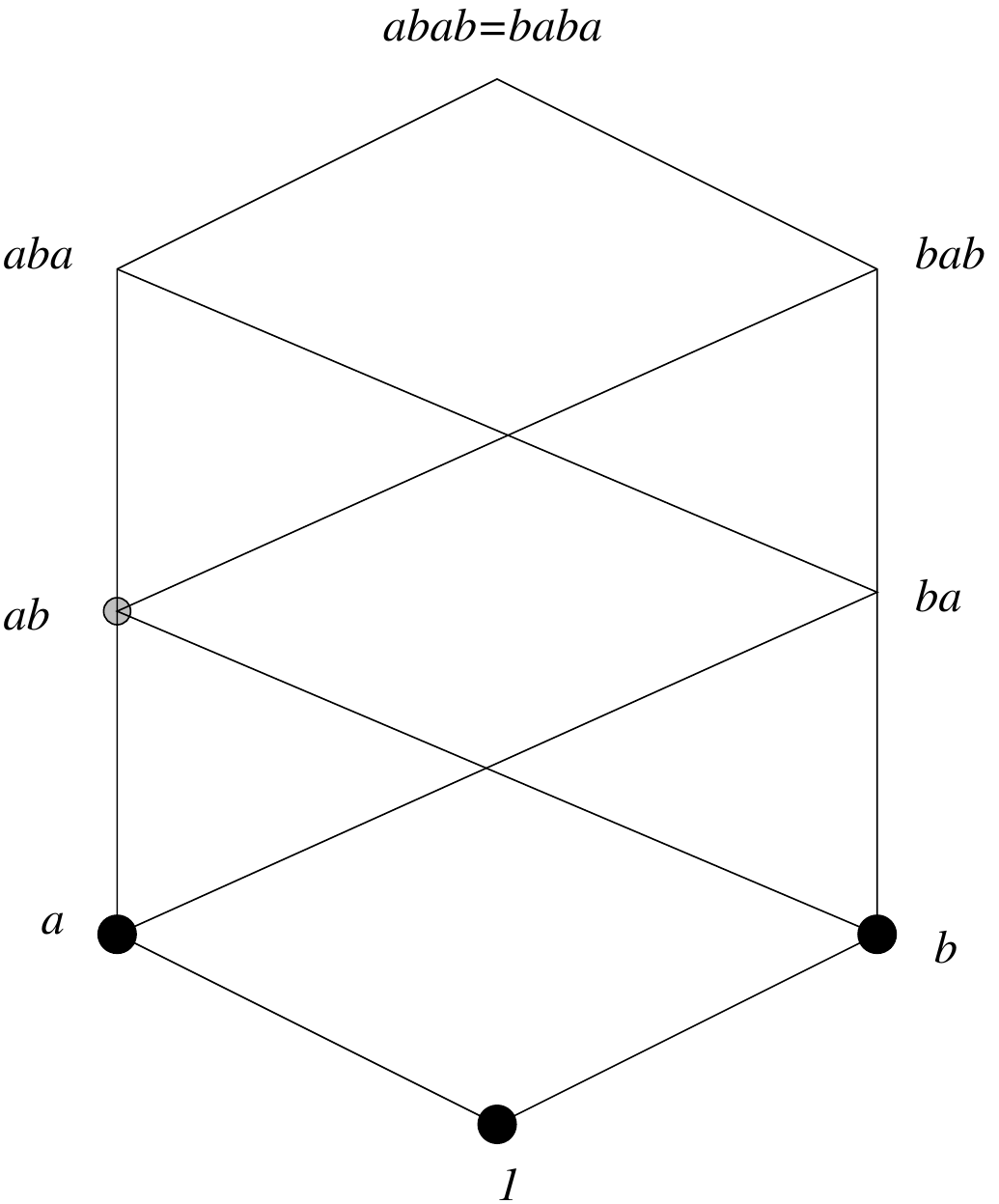}
\caption{The poset diagram of the Bruhat order for $W=B_2$. The larger permutations are higher in the figure. 
The circled nodes of the diagram constitute one of the two balanced thickenings.}
\label{B2.fig}
\end{figure}

\begin{figure}[tbh]
\includegraphics[width=90mm]{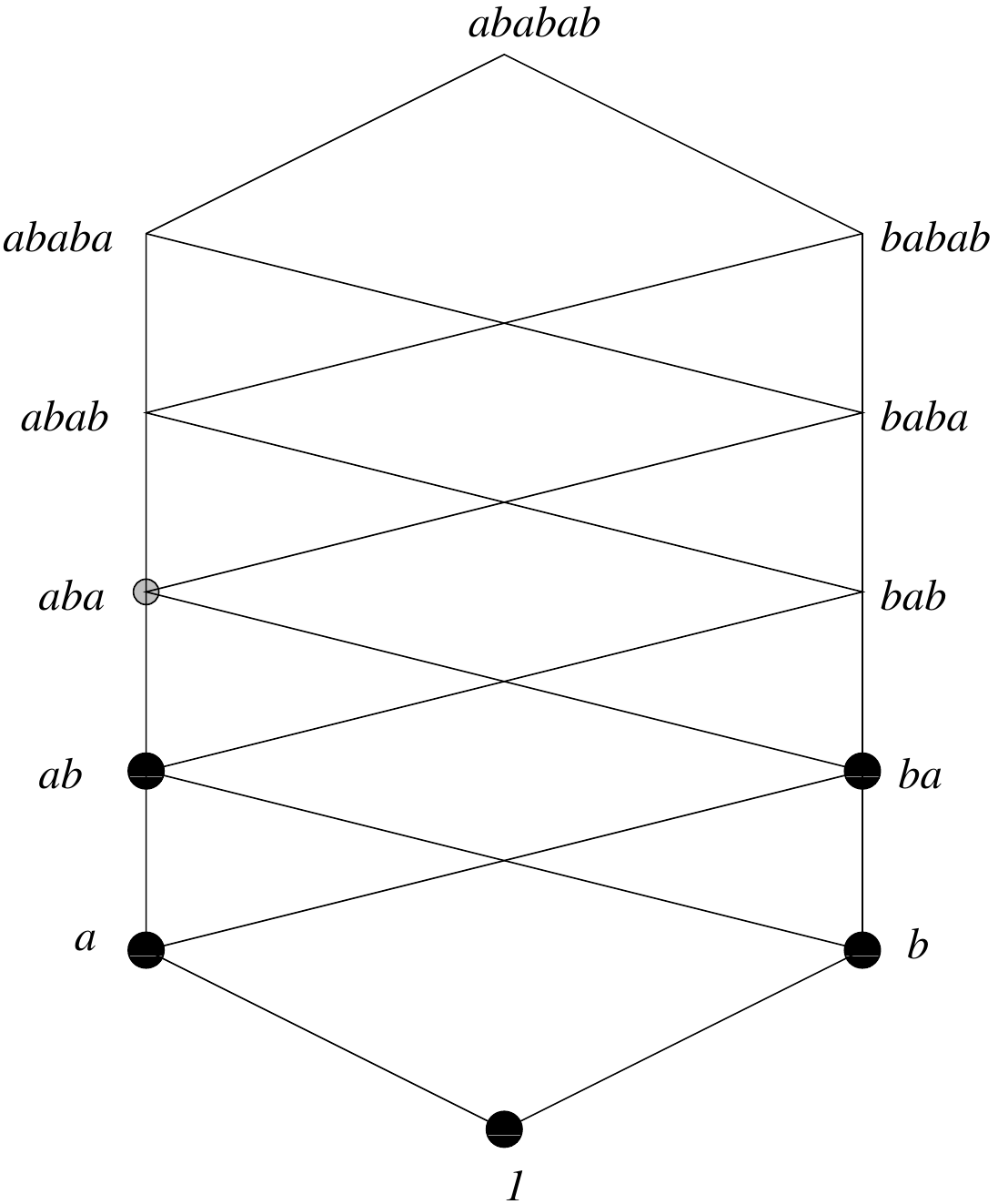}
\caption{The poset diagram of the Bruhat order for $W=G_2$ which is generated by two simple reflections $a, b$. 
The action of $w_0=(ab)^3$ is by the 180 degree rotation of the diagram. The larger permutations are higher in the figure.  The subset $\{1, a, b, ab, ba\}$ is contained 
in every fat thickening. The circled nodes of the diagram constitute one of the two balanced thickenings.}
\label{G2.fig}
\end{figure}

\begin{figure}[tbh]
\centerline{\epsfxsize=4.5in \epsfbox{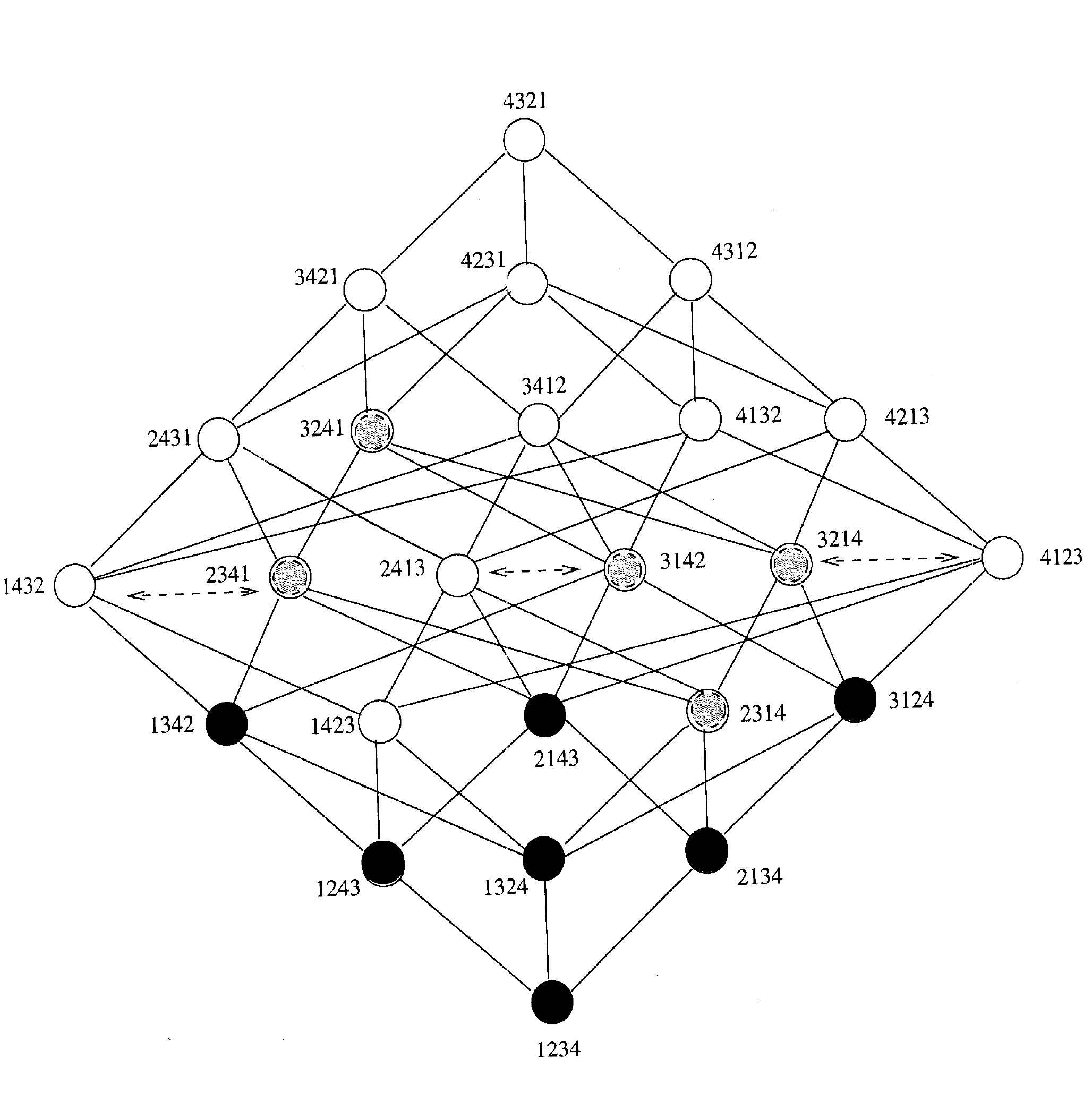}}
\caption{The poset diagram of the Bruhat order for $W=S_4$, see \cite[page 31]{BB}. The larger permutations are higher in the figure. The involution $w_0$ (marked by the 2-sided arrows) acts on this diagram by reversing the order of the labels. For instance, $3241 \stackrel{w_0}{\longleftrightarrow} 1423$. 
Each balanced thickening contains exactly one vertex of each pair $(1432, 2341), (2413, 3142), (3214, 4123)$ since the members of each pair are swapped by $w_0$. 
It follows that all balanced thickenings contain the vertices $1234, 1243, 1324, 2134, 1342, 2143, 3124$ (marked in  solid black). We describe the balanced thickenings by what other vertices they contain. 
(i)  Including in addition both vertices $1423, 2314$ results in a ``metric'' balanced thickening 
(see \cite[\S 4.4]{coco13} or \cite[\S 3.4.1]{coco15} for the definition). 
There are 8 such thickenings. 
(ii)  There are two ``nonmetric'' balanced thickenings, each determined by whether the vertex $3241$ or $4132$ is chosen. For instance, choosing $3241$ forces the thickening to contain the vertices $2341, 3142, 3214$ and $2314$ (these vertices are marked in grey). In total, there are 10 balanced thickenings.}
\label{S4.fig}
\end{figure}

\medskip
{\bf Geometric interpretation of the Bruhat order as the folding order,}
see \cite[\S 4.2+3]{coco13}.
Fix a reference chamber $\simod\subset\amod$.
It will represent the identity element $1\in W$.
For each chamber $\bar\si\subset\amod$ we have a unique $w\in W$ such that $\bar\si=w\simod$. 
Thus, we will identify $W$ with the set of chambers in the model apartment $a_{mod}$.  Given a reflection $s_H\in W$ whose fixed hyperplane (wall) $H$ separates a chamber $\bar\si$ from $\si_{mod}$, we set 
$$
s_H\bar\si< \bar\si.
$$
Now extend this order by transitivity to the entire set of chambers in $a_{mod}$. The result is the Bruhat order.

It is useful to further redefine the Bruhat order non-recursively in terms of ``foldings'' of the model apartment onto itself.
By a {\em folding map} $\amod\to\amod$,
we mean a type preserving continuous map which sends chambers isometrically onto chambers. 
In particular, such maps are 1-Lipschitz.
Intuitively, a folding map fixing the reference chamber $\simod$ moves the other chambers in $\amod$ ``closer'' to $\simod$.

The simplest examples of folding maps fixing $\simod$ are obtained as follows:
A wall $m\subset\amod$ splits $a_{mod}$ into two (simplicial) hemispheres,
the inner hemisphere $h^+$ containing $\si_{mod}$
and the outer hemisphere $h^-$.
This decomposition  gives rise to the folding map
which fixes $h^+$ and reflects $h^-$ onto it.
We call a composition of such folding maps at walls $m_i$
a {\em special folding}.
The above geometric interpretation of the Bruhat order can thus be rephrased:
Two chambers $\bar\si_1,\bar\si_2\subset\amod$ satisfy $\bar\si_1\leq\bar\si_2$ 
iff there exists a special folding moving $\bar\si_2$ to $\bar\si_1$.

In general, not all foldings are special. 
Nevertheless, there is no need of recognizing whether or not a folding is special.
Indeed, one can show that it makes no difference to the order
whether one uses all foldings fixing $\simod$ or only the special ones:

\begin{thm}[{\cite[Cor 4.5]{coco13}}]
For chambers $\bar\si_1,\bar\si_2\subset\amod$
it holds that $\bar\si_1\leq\bar\si_2$ 
iff there exists a folding map $\amod\to\amod$ fixing $\simod$ and sending $\bar\si_2\mapsto\bar\si_1$.
\end{thm}

\medskip
{\bf Thickenings.} We will use special subsets $\Th\subset W$, called ``thickenings of $1\in W$ inside $W$''. 
They are defined as unions of sublevels of the Bruhat order,
i.e.\ they satisfy the property:
$$ v\in\Th\;\hbox{~~and~~}\; u<v\Ra u\in\Th .$$
In particular, $1\in \Th$ for every {nonempty thickening} $\Th$. One can think of {thickenings} as being starlike with respect to $1\in W$ (and the Bruhat order defining intervals). Simple examples of thickenings are given by the ``closed balls''
$$
B(1, r)=\{w\in W: w\le r\}, 
$$
where one can think of $r\in W$ as the ``radius'' of the ball. General thickenings are unions of such ``balls''. 

\begin{rem}
In the theory of posets, thickenings are called {\em (lower) ideals} and thickenings of the form $\{w\le r\}$ are called {\em principal ideals}. 
\end{rem}

A thickening is {\em proper} if it is nonempty and is different from the entire $W$. 
The latter condition is equivalent to the requirement that the longest element 
$w_0$ is not in the thickening.

\begin{definition}\label{defn:thickenings}
1. A thickening $\Th\subset W$ is {\em slim} if $\Th\cap  w_0 \Th=\emptyset$.

2. A thickening $\Th\subset W$ is {\em fat} if $\Th\cup w_0 \Th=W$. 

3. A thickening $\Th\subset W$ is {\em balanced} if it is both slim and fat. 
\end{definition} 

Thus, for each $w\in W$,
a slim (fat, balanced) thickening contains at most (at least, exactly) one 
of the pair of complementary elements $w$ and $w_0w$.
In particular, a balanced thickening consists of precisely half of the elements of $W$.

\begin{theorem}
\cite{coco13, coco15}. Each finite Weyl group $W$ admits at least one balanced thickening.
\end{theorem}

\medskip 
{\bf Examples.} Among the rank 2 Weyl groups, $A_2$ admits exactly one balanced thickening,
whereas $B_2, G_2$ and $A_1\times A_1$ admit exactly two!  

\begin{example}
Consider $W= \Z_2 \times \Z_2$,  the Coxeter group of the type $A_1\times A_1$, see Figure \ref{A11.fig}. 
We let $a, b$ denote the generators of the direct factors of $W$; these are the simple reflections and the only reflections in $W$ (since $W$ is abelian). The poset graph of $W$ is completely described by the inequalities $w_0=ab > a > 1$ and $w_0> b>1$. 
The action of  the involution $w_0$ swaps the nodes $a, b$ as well as the nodes $1, w_0$. Therefore, the thickenings $B(1, a)$, $B(1, b)$ are balanced. The only other two proper thickenings are $\{1\}$ and $W - \{w_0\}$; these are respectively slim and fat. In particular, the group $W= \Z_2 \times \Z_2$ has exactly two balanced thickenings. 
\end{example}

\begin{example}
Consider $W=S_3$, the Coxeter group of the type $A_2$. We refer the reader to 
Figure \ref{figure13.fig} for the description of the poset graph of $W$. 
We will describe all proper thickenings in $W$.  Every nonempty thickening contains $1\in W$. Each thickening different from $\{1\}$ also contains at least one of the transpositions $(12)$ and $(23)$. 
The slim thickenings consisting of two elements are $B(1, (123))$ and $B(1, (231))$. 
Since the involution $w_0$ swaps the upper and lower halves of the poset diagram, the thickening $I:= \{(123), (132), (213)\}$ is balanced. 
A thickening containing  either $(231)$ or $(312)$ 
also contains $I$.
Hence $I$ is the only balanced thickening. 
The fat thickenings consisting of four elements are $B(1, (231))$ and   $B(1, (312))$.
\end{example}

\begin{example}\label{ex:B2}
Consider the Coxeter group $W$ of the type $B_2$. We refer the reader to 
Figure \ref{B2.fig} for the poset graph of $W$. The action of $w_0=(ab)^2$ swaps $ab$ and $ba$, $a$ and $bab$, $b$ and $aba$. Since balanced thickenings consist of exactly four elements, they can contain neither $aba$ nor $bab$. 
Also, a balanced thickening has to contain either $ab$ or $ba$ but not both. From this, we conclude that 
the only two balanced thickenings in $W$ are $B(1, ab)$ and $B(1, ba)$. 
\end{example}

\begin{example}
Consider the Coxeter group $W$ of the type $G_2$, see  
Figure \ref{G2.fig} for the poset graph of $W$. The only two balanced thickenings in $W$ are $B(1, aba)$ and 
$B(1, bab)$. 
\end{example}

\begin{rem}
A  subset $R\subset W$ determines the thickening  
$$
\Th_R:=\bigcup_{r\in R} B(1, r). 
$$
Every thickening has this form. For instance, in Example \ref{ex:B2} we can take $R=\{ab, ba\}$ and hence obtain the fat unbalanced thickening 
$$
\Th_R= \{1, a, b, ab, ba\}. 
$$ 
In order to simplify the notation we omit the symbol $R$ in the notation $\Th_R$ from now on: $\Th$ will always denote {\em a certain thickening}. 
\end{rem}

\begin{rem}
In { our work \cite{morse, coco15, bordif}}, we need the folding order (and thickenings) 
more generally also for 2-sided quotients $W_P\backslash W/W_Q$, where $W_P$ and $W_Q$ are Coxeter subgroups generated by subsets of the set of simple reflections.
However, {for the sake of simplicity}, we will not discuss these here. 
\end{rem}

\begin{question}
Is there a reasonable combinatorial classification of balanced thickenings for a given finite Coxeter group 
$W$? What are the asymptotics of the numbers of balanced thickenings in the Weyl groups of the types
$$A_n, \quad B_n, \quad D_n, \quad \underbrace{A_1\times \ldots \times A_1}_{n\; \hbox{\scriptsize times}}$$ 
as $n\to\infty$? 
\end{question}

\medskip
We now turn to discussing flag manifolds. For each chamber $\si\in\DF X\cong G/B$ and $r\in W$ we define the ``combinatorial $\delta$-sphere'' of  the ``combinatorial radius'' $r$,  
$$
S(\si,r)=\{\si'\in G/B: \delta(\si', \si)= r\},
$$
also known as a {\em Schubert cell} in $G/B$, and the ``combinatorial $\delta$-ball'' 
$$
B(\si, r)=\{\si'\in G/B: \delta(\si', \si) \le r\}, 
$$
{also known as a {\em Schubert cycle}. The following is a basic fact of the theory of Lie groups that 
plays a critical role in our analysis of discontinuity domains.
It expresses that the Bruhat order on $W$ corresponds to the inclusion order on Schubert cycles in $G/B$
(with respect to a fixed reference chamber $\si$, respectively, minimal parabolic subgroup $B$):

\begin{theorem}\label{thm:basic}
The distance $\delta$ is lower  semicontinuous with respect to the manifold topology.
Moreover, 
\begin{equation}
\overline{S(\si, r)}= B(\si, r),
\end{equation}
where the closure is taken in the manifold topology of $G/B$. \end{theorem}

Consequently, $S(\si, r')\subseteq\overline{S(\si, r)}$ iff $r'\leq r$,
and otherwise $S(\si, r')\cap\overline{S(\si, r)}=\emptyset$.

In the case of complex Lie groups, this theorem goes back to work of Chevalley in the 1950s \cite{Chevalley}, see 
also \cite{BB}; for the proofs in the general case (including reductive groups over local fields), see \cite{Borel-Tits} as well as \cite{Mitchell, Mitchell2008}. The most general case dealing with subsets of partial flag manifolds is established in \cite{coco15}.}

\medskip
We next use the thickenings of the neutral element inside $W$ 
to produce corresponding thickenings of (sets of) chambers inside $\DF X\cong G/B$.

Let $\Th\subset W$ be a thickening. Given a chamber $\si\in \DF X$, define its {\em thickening} inside $\DF X$ by 
\begin{equation}\label{eq:FTH}
\Th(\si)= \{\si'\in \DF X: \delta(\si', \si)\in \Th\}. 
\end{equation}
For a subset $\La\subset \DF X$ we define its {\em $\Th$-neighborhood} or  {\em thickening} as
\begin{equation}\label{eq:thickening}
\Th(\La):= \bigcup_{\la\in \La} \Th(\la)= \{\si'\in \DF X : \exists \la\in \La, \delta(\si', \la)\in \Th\}. 
\end{equation}
It is clear ({from the $G$-invariance of $\delta$}) that thickenings are $G$-invariant:
$$
\Th(g \La)= g \Th(\La), g\in G.  
$$
Thus, if $\Ga<G$ is a subgroup preserving $\La$, it also preserves $\Th(\La)$. 

Our motivation for introducing the notion of slimness is the observation,
that the slimness of a thickening in $W$
is equivalent to the disjointness of the corresponding thickenings 
of any two antipodal chambers in $\DF X$:
\begin{lem}
[{\cite{coco15}}]
\label{lem:sldisj}
Let $\Th\subset W$ be a slim thickening.
Then for any two antipodal chambers $\si,\hat\si\in\DF X$
it holds that 
$$\Th(\si)\cap\Th(\hat\si)=\emptyset .$$
\end{lem}
\proof
This follows from the definition of slimness and the triangle inequality\footnote{The inequality
can be regarded as a triangle inequality in $G/B$ for the $W$-valued combinatorial side lengths
of the triangle with vertices $\si,\si'$ and $\hat\si$.}
$$ \delta(\si',\hat\si)\geq \hbox{c-}\delta(\si',\si)$$ 
for chambers $\si'\in\DF X$.

Indeed, suppose that $\si'\in\Th(\si)\cap\Th(\hat\si)$.
Then $\delta(\si',\hat\si), \delta(\si',\si)\in\Th$.
Due to the inequality,
also $\hbox{c-}\delta(\si',\si)\in\Th$,
equivalently, $\delta(\si',\si)\in w_0\Th$. 
It follows that $\Th\cap w_0\Th\neq\emptyset$, contradicting slimness.

To verify the inequality,
consider the apartment $a\subset\tits X$ containing $\si,\hat\si$
and a folding retraction $r:\tits X\to a$,
i.e.\ a type preserving continuous map which fixes $a$ pointwise. 
Such a retraction is given e.g.\ by the natural projection
$\tits X\to\tits X/B_{\si}\cong a$ 
where $B_{\si}$ denotes the minimal parabolic subgroup fixing $\si$.
Then $ \delta(\si',\hat\si) \geq \delta(r\si',\hat\si) = \hbox{c-}\delta(r\si',\si)=c-\delta(\si',\si)$. 
\qed

\medskip
The importance of slimness comes therefore from the following fact.
Suppose $\Th\subset W$ is a slim thickening and $\La\subset \DF X$ is an {\em antipodal subset}, 
i.e.\ a subset where any two distinct elements are antipodal. 
Then for every $\si\in \Th(\La)$ there exists a unique $\la=\la_\si\in \La$ such that $\delta(\si, \la)\in\Th$. 
Thus, we obtain a natural projection 
$$
\pi: \si\mapsto \la_\si, \Th(\La)\to \La. 
$$
As an exercise, let us prove continuity (in the subspace topology induced from $\DF X$ ) 
of this projection, provided that $\La$ is closed: 
Let $\si_i\to \si$ in $\Th(\La)$, where $\si_i\in \Th(\la_{\si_i}), \si\in \Th(\la_{\si})$.
After extraction, $\la_{\si_i}\to \la\in \La$, since $\La$ is compact. 
By the semicontinuity of 
$\delta$ and the fact that $\Th$ is a thickening, we obtain that 
$\delta(\si, \la)\in \Th$ and hence $\la=\la_{\si}$,
establishing continuity.
One verifies further that the projection $\pi$ is a fiber bundle over $\La$ with fibers homeomorphic to $\Th(\la)$, $\la\in \La$ (see \cite{coco15}). 

\medskip
The fatness of a thickening, in turn, does not imply that the union 
$\Th(\si)\cup \Th(\hat\si)$ is the entire $\DF X$, but it does imply that $\Th(\si)\cup \Th(\hat\si)$ covers the entire chamber set of the unique apartment $a\subset \tits X$ containing $\si, \hat\si$.  
The importance of the notion of fatness is less immediate. 
The proof of Theorem~\ref{thm:proper} below shows why 
it is useful for the proof of proper discontinuity of discrete group actions on certain domains in flag manifolds.

\medskip 
Here is how one can think of thickenings $\Th(\la)\subset \DF X$ of points in the Furstenberg boundary.
First of all, if $\Th=\Th_r, r\in W$, then one can think of 
$\Th_r(\la)$ as the ``combinatorial $r$-neighborhood'' of $\la$ in $ \DF X$, 
as it consists of all $\si\in \DF X$ which are within $\delta$-distance $\le r$ from $\la$. 
However, caution is needed here, since for any proper thickening $\Th\subset W$, 
the corresponding thickening $\Th(\la)\subset \DF X$ is nowhere dense in the visual topology 
of   $\DF X$. Therefore, a better way to think of thickenings of subsets in $\DF X$ is as follows. 
The choice of $\Th$ describes the ``degree of nongenericity'' of the relative position of chambers $\si\in  \DF X$ with respect to $\la$. For instance, for $\Th=\Th_r$ ($r\in W$), 
the larger the element $r$, the more generic the relative position we allow for points in 
$B(\la, r)\subset \DF X$. The most generic relative position is achieved for points 
in the open Schubert cell $S(\la, w_0)=\la^{opp}$, consisting of all chambers $\si$ in $\DF X$  opposite of 
$\la$. The closure of this open Schubert cell is the entire Furstenberg boundary, the metric ball $B(\la, w_0)$. 
One implication of Theorem \ref{thm:basic} is that taking limits of sequences 
$\si_i$ as $i\to \infty$, can only result in the decrease of genericity (the limit of a sequence 
can only be ``more special'' with respect to $\la$, not ``less special'').

\medskip
The next lemma implies the important fact that thickenings of compact subsets are compact.

\begin{lemma}
\label{lem:closth}
For every thickening $\Th\subset W$ and every subset $\La\subset \DF X$, we have 
$$
\overline{\Th(\La)}= \Th(\overline{\La}). 
$$
\end{lemma}
\proof Suppose that $(\la_i)$ is a sequence in $\La$ and $\si_i\in \Th(\la_i)$
such that $\si_i\to\si$ in $\DF X$.
After extraction, $\la_i\to\la\in \overline{\La}$. 
In view of semicontinuity (Theorem \ref{thm:basic}), 
we have 
$\delta(\si, \la)\leq\delta(\si_i, \la_i)$
for large $i$. 
Since $\delta(\si_i, \la_i)\in \Th$,
it follows that $\delta(\si, \la)\in\Th$. 
This shows that $\overline{\Th(\La)}\subset \Th(\overline{\La})$.

Conversely, consider a sequence $\la_i\in \La$ converging to $\la\in\overline{\La}$ and let $\si\in \Th(\la)$. 
Then, since $\DF X= G/B$, there exists a sequence $g_i\to1$ in $G$ such that $\la_i=g_i\la$. 
Then $\si$ is the limit of the $\si_i:= g_i\si\in\Th(\la_i)\subset \Th(\La)$ 
and, therefore, belongs to $\overline{\Th(\La)}$. \qed

\begin{example}
\label{ex:ellth}
Consider $G=SL(3,\R)$ and the unique balanced 
thickening $$\Th=B(1, (12)) \cup B(1, (23)).$$ 
Let $\la$ be a chamber 
in the Tits building $\tits X$ of the symmetric space $SL(3,\R)/SO(3)$, which is the incidence graph of the real projective plane. 
We will think of $\la$ as a flag $(p,l)$ (where $p$ is a point and $l$ is a line in projective plane). 
The  thickening $\Th(\la)$ consists of all chambers in $\tits X$ which share a vertex with $\la$. In other words,  
the condition that a flag $(p',l')$ belongs to $\Th(\la)$ means that either 
$p'\in l$ or $p\in l'$.
Clearly,  $\Th(\la)$ is a closed subset of the flag manifold, homeomorphic to the wedge of two projective lines. It equals the closure of the union of two 
``combinatorial spheres'' $S(\la, (12))$ and $S(\la, (23))$; this union consists of chambers sharing a vertex with $\la$ but different from $\la$.   

Next, take an  ellipse $E\subset \R^2\subset \R P^2$. The (projectivized) tangent bundle $PE$ 
of $E$ defines a lift $\tilde{E}$ of $E$ to the flag manifold  $\Flag(\R^3)$, the full flag manifold of $\R^3$. 
It consists of the {\em tangent flags} $(p,l)$, $p\in E$, $l$ is the tangent line to $E$ at $p$; clearly, $\tilde E$ is homeomorphic to $E$.
Now, let $\Th$ be the unique balanced thickening in the group $W=S_3$. 
Then the corresponding thickening $\Th(\tilde E)$ 
consists of the flags $(q, m)$, where either $q\in E$ (and $m$ is any line through $q$) or $m$ is a line tangent to $E$ and $q$ is any point on $m$. 
Topologically speaking, $\Th(\tilde E)$ is the trivial bundle over $\tilde E$ with fibers homeomorphic to $S^1\vee S^1$,
and $\tilde E$ is (the image) of its distinguished section
with values in the singular points of the fibers.
The projection $\pi: \Th(\tilde E)\to \tilde E$ sends 
each flag $(p,m)$, $p\in E$, to the corresponding tangent flag $(p,l)$; 
and it sends each flag $(q,l)$, $l$ is tangent to $E$, to the corresponding tangent flag $(p,l)$.
\end{example}

\begin{figure}[tbh]
\includegraphics[width=90mm]{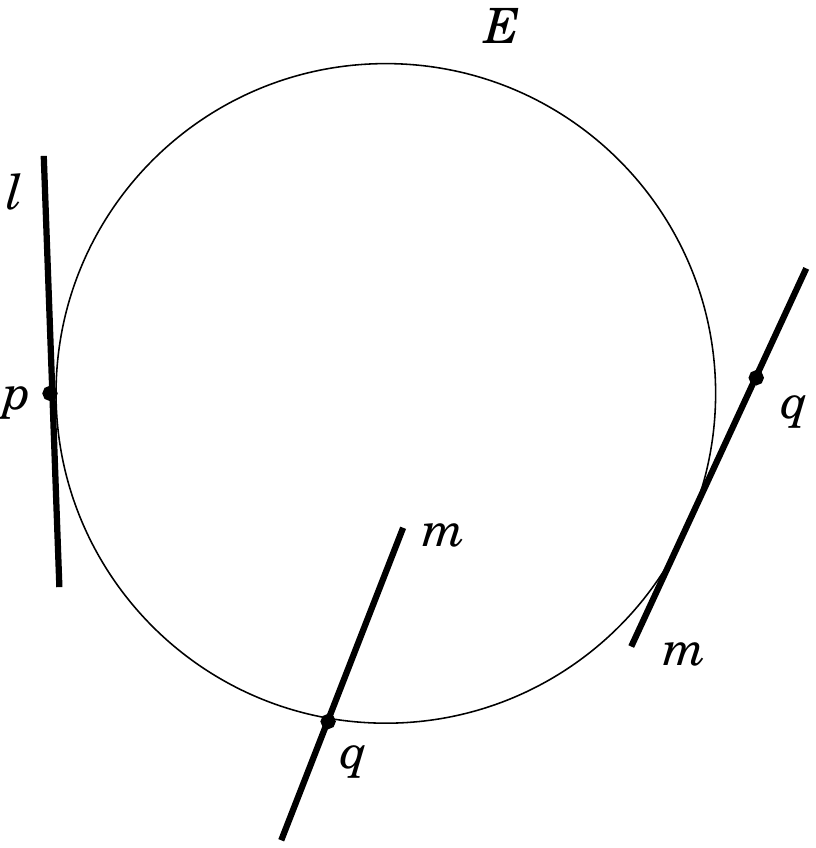}
\caption{The balanced thickening of $\tilde E$}
\label{figure14.fig}
\end{figure}

\subsection{Domains of proper discontinuity, cocompactness and nonemptiness}
\label{sec:dd}

We start this section
by reviewing some basic notions from topological dynamics.
We consider topological actions $\Ga\acts Z$ of discrete groups on metrizable locally compact topological spaces. 

Recall that the $\Ga$-action on an invariant open subset $\Om\subset Z$
is {\em properly discontinuous} if for each compact $K\subset \Om$ it holds that 
$$\ga K\cap K=\emptyset$$
for all but finitely many $\ga\in\Ga$. 
A weaker condition is the {\em discontinuity} of the action. A point $z\in Z$ is said to be {\em wandering} for the $\Ga$-action
if there exists a neighborhood $U$ of $z$ such that $$\ga U\cap U=\emptyset$$  for all but finitely many $\ga\in \Ga$. 
An action is called {\em discontinuous} if each point is wandering. 
The {\em domain of discontinuity} $\Om_{disc}\subset Z$ for the action $\Ga\acts Z$ is the set of wandering points. 
This set is clearly open and invariant; in general, however, the action on the domain of discontinuity  is {\em not proper}. 

\begin{example}\label{ex:cyclic}
(Compare Example~\ref{bex:newst})
$\ga\in SL(3,\R)$, $\ga=Diag(\la, 1, \la^{-1})$, $\la>1$. 
Then $\ga$ has on $\R P^2$ the three fixed points $e_1=[1:0:0], e_2=[0:1:0], e_3=[0:0:1]$. 
The point $e_1$ is attractive, $e_3$ is repulsive and $e_2$ is hyperbolic for the action of $\ga$ on $\R P^2$. Denoting $L_{ij}$ the projective line through $e_i, e_j$, $i<j$, we obtain that the domain of discontinuity for the action of $\Ga=\<\ga\>$   on $\R P^2$ is the complement to its fixed point set, $\{e_1, e_2, e_3\}$.  
However, the action of $\Ga$ on this domain is not proper. In order to get a maximal domain of proper discontinuity, 
one removes from $\R P^2-\{e_1, e_2, e_3\}$ either the entire line $L_{12}$ or $L_{23}$. 
\end{example}

\begin{figure}[tbh]
\includegraphics[width=90mm]{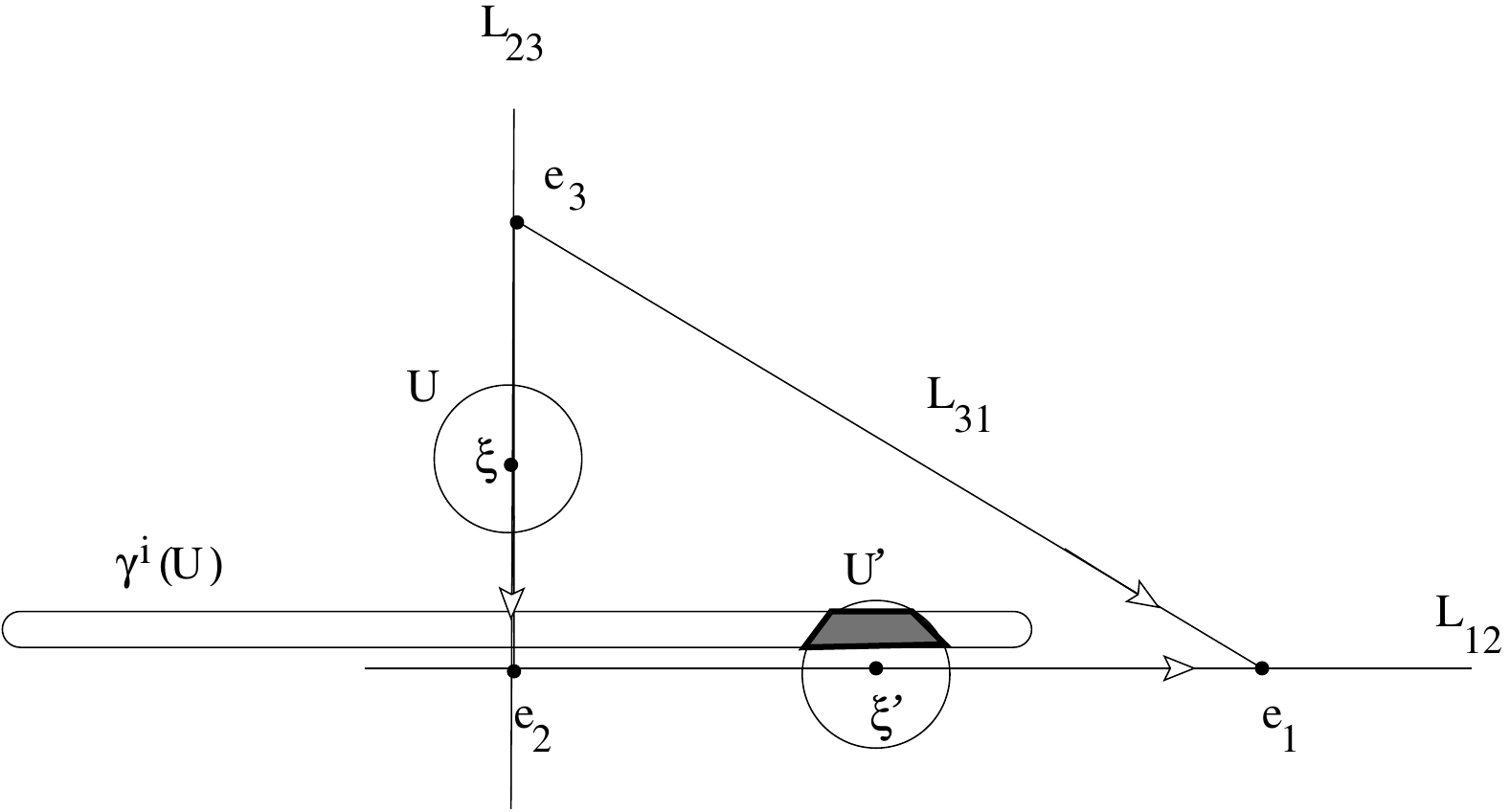}
\caption{Dynamics of a cyclic subgroup $\Ga=\<\ga\>$ on the projective plane. The points $\xi, \xi'$ are dynamically related.}
\label{figure15.fig}
\end{figure}

\medskip 
Our arguments for proving proper discontinuity will be based on the fact that 
it is equivalent to the {\em absence of dynamical relations} between points 
relative to the $\Ga$-action.

\begin{defn}
Two points $\xi, \xi'\in Z$ are {\em $\Ga$-dynamically related},
$\xi\stackrel{\Ga}{\sim}\xi'$, 
 if for each pair of neighborhoods $U, U'$ of $\xi, \xi'$ respectively, there are infinitely many elements $\ga\in \Ga$ such that
$$
\ga U\cap U'\ne \emptyset. 
$$ 
\end{defn}

Note that a point is wandering iff it is not dynamically related to itself.

It is straightforward that, since the space $Z$ is Hausdorff and 1st countable, 
$\Ga$-dynamical relation $\xi\stackrel{\Ga}{\sim}\xi'$ can be reformulated as follows: 
There exists a sequence of distinct elements $\ga_n\in \Ga$ and a sequence $\xi_n\to \xi$ in $Z$ such that $\ga_n\xi_n\to \xi'$. We will write then more precisely
$$\xi\stackrel{(\ga_n)}{\sim}\xi'.$$

\begin{lemma}
The $\Ga$-action on an open invariant subset $\Om\subset Z$
is properly discontinuous iff no two points of $\Om$ are $\Ga$-dynamically related. 
\end{lemma}
\proof 
Suppose first that the action $\Ga\acts\Om$ is not properly discontinuous. 
Then there exists a compact $C\subset \Om$ 
and elements $\ga_n\to\infty$ in $\Ga$ such that $\ga_n C\cap C\ne \emptyset$.
Hence there are points $x_n\in C$ such that also $\ga_nx_n\in C$.
By compactness,
after extraction, we have convergence $x_n\to x\in C$ and $\ga_nx_n\to x'\in C$
and it follows that $x,x'$ are $\Ga$-dynamically related.

Conversely, suppose that the points $x, x'\in \Om$ are dynamically related. 
Taking $U, U'$ to be relatively compact neighborhoods of $x, x'$ respectively, we obtain that 
$\ga U\cap U'\ne \emptyset$
for infinitely many $\ga\in \Ga$. 
For these $\ga$ and the compact $C=\ol{U}\cup \ol{U'}$ it holds that $\ga C\cap C\ne \emptyset$.
\qed 

\medskip
We now take up the discussion of the topological dynamics of 
discrete group actions on flag manifolds.
We will restrict to the case of actions on the full flag manifold $\DF X\cong G/B$.

Regarding proper discontinuity,
the connection between dynamical relation
and the Bruhat order and thickenings
comes in the form of the following key lemma 
which provides a relative position inequality for dynamically related points.
Roughly speaking,
it says that they cannot be both far, in the sense of the combinatorial distance $\delta$, from the chamber limit set $\La_{ch}(\Ga)$.
The lemma, in turn, is derived from the {\em higher rank convergence property} for the action $\Ga\acts G/B$ 
discussed in section~\ref{sec:conprop}. 

\begin{klemma}[\cite{coco13, coco15}]
Suppose that $\Ga< G$ is regular and $\xi, \xi'\in G/B$ are $\Ga$-dynamically related. 
Then there exist (not necessarily distinct) limit chambers $\la, \la'\in \La_{ch}(\Ga)$ such that
$$
\delta(\xi',\la') \le {\hbox{c-}}\delta(\xi,\la). 
$$
\end{klemma}
\proof  Suppose first that we have a dynamical relation $\xi\stackrel{(\ga_n)}{\sim}\xi'$, $\ga_n\in \Ga$.
Then, by the definition of dynamical relation, there exists a sequence $(\xi_n)$ in $G/B$ such that 
$\xi_n\to\xi$ and $\ga_n\xi_n\to\xi'$. 
The regularity of the subgroup $\Ga<G$ translates via Theorem~\ref{thm:regconv}
into the higher rank convergence property for the action $\Ga\acts G/B$.
Hence, 
after extraction, there exists a pair of limit chambers $\si_\pm\in\La_{ch}(\Ga)$ such that 
$\ga_n$ converges to $\si_+$ uniformly on compacts in the open Schubert cell $\si_-^{opp}$. 
Let $a\subset\geo X$ be an apartment containing $\si_-$ and $\xi$. 
Nearby apartments $a_n$ containing $\xi_n$ 
can be obtained by using small isometries $h_n\to 1$ in $G$, 
with $\xi_n=h_n\xi$ and putting $a_n=h_na$.
Let $\hat\si_-\subset a$ be the chamber opposite to $\si_-$,
and let $\si_n=h_n\hat\si_-\subset a_n$.
Then $\si_n\to\hat\si_-$.
Since $\hat\si_-\in \si_-^{opp}$,
the locally uniform convergence of $\ga_n$ to $\si_+$ 
implies that $\ga_n\si_n\to\si_+$. 
We obtain
\begin{equation*}
\delta(\xi',\si_+)\le \delta(\ga_n\xi_n,\ga_n\si_n)
=\delta(\xi_n,\si_n)
=\delta(h_n\xi,h_n\hat\si_-)
=\delta(\xi,\hat\si_-)={\hbox{c-}}\delta(\xi,\si_-) 
\end{equation*}
for large $n$, 
where the first inequality follows from 
the semicontinuity of $\delta$, see Theorem \ref{thm:basic}. 
Putting $\la=\si_-$ and $\la'=\si_+$ yields the assertion. 
\qed 

\medskip
As a consequence, 
no dynamical relations occur in domains which are far enough from the chamber limit set
in the combinatorial sense of relative position,
i.e.\ which avoid a sufficiently large thickening of it. 
(Recall that this means that the points in these domains have sufficiently generic position with respect to all limit chambers.)
We obtain:

\begin{thm}
[Proper discontinuity \cite{coco13}]
\label{thm:proper} 
Suppose that $\Ga< G$ is regular and $\Th\subset W$ is a fat thickening. 
Then no two points in the domain
$$\Om_{Th}(\Ga):=G/B-\Th(\La_{ch}(\Ga))$$ are $\Ga$-dynamically related.\footnote{We recall that 
$\Th(\La_{ch}(\Ga))$ is compact, because $\La_{ch}(\Ga)$ is, cf.\ Lemma~\ref{lem:closth}.}
In other words, the action $\Ga\acts\Om_{Th}(\Ga)$ is properly discontinuous. 
\end{thm}
\proof Suppose that $\xi, \xi'\in G/B$ are dynamically related. 
Then, by the lemma, there exist $\la, \la'\in\La_{ch}(\Ga)$ such that
$$
\delta(\xi',\la') \le \hbox{c-}\delta(\xi,\la). 
$$ 
By the definition of fat thickening, for the relative position $w:=\delta(\xi,\la)$ either 
\begin{enumerate}
\item $w\in \Th$, or
\item $w_0 w\in \Th$.
\end{enumerate}
In the former case, $\xi\in \Th(\La_{ch})$.
In the latter case, $\hbox{c-}\delta(\xi,\la)\in \Th$, which, 
{\em by the definition of a thickening}, implies that $\delta(\xi',\la')\in \Th$,
and thus $\xi'\in \Th(\La_{ch})$.
Hence, $\xi, \xi'$ cannot be both in $\Om_{Th}(\Ga)$. \qed  

\medskip
Note that we are not assuming here that the chamber limit set $\La_{ch}$ is antipodal. 
Antipodality is used, in conjunction with {\em slimness} of $\Th$ and the {\em expansion} axiom, to ensure the {\em cocompactness} of $\Ga$-actions. 
It is an important fact that, 
for a slim thickening $\Th$ of an antipodal set $\La$, the natural projection $\Th(\La)\to \La$ sending $\xi\in \Th(\La)$ to the unique $\la\in \La$ with $\delta(\xi,\la)\in \Th$, is a topological fibration,
compare Lemma~\ref{lem:sldisj} and the discussion afterwards.
We use this fact for chamber limit sets of RA (regular and antipodal) subgroups.
For RCA subgroups,\footnote{Recall that, 
in addition to being regular and antipodal, RCA subgroups are also expanding at $\La_{ch}(\Ga)$,
compare \S \ref{sec:conical_convergence2}.
The RCA property is equivalent to the Anosov property, cf.\ Theorem~\ref{thm:main}.}
we have the following counterpart to the above proper discontinuity result:

\begin{thm}[Cocompactness \cite{coco13}]\label{thm:CC}
Suppose that $\Ga<G$ is an RCA subgroup
and $\Th\subset W$ is a slim thickening. 
Then the action $\Ga\acts\Om_{Th}(\Ga)$ is  cocompact.
\end{thm}

If one works with balanced thickenings which, by definition, are both fat and slim,
one can conclude both proper discontinuity and cocompactness for suitable classes of discrete subgroups:

\begin{cor}\label{cor:flPDCC}
If $\Ga<G$ is an RCA subgroup
and $\Th\subset W$ is a balanced thickening,
then the action $\Ga\acts\Om_{Th}(\Ga)$ is properly discontinuous and cocompact.
\end{cor}

We obtain such results more generally 
for $\taumod$-RCA subgroups acting on partial flag manifolds $G/P_{\taumod}$,  
see \cite{coco15}.

\begin{example}
Let $\Ga< PO(2,1)$ be a cocompact Fuchsian group. 
Then $\Ga<PO(2,1)< PGL(3,\R)$ is Morse, preserves  the Klein model $\H^2$ of the hyperbolic plane  in $\R P^2$.
The hyperbolic plane  $\H^2\subset \R P^2$ is bounded by an ellipse $E$, cf. Example~\ref{ex:ellth}. The group $\Ga$ 
acts properly discontinuously on $\H^2$ and ergodically on the complement, since the latter is $\Ga$-invariantly isomorphic to the space of unparameterized geodesics in the hyperbolic plane. In particular, $\Ga$ does not act properly discontinuously on the complement of 
$E$ in the  projective plane. However, $\Ga$ acts properly discontinuously and cocompactly on the complement $\Flag(\R^3) - Th(\tilde{E})$, where $\tilde{E}=\La_{ch}(\Ga)$ 
is the lift of $E$ to the flag manifold of $PGL(3,\R)$ and $\Th$ is the unique balanced thickening, see Example ~\ref{ex:ellth}. 
\end{example}

Given the above results,
the question arises if and when $\Om_{Th}(\Ga)$ is nonempty. 
Note that in rank 1, the domain of discontinuity in $G/B=\geo X$ is empty in the case of lattices $\Ga<G$. 
In contrast,
it turns out that in higher rank 
our domains $\Om_{Th}(\Ga)$ 
for RA subgroups $\Ga$ and balanced thickenings $\Th$ 
have a tendency to be nonempty.
Intuitively,
the reason is that the emptiness of such a domain 
would imply the existence of certain ball packings at infinity,
e.g.\ of a packing of $\DF X$ by the combinatorial ``balls'' $\Th(\la)$ for $\la\in\La_{ch}(\Ga)$,
and such packings do not exist for many Weyl groups. 
We show:

\begin{thm}[Nonemptiness \cite{coco13,coco15}]
\label{thm:nonempt}
Suppose that $X$ has at least one de Rham factor not of the type $A_1, B_2$ or $G_2$.
Then for each RA subgroup $\Ga< G$, 
there exists a balanced thickening $\Th\subset W$ for which $\Om_{Th}(\Ga)$ is nonempty.  
\end{thm}
 
\begin{rem}
For some Lie groups $G$ of type $B_2$, we can still prove nonemptiness of  $\Om_{Th}(\Ga)$ 
for some balanced thickenings $\Th$ 
(independent of the discrete group $\Ga<G$).
This includes $O(n,2)$ with $n$ odd. 
See \cite{coco15}.
\end{rem} 

Now we can explain the analogy with GIT: For any balanced thickening $\Th$, 
the domain $\Om_{Th}(\Ga)$ serves as the set of stable points for the $\Ga$-action on $G/B$, 
while the thickening of the limit set $\Th(\La_{ch}(\Ga))$ plays the role of the set of unstable points,
and the limit set $\La_{ch}(\Ga)$ itself of the set of maximally unstable points. 

\begin{rem}
Comparison of our discontinuity and cocompactness results with that of \cite{GW}: 

1. Our treatment of domains of discontinuity is intrinsic, 
while in \cite{GW} 
first a theory for $P$-Anosov subgroups of $Aut(F)$ is developed 
(where the $F$'s are certain bilinear and hermitian forms) 
and then general semisimple Lie groups are embedded into groups $O(p,q)$. 

2. Due to the intrinsic nature of our construction, we gain much better control of the nature of domains of proper discontinuity which allows us to get them in $G/B$ (and other flag manifolds) instead of $G/AN$ as in \cite{GW}, for general semisimple Lie groups. (Note, however, that in the case of ``classical''  Lie groups, \cite{GW} also obtain a domain of cocompactness and proper discontinuity inside $G/B$.) 

3. While for some Lie groups of types $A_2,B_2$
the outcomes of the two constructions are the same, it appears that our construction is more general. For instance, we expect that  discontinuity domains constructed via the  
two non-metric balanced thickenings for $SL(4,\R)$,
see Figure \ref{S4.fig}, cannot be obtained via the construction in \cite{GW}. 

4. Theorem~\ref{thm:nonempt} is both weaker and stronger than the
nonemptiness results in \cite[Thms.\ 1.11, 1.12 and 9.10]{GW}.
It is stronger in the sense that it applies to hyperbolic groups $\Ga$ without assumptions on their cohomological dimension,
unlike the results in \cite{GW} which require small cohomological dimension. 
On the other hand, it is weaker in the sense that it addresses only the $\simod$-regular case.
We also note that some examples of Anosov subgroups for which some discontinuity domains are empty are given in \cite[Remark 8.5]{GW}.
\end{rem}

\subsection{Example: Thickenings in visual boundaries of products of rank one spaces}

In this section we work out in detail the case when $W=\Z_2^n$. We identify $\Z_2$ with the multiplicative group $\{-1, 1\}$. Elements of $W$ are identified with $n$-tuples of $\pm 1$'s. The model flat is $\R^n$ and the generators of $W$ act via reflections in the coordinate hyperplanes (walls). We choose the fundamental chamber $\Delta$ to be the orthant 
given by the inequalities $x_i\ge 0$, $i=1,\ldots,n$ (it is clearly a fundamental domain for the action of $W$ on $\R^n$). 
The {\em central direction} in $\Delta$ is given by the  vector 
$$
\bar\zeta= (1,\ldots, 1). 
$$ 
The longest element $w_0=(-1,\ldots, -1)$ acts as $-\id$.

The Bruhat order is given by
$$
w=(\eps_1,\ldots,\eps_n)\le w'= (\eps_1',\ldots,\eps_n') \quad \iff \quad \eps_i\ge \eps'_i \;\forall i.  
$$

Examples of thickenings are given by strict and nonstrict linear  inequalities as follows. Let $a=(a_1,\ldots,a_n)$ be a vector with (strictly) positive entries. The subsets  
$$
{\Th}_a=\{w\in W: a\cdot w> 0\}, \overline{\Th}_a=\{w\in W: a\cdot w\ge 0\}
$$
are {\em metric thickenings}. The former thickening is slim while the latter is fat. A thickening $\Th_a$ is balanced iff $a$ does not satisfy an equation
$$
\sum_{i\in I} a_i = \sum_{j\notin I} a_j,
$$
for any subset $I\subset \{1,\ldots,n\}$.  Hence, for ``generic'' values of $a$, $\Th_a$ is balanced. 

Consider now a rank one symmetric space $Y$ (e.g.\ $Y=\H^2$) and $S=\geo Y$. 
Let $X=Y^n$, the $n$-fold product of $Y$. 
Then $Z:=\DF X=S\times ... \times S$, the $n$-fold product of $S$. 
Moreover, let $D\subset Z$ denote the diagonal
$$D=\{(s, \ldots,s): s\in S\}. $$
We will think of elements of $Z$ as configurations of points in $S$.

The relative position of two configurations $z=(s_i)$ and $z'=(s'_i)$ equals $\de(z',z)=(\eps_i)$
iff $$ s'_i=s_i\iff\eps_i=+1 ,$$
i.e.\ $\delta$ records the entries $i$ where $z'$ agrees with $z$. 
Consequently, 
$\de(z',z)\leq(\eps_i)$
iff $s'_i=s_i$ whenever $\eps_i=+1$.

The vector $a$ assigns {\em weights} $a_i$ to the $i$-th members of the configuration. Each weighted configuration $z=(s_1,\ldots, s_n)$ thus gives rise to a finite measure $\mu$ on $S$,
$$
\mu_z= \sum_{i=1}^n a_i \delta_{s_i},
$$ 
where $\delta_s$ is the probability measure on $S$ supported at the point $s$ (masses add when points $s_i$ ``collide''). The total mass of $\mu_z$ equals
$$
M=a_1+\ldots+a_n. 
$$
A weighted configuration $z$
is called  {\em stable} if $\mu_z(s)< M/2$ for all points $s\in S$, 
and {\em semistable} if $\mu_z(s)\le M/2$ for all $s\in S$. 
In the balanced case, these notions agree: ``stable=semistable.'' It is then immediate that 
$$
z\in \Th_a(D) \iff \mu_z \hbox{~~is not semistable}
$$
and 
$$
z\in \overline{\Th}_a(D) \iff \mu_z \hbox{~~is not stable}.  
$$
The sets of stable and semistable weighted configurations are denoted $Z_{st}$ and $Z_{sst}$. 
They of course depend on $a$. 
For instance, if $a_i> M/2$ for some $i$,
then $Z_{sst}=\emptyset$. 
On the other hand,
if $a_i< M/2$ for all $i$,
then $Z_{st}\neq\emptyset$;
e.g.\ all configurations of pairwise distince points $s_i$ are stable. 

\medskip
Assume now that $H$ is the isometry group of $Y$ acting diagonally on $X$ and, hence, on $Z$. 
The latter action preserves the diagonal $D$,
which we can regard as the chamber limit set of the Lie subgroup $H<G:=\Isom(X)\cong H^n$,
$D=\La_{ch}(H)$.
Mumford's GIT defines the {\em Mumford quotient} 
$$
Z//_a H= Z_{sst}//H.
$$  
In the balanced case, we simply have
$$
Z//_a H= Z_{sst}//H= Z_{st}/H. 
$$
A nice exercise is to prove directly that $Z_{sst}//H$ is compact and Hausdorff in this case. For instance, if $H=PSL(2,\R)$, $Y=\H^2$, $n=3$ and $a=(1,1,1)$ then $Z//_a H$ consists of exactly two points represented by configurations of three distinct points on the circle with different cyclic orders. Continuing with $Y=\H^2$ and letting $n=4$, one verifies that for $a=(2,1,1,1)$ the Mumford quotient  is homeomorphic to $S^1$, while for $a=(5,4,3,1)$ the 
Mumford quotient is homeomorphic to the disjoint union of two circles. Taking $n=5$, one obtains that for 
$a=(1,1,1,1,1)$ the Mumford quotient is the genus 4 oriented surface, while for $a= (5,4,1,1,1)$ the quotient is the disjoint union of two 2-spheres. Thus, we see that quotients can be non-homeomorphic for distinct choices of $a$. 
We refer the reader to \cite[Theorem 2]{KM} for  proofs of these descriptions of Mumford quotients using their identification with polygon spaces. 

\begin{rem}
The hyperplanes $\sum_{i\in I} a_i = \sum_{j\notin I} a_j$ (called {\em walls}), 
where $I$ runs through the subsets of $\{1,\ldots,n\}$, 
partition the space
$$
A=\{(a_1,\ldots,a_n): a_i>0\}
$$ 
into open convex subsets called {\em chambers} (they are not fundamental domains for the $S_n$-action!). 
The topology of $Z//_a H$ does not change as long as $a$ varies in a single chamber; 
permuting the chambers does not change the topology either; 
however, {\em crossing through a wall} amounts to a certain Morse surgery on the manifold  $Z//_a H$. 
This can be seen  by identifying the quotients $Z//_a H$ with certain spaces of polygons with fixed side lengths: 
In the case when $H=PSL(2,\R)$, these are polygons in the euclidean plane, see \cite{KM}. 

It was conjectured by Kevin Walker (in his undergraduate thesis written in 1986 under Bill Thurston; Walker was working with euclidean polygons) that, for $n\ge 5$,  if $a, a'$ belong to chambers in distinct $S_n$-orbits, then the Mumford quotients are not homeomorphic. This conjecture was proven 20 years later in ``most'' cases by Farber, Hausmann and Sch\"utz \cite{FHS} and in full generality by  Sch\"utz \cite{Schutz}. Similar results hold when the circle is replaced by a $k$-sphere; in fact, different quotients are distinguished by their intersection cohomology rings, see \cite{Schutz1, Schutz2}. 
\end{rem}

Now, suppose that $\Ga< H$ is a uniform lattice.
Then $\La_{ch}(\Ga)=D$. 
The subgroup $\Ga$, diagonally embedded in $G$,  is uniformly regular in $G$: $\Ga$ preserves the diagonally embedded copy of $Y$ in $X$, 
and any geodesic segment in it has $\De$-length contained in the diagonal of $\De\cong[0,\infty)^n$.
We conclude that $\Ga$ is $\Theta$-regular with $\Theta$ consisting of a single point,
namely the center of the model spherical chamber of $X$, represented by the unit vector 
$$
\frac{1}{\sqrt{n}}(1,\ldots, 1). 
$$
The group $\Ga$ is quasiisometrically embedded in $H$ and hence in $G$. Thus, $\Ga<G$ is URU.

Given a balanced metric thickening $\Th=\Th_a$, the domain of discontinuity $\Om_{Th}(\Ga)$ equals the set $Z_{st}$ 
of stable weighted $n$-point configurations  in $S$ (stability is, of course, defined with respect to $a$). 

We now specialize to the case when 
$H=PSL(2,\R)$ or $PSL(2,\C)$
and $\Ga$ is torsion-free.
Then the group $H$ 
acts freely and properly on $Z_{st}$, and we have a principal $H$-bundle
$$
H\to Z_{st}\to Z_{st}/H= Z//_a H.  
$$
Dividing $Z_{st}$ by $\Ga$ instead of $H$, we obtain a fiber bundle
$$
F \to Z_{st}/\Ga\to Z_{st}/H, 
$$
with fiber $F=\Ga\backslash H$, the oriented orthonormal frame bundle over the manifold $Y/\Ga$.  In particular, by taking non-homeomorphic Mumford quotients $Z_{st}/H$, 
we may obtain non-homeomorphic quotients $\Om_{Th}/\Ga=Z_{st}/\Ga$. 
For instance, 
taking $H=PSL(2,\R)$ and $n=4$, we obtain three distinct topological types of quotients: The empty quotient, 
a connected nonempty quotient (a bundle over the circle with fiber $F$, the unit tangent bundle of a hyperbolic surface) and a disconnected quotient 
(an $F$-bundle over $S^1 \sqcup S^1$).

\subsection {Finsler bordifications of  locally symmetric spaces} \label{sec:bordif}

For a regular subgroup $\Ga< G$ and a thickening $\Th\subset W$, we define the {\em Finsler thickening} 
of the chamber limit set $\La_{ch}(\Ga)$ as follows. 
First, 
recall the definition \eqref{eq:FTH} of the thickening $\Th(\si)\subset \DF X$ of a chamber $\si\in \DF X$ 
inside the Furstenberg boundary,
and the definition \eqref{eq:star} of the star $\st(\tau)$ of a simplex $\tau$.
We then introduce the {\em Finsler thickening} of the chamber $\si$ as the union of small strata
$$\Th^{Fins}(\si)= \bigcup\bigl\{ S_{\tau}: \st(\tau)\subset \Th(\si)\bigr\} \subset \geo^{Fins}X .$$ 
Finsler thickenings of antipodal chambers are disjoint if $\Th$ is slim. 

We obtain the {\em Finsler thickening} of the chamber limit set $\La_{ch}(\Ga)$
by taking the union of the Finsler thickenings of all limit chambers, 
$$
\Th^{Fins}(\La_{ch}(\Ga))=\bigcup_{\si\in \La_{ch}(\Ga)} \Th^{Fins}(\si)\subset \geo^{Fins}X. 
$$
This subset is closed, $\Ga$-invariant and {\em saturated}, i.e. a union of small strata $S_{\tau}$.
We consider the domain at infinity
$$\Om_{Th}^{Fins}(\Ga)= \geo^{Fins} X - \Th^{Fins}(\La_{ch}(\Ga))$$
and the domain
$$  X\sqcup \Om_{Th}^{Fins}(\Ga) = \ol{X}^{Fins}- \Th^{Fins}(\La_{ch}(\Ga)).$$
Recall from section~\ref{Finsler-com} 
that the Furstenberg boundary sits inside the Finsler boundary (as a big stratum),
$\DF X\subset \geo^{Fins}X$,
and note that our domains in the latter extend the domains in the former,
$$ \Om_{Th}^{Fins}(\Ga) \cap\DF X = \Om_{Th}(\Ga),$$ 
because 
$\Th^{Fins}(\si)\cap\DF X=\Th(\si)$.

Theorems~\ref{thm:finsPD}, \ref{thm:Fcocom}
and Corollary~\ref{cor:finsPDCC}  below are 
{\em Finsler extensions} 
of Theorems \ref{thm:proper}, \ref{thm:CC}
and Corollary \ref{cor:flPDCC} about discrete group actions on the Furstenberg boundary $\DF X\cong G/B$.

\begin{thm}[Finsler domains of proper discontinuity {\cite[Theorem 9.13]{bordif}}]
\label{thm:finsPD}
Suppose that $\Ga< G$ is regular
and $\Th\subset W$ is a fat thickening.
Then the action 
$$ \Ga\acts X\sqcup \Om_{Th}^{Fins}(\Ga) $$
is properly discontinuous.
\end{thm}
We note that our construction of domains provides, more generally,
domains of proper discontinuity for the action 
of {\em arbitrary} discrete subgroups $\Ga< G$ on $\ol{X}^{Fins}$,
not only for subgroups which are $\taumod$-regular 
for some $\taumod$ 
(see Theorems 9.16 and 9.18 of \cite{bordif}).
These more general domains involve complements to unions of Finsler thickenings of $\taumod$-limit sets 
of the subgroups $\Ga$ with $\taumod$ running through all the faces of $\simod$. 

\begin{thm}[Nonemptiness {\cite[Prop.\ 9.20]{bordif}}]\label{thm:finsNE}
Suppose that $\Ga< G$ is an RA subgroup, 
$\Th\subset W$ is a slim thickening and $\rank(X)\geq2$.
Then $\Om_{Th}^{Fins}(\Ga)$ is nonempty.
\end{thm}

Note that, unlike Theorem \ref{thm:nonempt}, this result does not exclude products of 
symmetric spaces of type $B_2$ and $G_2$. It is also not limited to $\simod$-regular subgroups, but holds for all 
$\taumod$-regular antipodal subgroups.

In order to address cocompactness,
we convert the action $\Ga\acts \ol{X}^{Fins}$ to a topological convergence group action via the following $\Ga$-invariant collapsing procedure: 
Form a quotient of $\ol{X}^{Fins}$ 
by simultaneously collapsing the thickenings $\Th^{Fins}(\si)$ for all $\si\in \La_{ch}(\Ga)$ 
to points. Let $Z$ denote the resulting quotient space, 
and $\La$ the projection of 
$\Th^{Fins}(\La_{ch}(\Ga))$ to $Z$. Then $\La$ is equivariantly homeomorphic to $\La_{ch}(\Ga)$. 

\begin{thm}[{\cite[Corollary 11.7, Lemma 11.9]{bordif}}]\label{thm:ca}
If $\Ga<G$ is an RA subgroup and $\Th\subset W$ is a balanced thickening,
then the (obviously compact) quotient space $Z$ is metrizable 
and $$\Ga \acts Z$$ is a convergence group action with limit set $\La$.  
\end{thm}

The last theorem is yet another indication of the ``rank 1 nature'' of RA subgroups $\Ga< G$. 
It is used in \cite{bordif} to prove:

\begin{thm}[Finsler cocompactness {\cite[Theorem 11.11]{bordif}}]
\label{thm:Fcocom}
Suppose that $\Ga<G$ is an RCA subgroup and $\Th\subset W$ is a slim thickening.
Then the action $\Ga\acts X\sqcup \Om_{Th}^{Fins}(\Ga)$
is cocompact. 
\end{thm}

Combining Theorems \ref{thm:finsPD} and \ref{thm:Fcocom} we obtain: 

\begin{cor}\label{cor:finsPDCC}
If $\Ga<G$ is an RCA subgroup
and $\Th\subset W$ is a balanced thickening,
then the action $\Ga\acts X\sqcup \Om_{Th}^{Fins}(\Ga)$ is properly discontinuous and cocompact.
\end{cor}

Note that in this result in the $\simod$-regular case one does not need antipodality of the limit set to conclude proper discontinuity and cocompactness, provided that $\Th$ is a {\em metric thickening} associated with a {\em nearly root element} $\bar\theta$ { (see the May 2015 version of the preprint \cite{bordif} for the details)}. 
The RCA assumption is, however, needed in the general $\taumod$-regular case. 

We apply our construction of domains to obtain bordifications and compactifications of locally symmetric spaces: 

\begin{cor}
\label{cor:lcsm}
1. For each regular subgroup $\Ga< G$, the locally symmetric orbifold 
$X/\Ga$ admits a real-analytic bordification as an orbifold with corners\footnote{{See Appendix \ref{sec:corners} for the precise definition.}} 
$$
\left(X\sqcup \Om_{Th}^{Fins}(\Ga)\right)/\Ga, 
$$
provided that $\Th\subset W$ is fat. When this quotient is treated as an orbifold with boundary, the boundary of this 
orbifold is  $(\Om_{Th}^{Fins}(\Ga))/\Ga$. 

2. If $\Ga$ is RCA and $\Th$ is balanced, then this bordification of $X/\Ga$ is a compact orbifold with corners. \end{cor}

\begin{rem}
This corollary implies the {\em topological tameness} of the orbifold $X/\Ga$. However, topological tameness is a weaker property than the existence of a compactification given by the corollary. For instance, considering finitely generated discrete subgroups $\Ga< PSL(2,\C)$, all quotient spaces $\H^3/\Ga$ of such groups are topologically tame, but for many groups $\Ga$ the bordification
$$
(\H^3\cup \Om(\Ga))/\Ga
$$ 
is not compact. The latter happens, for instance, for singly degenerate groups. 
\end{rem}

We show furthermore a converse to the cocompactness part of  Theorem \ref{thm:Fcocom}, 
implying that Anosov subgroups are characterized among uniformly regular subgroups by 
the cocompactness of their action on complements to balanced Finsler thickenings.
More generally, we consider the following property 
of admitting cocompact domains of proper discontinuity in the Finsler compactification:

\begin{dfn}
[$S$-cocompact {\cite[Def.\ 12.4]{bordif}}]
We say that a discrete subgroup $\Ga<G$ is {\em $S$-cocompact} 
if there exists a $\Ga$-invariant saturated open subset $\Om_{\infty}\subset\geo^{Fins}X$ 
such that the  action 
$\Ga\acts X\sqcup\Om_{\infty}$ 
is properly discontinuous and cocompact. 
\end{dfn}

\begin{remark}
{The terminology $S$-cocompact comes from ``saturated'', 
although the letters S from ``Satake'' and  ``stratified'' also appear naturally in this context.}  
\end{remark} 

A useful implication of $S$-cocompactness is given by:
\begin{thm}[Cocompactness implies retract {\cite[Thm. 12.5]{bordif}}]
\label{thm:S->retract}
$S$-cocompact discrete subgroups $\Ga<G$ are coarse retracts\footnote{cf.\ Definition \ref{defn:retract}}. 
In particular, they are undistorted.
\end{thm}

Combining this theorem with the fact that URU subgroups are Anosov, we obtain:
\begin{thm}[Cocompactness implies Anosov {\cite[Thm.\ 1.9]{bordif}}]
\label{thm:URS->A}
$S$-co\-com\-pact uniformly regular subgroups $\Ga<G$ are  Anosov.
\end{thm}

We conclude:

\begin{cor}
[Dynamical characterizations of Anosov subgroups II: actions on Finsler compactifications 
{\cite[Cor.\ 1.10]{bordif}}]\label{cor:S-coco}
For uniformly regular subgroups $\Ga<G$,
the following properties are equivalent:

(i) Anosov

(ii) $S$-cocompact

(iii) coarse retract
\end{cor}

Combining Corollary \ref{cor:S-coco} with Theorem \ref{thm:main}, we obtain 
a higher rank analogue of the Equivalence Theorem for convex cocompact groups of isometries of rank 1 symmetric spaces, see Theorem \ref{thm:main1}, with the conditions CC0, CC1 and CC8  excluded as inappropriate in higher rank. 

\begin{figure}[tbh]
\includegraphics[width=90mm]{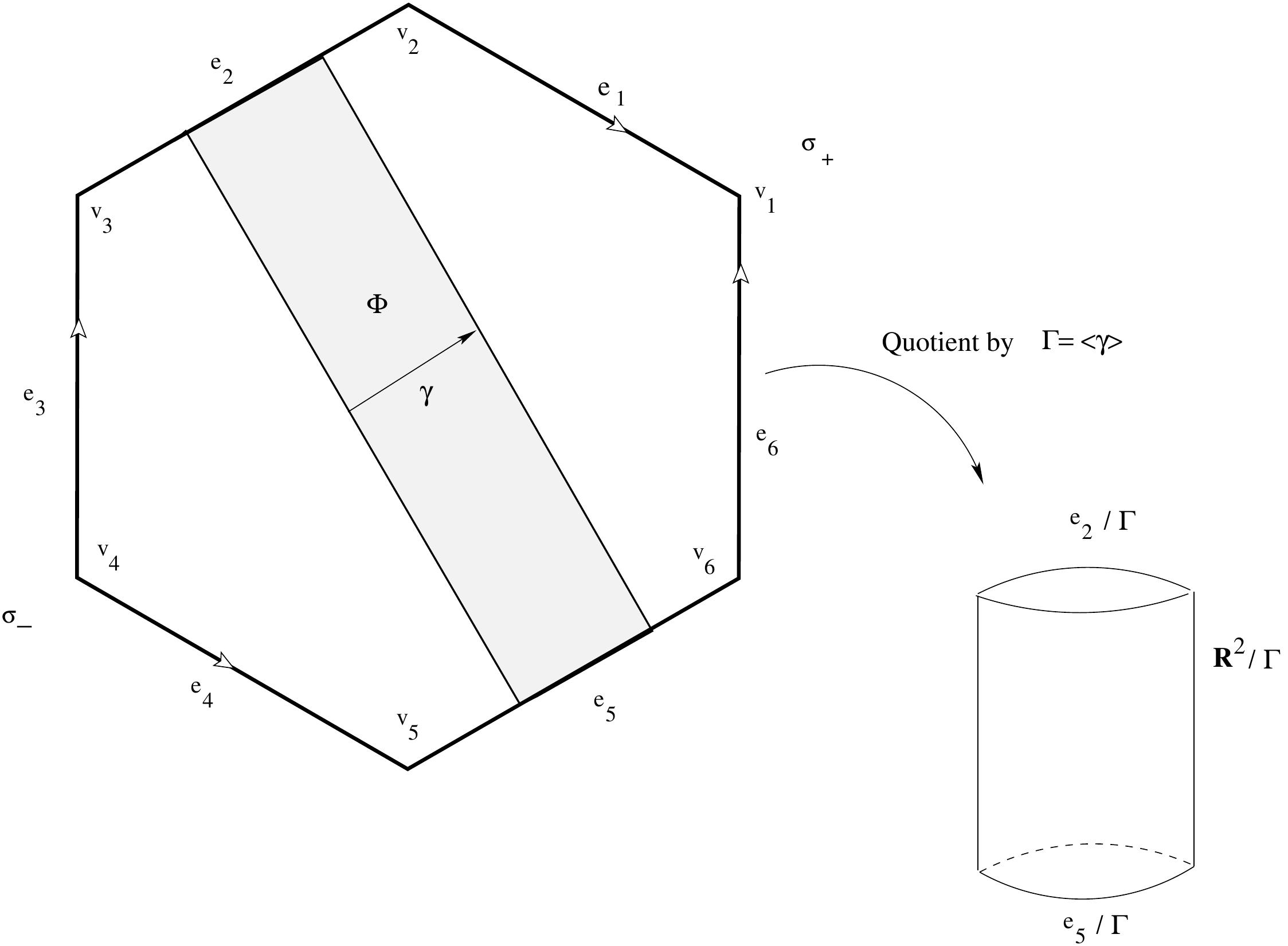}
\caption{The action of a cyclic subgroup $\Ga=\<\ga\>$ on the Finsler compactification $\ol{F}^{Fins}_{mod}$ 
of the model flat and the quotient space of $\ol{F}^{Fins}_{mod}- \Th^{Fins}(\La_{ch})$ by the $\Ga$-action.}
\label{figure16.fig}
\end{figure}

\begin{example}
We now work out an example illustrating these results. 
Consider an infinite cyclic subgroup $\Ga=\<\ga\><PGL(3,\R)$ 
generated by a regular hyperbolic isometry $\gamma$. 
For simplicity, we only describe the action on the Finsler compactification 
of the unique invariant maximal flat $F\subset X$.
The Finsler compactification $\ol F^{Fins}$ is a hexagon with vertices $v_1,\ldots,v_6$ and edges $e_1,\ldots,e_6$. 
The vertex set equals the Furstenberg boundary,
$\DF F=\{v_1,\ldots,v_6\}$.
We label the vertices so that $v_1$ and $v_4$ correspond to the repulsive and attractive chambers $\si_-,\si_+\in\DF F$. 
The vertices are fixed by $\ga$, but $\ga$ has nontrivial dynamics on the edges:
The interior points of each edge $e_i=[v_i, v_{i+1}]$ are moved by $\ga$ towards one of the two endpoints of $e_i$, 
namely to the one which corresponds to the chamber in $\DF F$
whose position relative to the attractive chamber $\si_+$ 
is smaller in the Bruhat order. 
This is in stark contrast with the action of $\ga$ on the visual boundary of $\geo F$ (with respect to the flat metric), 
which is fixed pointwise. 
The chamber limit set $\La_{ch}(\Ga)\subset\DF X$ 
is the 2-point set $\{\si_-, \si_+\}=\{v_1, v_4\}\subset\DF F$. 
The balanced thickening of  $\La_{ch}(\Ga)$ inside $\geo^{Fins} F$  is the union (of closed edges)
$$
\Th^{Fins}(\si_-) \cup \Th^{Fins}(\si_+)= \left( e_3 \cup e_4\right) \cup  \left( e_1 \cup e_6\right)
$$
The intersection
$$
\Om=\Om(\Ga)=  \Om_{Th}^{Fins} (\Ga) \cap \geo^{Fins}F
$$
is the union of the interiors of the edges $e_2$ and $e_5$. 
The rectangle $\Phi$ in Figure \ref{figure16.fig} is a (compact) fundamental domain for the action of $\Ga$ on $F\cup \Om$. 
The quotient $\Om/\Ga$ is homeomorphic to the cylinder $S^1 \times [-1,1]$. 
Now, let us collapse each thickening 
$\Th^{Fins}(\si_-), \Th^{Fins}(\si_+)$ to a point. The result is a convergence action of $\Ga$ on the quotient space $Q$, 
homeomorphic to the closed 2-disk $D^2$. {Note that collapsing is  natural here since, before the collapse,
the mapping $\gamma$ has too many fixed points in $\geo^{Fins}F$, namely all vertices $v_1,...,v_6$, while 
an infinite cyclic group acting as a discrete convergence group can have at most two fixed points \cite{Tukia1994}.  
After the collapse only two fixed points are left, namely the projections (still denoted $\si_+, \si_-$) of $v_1$ and $v_4$.
On the quotient space $Q$ we recover the familiar attractive-repulsive dynamics of hyperbolic isometries $\ga$ of $\H^2$ acting on the visual compactification of $\H^2$: The point $\si_+$ is the attractive point and the point $\si_-$ is the repulsive point for the action of $\ga$: 
$$
\lim_{n\to\infty} \ga^n= \si_+, 
$$
uniformly on compacts in $Q - \{\si_-\}$, and 
$$
\lim_{n\to-\infty} \ga^{-n}= \si_-, 
$$
uniformly on compacts in $Q - \{\si_+\}$.} 
\end{example}

\begin{example}
[A product example] We continue with Example \ref{ex:product-case}
of a cyclic isometry subgroup $\Ga$ of the product $X=X_1\times X_2$ of two hyperbolic spaces. 
The Finsler compactification of $X$ is naturally homeomorphic to $\ol{X}_1\times \ol{X}_2$.  
Assume that $g=(g_1, g_2)$ where $g_1$ is hyperbolic (with the fixed points $\la_1^+, \la_1^-$) and $g_2$ is parabolic (with the fixed point $\la_2$). As we noted in  Example \ref{ex:product-case}, the group $\Ga=\<g\> < \Isom(X)$ is regular but not uniformly regular.  Therefore, it is not Anosov. On the other hand, it  is $S$-cocompact. Namely, it acts properly discontinuously and cocompactly on 
$$
(\ol{X}_1 - \{\la_1^-, \la_1^+\} ) \times \ol{X_2}. 
$$
In particular, 
$\Ga$ is a coarse retract, and hence undistorted.  Thus, uniform regularity cannot be weakened to regularity 
in Theorems \ref{thm:main} (item 8), \ref{thm:URS->A} and Corollary~\ref{cor:S-coco}.
\end{example}

\begin{figure}[tbh]
\includegraphics[width=90mm]{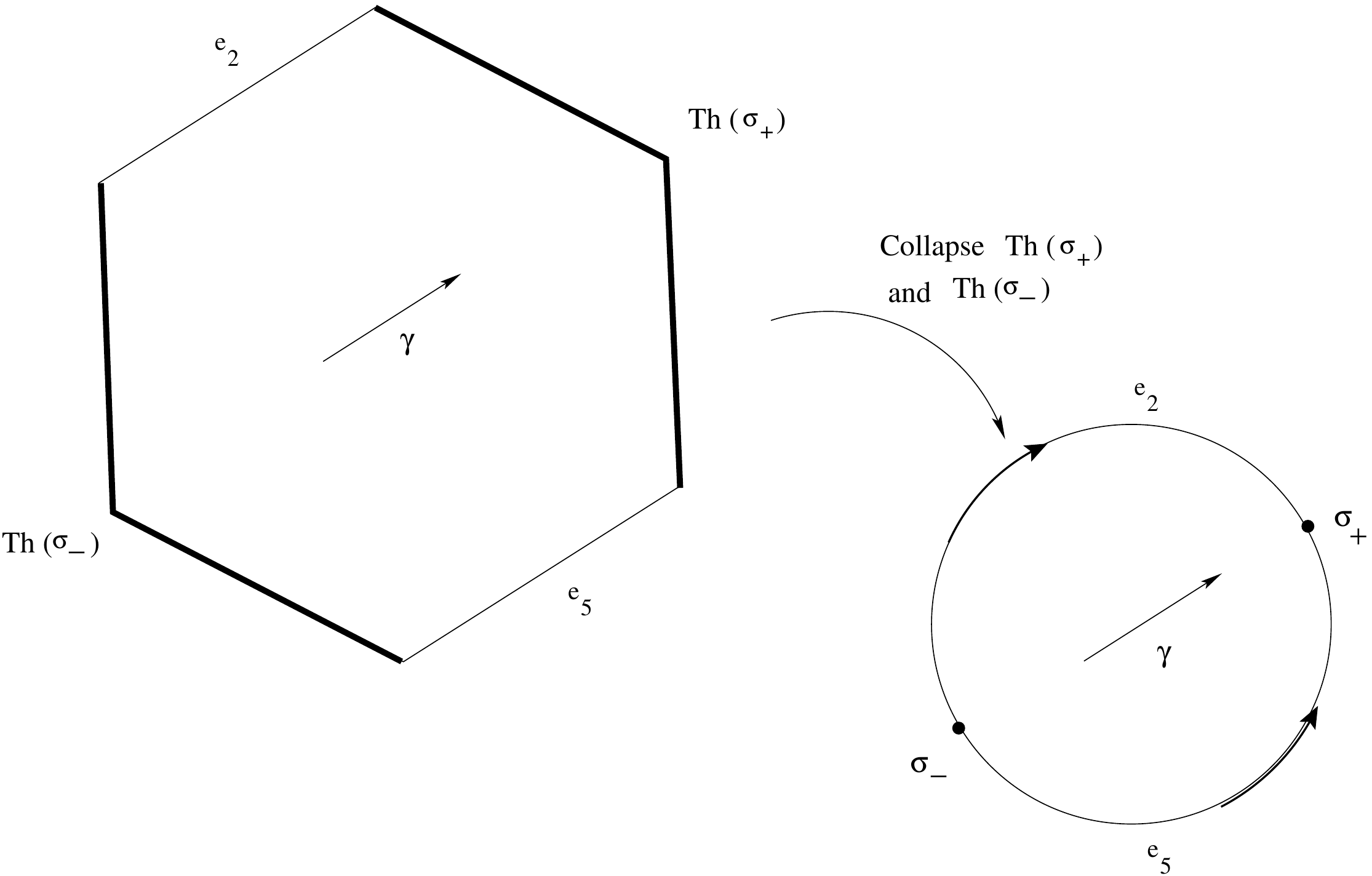}
\caption{Collapsing thickened limit set for the action of a cyclic subgroup $\Ga=\<\ga\>$ on the Finsler compactification of the model flat.}
\label{figure17.fig}
\end{figure}

\begin{rem}[Relation with the work \cite{GW,GGKW2}]
(i)
We note that
the existence of an orbifold with boundary compactification of locally symmetric quotients by
Anosov subgroups of some special classes of simple Lie groups (namely, $Sp(2n,\R), SU(n,n), SO(n,n)$)
appeared in \cite{GW}.

(ii)
The main results of this section 
dealing with Finsler bordifications of locally symmetric spaces 
(Theorems \ref{thm:finsPD}, \ref{thm:Fcocom}
and Corollaries~\ref{cor:finsPDCC}, \ref{cor:lcsm})
are contained in the second version of \cite{bordif}.
After that work had been completed, the e-print \cite{GGKW2} was posted,
also addressing the compactification of locally symmetric spaces.
Theorem~1.2 there provides orbifolds with corners compactifications of $X/\Ga$
via the maximal Satake compactification of $X$
for $\taumod$-Anosov subgroups $\Ga<G$ of special face types $\taumod$.
However,
Theorem~1.1 of \cite{GGKW2} dealing with the general case
still lacks a complete proof.
It remains unclear
whether the compactifications constructed there
are orbifolds with corners.
Namely,
the approach uses ``generalized'' Satake compactifications,
but the proof of Lemma A9, establishing that the latter are manifolds with corners,
lacks details.
Note also that in the first version of \cite{GGKW2}
there was a {basic} mistake in the cocompactness argument.
It was corrected in the third version
using methods from \cite{coco13}.
\end{rem}

\bigskip 
\section{Future directions}\label{sec:5}

{\bf Regular antipodal subgroups.} One can think of the class RA of {\em regular antipodal} (or URA: uniformly regular antipodal) subgroups of $G$ as {\em discrete subgroups exhibiting rank 1 behavior}: We saw several examples of this at work. Dropping conicality, we obtain (ignoring the issue of parabolic elements) an analogue of (rank 1) Kleinian groups without any geometric finiteness assumptions. Quite likely, this class by itself deserves some attention.  More  generally, one can define {\em uniformly rank $k$} discrete subgroups $\Ga< G$ as those whose visual limit sets $\La(\Ga)\subset \geo X$ are disjoint from the codimension $k$ skeleton of the Tits boundary of $X$. 
A {\em limit simplex} of $\Ga$ is a simplex $\tau\subset\geo X$ which contains a limit point.  

\medskip 
{\em Geometric finiteness} by no means should be limited to subclasses of uniformly regular 
(and, more generally, uniformly rank $k$) 
discrete subgroups. Here is a {\em wish list} for a good 
notion of geometric finiteness in higher rank:

\medskip 
\noindent {\bf A. Conjectural properties of geometrically finite subgroups:} 

1. Geometrically finite groups should be stable under small deformations, provided that we have the {\em right algebraic restrictions} (yet to be determined) on the deformations of representations. In rank 1, such restrictions amount to a certain control on the deformations of maximal parabolic subgroups, see \cite{Bowditch-stab}. 

2. The  locally symmetric quotient spaces $X/\Ga$ of geometrically finite groups should admit geometrically natural compactifications, 
by attaching quotients of domains of discontinuity at infinity and compactifying ``cusps'', incorporating (as a special case) the Borel--Serre compactifications (see e.g. \cite{Borel-Ji}) 
of locally-symmetric spaces of finite volume. In particular, such groups should be finitely presented. 

3. Algebraically speaking, geometrically finite groups should be (suitably relativized) semihyperbolic groups. 
Note that semihyperbolic groups (cf. \cite[III.$\Gamma$.4]{BH}) represent a coarsification of the notion of CAT(0) groups. Any notion of relative semihyperbolicity should include relatively hyperbolic groups as a special case. 
We refer the reader to \cite{KR} for one possible definition of relatively semyhyperbolic groups. 

4. Geometric finiteness should be semidecidable. 

5. Geometric finiteness should be stable under embeddings of the ambient Lie groups: 
If $\Ga< G_1$ is geometrically finite and $\phi: G_1\to G_2$ is an embedding of semisimple (or reductive) Lie groups, then $\phi(\Ga)< G_2$ is again geometrically finite. (Note that this fails, in general, for Anosov subgroups.) 

6. Geometric finiteness should be stable under taking finite index subgroups: If $\Ga_1< \Ga_2$ is a finite index subgroup, then $\Ga_2< G$ is geometrically finite iff $\Ga_1< G$ is geometrically finite. 

  \bigskip
  
\noindent {\bf B. A conjectural class of geometrically finite subgroups should include:} 

1. Direct products of geometrically finite groups: If $\Ga_i< G_i$ are geometrically finite, $i=1, 2$, then 
$\Ga_1\times \Ga_2< G_1 \times G_2$ is also geometrically finite.  

2. All geometrically finite groups in rank one (allowing parabolics).  

3. Standard representations of Coxeter groups \cite[Sect. V.4]{Bourbaki}. 

4. Groups of projective transformations acting properly discontinuously on bounded convex domains in the affine space such that the quotient has finite volume with respect to the Hilbert metric. 

5. Lattices in semisimple algebraic Lie groups (of any rank). 

\medskip
At this point, the following two definitions of geometric finiteness, in the setting of discrete groups containing no parabolic elements, appear to be most promising:

\begin{definition}
$\Ga$ is geometrically finite if it is $S$-cocompact. 
\end{definition}

For instance, every convex cocompact subgroup $\Ga< G$ in the sense of \cite{convcoco, Quint} is $S$-co\-com\-pact. 
Every Anosov subgroup $\Ga< G$ is $S$-cocompact as well.

\begin{definition}
$\Ga<G$ is geometrically finite if it is an equivariant coarse retract (cf.\ Definition~\ref{defn:retract}).
\end{definition}

Note that, according to Theorem \ref{thm:S->retract}, $S$-cocompactness implies
being an equivariant coarse retract. 
In particular, either definition implies semihyperbolicity of $\Ga$. 
Moreover, the definition using retractions is clearly stable under embeddings of Lie groups as mentioned above.

\section{Appendix. Horofunction compactification} 
\label{sec:horoboundary}

Let $(Y,d)$ be a locally compact geodesic metric space. For each $y\in Y$ define the 1-Lipschitz function 
$d_y= d(y, \cdot)$ on $Y$.  This leads to the embedding $\kappa: Y\to C(Y)=C(Y,\R)$, $y\mapsto d_y$.  We let $\R\subset C(Y)$ denote the linear subspace of constant functions.  Composing the embedding $\kappa$ with the projection $C(Y)\to C(Y)/\R$ we obtain the {\em Kuratowski  embedding} of $Y$,
$$
Y\embed  C(Y)/\R. 
$$ 
Then $\ol{Y}$, the closure of $Y$ in ${C(Y)}/\R$, is the {\em horofunction compactification} of $Y$. 
Functions representing points in $\geo Y= \ol{Y}- Y$ are the {\em horofunctions} on $Y$. In other words, horofunctions on $Y$ are limits (uniform on compacts in $Y$) of sequences of normalized distance functions $d_{y_i} - d_{y_i}(o)$, where $y_i\in Y$ are divergent sequences  in $Y$. Each geodesic ray $r(t)$ in $Y$ determines a horofunction in $Y$ called a {\em Busemann function} $b_r$, which is the subsequential limit
$$
\lim_{i\to\infty} d_{r(i)} - d_{r(i)}(o). 
$$
If $Y$ is  a CAT(0) space, then each limit as above exists (without passing to a subsequence). Furthermore, each horofunction is a Busemann function.  This yields a topological identification of the visual compactification of $Y$ and its horofunction compactification. Level sets of Busemann functions are called {\em horospheres} in $X$. The point $r(\infty)\in \geo Y$ is the {\em center} of the horosphere $\{b_r=c\}$. We refer the reader to \cite{Gromov_hypmfs, Ballmann} for further details and to \cite{bordif} for the detailed treatment of this construction in the case of nonsymmetric metrics.

\section{Appendix. Expanding and relatively expanding actions} 
\label{sec:expanding_actions}

Let $(Z,d)$ be a compact metric space. A map $f: Z\to Z$ is said to be {\em metrically expanding} at a point 
$z\in Z$ if there exists a neighborhood $U$ of $z$ in $Z$, a number $c > 1$ (an {\em expansion factor}) 
such that for all $z', z''\in U$,
$$
d(f(z'), f(z''))\ge c d(z', z'').  
$$ 
 A sequence of maps $f_n: Z\to Z$ is said to have {\em metrically diverging expansion} at $z\in Z$ if there exists a 
 system of neighborhoods $U_n$ of $z$ and expansion factors $c_n\to \infty$ such that each $f_n|_{U_n}$ expands 
 with the expansion factor $c_n$. 
 
 A topological action $\Ga\acts Z$ is said to be {\em expanding} at a $\Ga$-invariant subset 
$E\subset Z$ if for each $z\in E$ there exists $\ga\in \Ga$ which is expanding at $z\in Z$. 

The expansion concepts have infinitesimal versions in the case of diffeomorphisms and smooth group actions on Riemannian manifolds. Suppose that $M$ is a Riemannian manifold and $f: M\to M$ is a diffeomorphism. The {\em infinitesimal expansion factor} of $f$ at a point $x\in M$ is the number 
$$
\eps(f, x)= \inf_{u\in U_xM} |df(u)| 
$$
where $U_xM$ is the unit sphere in $T_xM$. 
A smooth map $f: M\to M$ is said to be {\em infinitesimally expanding} at $x$ if $\eps(f,x)>1$. It is easily seen that a smooth map is infinitesimally expanding at $x$ iff it is metrically expanding at $x$. A sequence of smooth maps $f_n: M\to M$ is said to have {\em  diverging infinitesimal expansion} at $x\in M$ if 
$$
\lim_{n\to\infty} \eps(f_n, x)=\infty. 
$$
{A group of diffeomorphisms $\Ga< Diff(M)$ is {\em infinitesimally expanding} at a 
subset $Z\subset M$ if for every $z\in Z$ there exists $\ga\in \Ga$ which is infinitesimally expanding at $z$.}

\medskip 
More generally, one defines {\em relatively expanding actions} of groups on metric spaces. Suppose again that 
$(Z,d)$ is a compact metric space, $\Ga\acts Z$ is a topological group action preserving a compact subset $E\subset Z$. 
Suppose, furthermore, that $\pi: E\to\La$ is a continuous map, which is equivariant with respect to actions 
$\Ga\acts E, \Ga\acts \La$. We let $E_\la=\pi^{-1}(\la)$ for $\la\in \La$. The action 
of $\Ga$ is said to be {\em relatively expanding} at $E$ with respect to $\pi$ 
(or {\em expanding relative to} $\pi: E\to \La$) 
if: 

For each $\la\in \La$ there exists a neighborhood $U_\la$ of $E_\la$ in $Z$, a number $c>1$ and an element $\ga\in \Ga$ such that for all $E_{\la'}\subset U$ and $z\in U$,
$$
d(\ga(z), \ga E_{\la'})\ge c d(z, E_{\la'}).  
$$ 
Here the distance $d(z, W)$ froma  point $z\in Z$ to a subset $W\subset Z$ is
$$
d(z, W):= \inf_{w\in W} d(z,w). 
$$

Such relatively expanding actions frequently appear with the sets $E_\la$ being {\em stable sets} of the action, when the dynamics of $\ga$ inside $E_{\la}$ is complicated (say, non-expanding), but is still relatively expanding with respect to $\pi$. One can think of this setting as {\em actions expanding transversally to the fibers of $\pi$}.

\begin{lemma}
[\cite{coco13, coco15}]\label{lem:exp_coco}
If $\Ga\acts Z$ is expanding relative to $\pi: E\to\La$ then $(Z- E)/\Ga$ is compact (not necessarily Hausdorff, of course).   
\end{lemma}

The idea of the proof is that, if $V$ is a sufficiently small neighborhood of $E$ in $Z$, then $V$ cannot contain the entire $\Ga$-orbit as some points of the orbit will be repulsed away from $E$ (into the complement of $V$) by an expanding element $\ga\in\Ga$.  

{
\begin{example}
Below are two examples of expanding actions with non-Hausdorff quotients:

1. $Z=S^1$, $\Z\cong \Ga < Isom(S^1)$,  $\La = E=\emptyset$. Then the action of $\Ga$ is expanding 
relative to $\pi: E\to\La$, but every $\Ga$-orbit is dense in $S^1$. In particular, $S^1/\Ga$ is infinite with trivial topology. 

2. A more interesting example is given by a cocompact Fuchsian subgroup $\Ga < PSL(2, \R)$ and its product action 
on $Z=S^1\times S^1$. We let $E=\La$ be the diagonal in $S^1\times S^1$ with the identity map $\pi: E\to \La$. 
The action $\Ga\acts Z$ is expanding relative to $\pi$; this can be seen, for instance, by observing that the action  
$\Ga\acts S^1$ is infinitesimally expanding.  On the other hand, $(Z - E)/\Ga$ is non-Hausdorff since the action 
$\Ga\acts Z$ is ergodic and, hence, almost every orbit is dense. 
\end{example}
}

\section{Appendix. Abstract convergence actions and groups}\label{app:congru}

Let $Z$ be a compact metric space which consists of at least three points. We define the space $TZ$ to be the subset of $Z^3$ consisting of triples of pairwise distinct points in $Z$. Every topological action $\Ga\acts Z$ induces a topological action $\Ga \acts TZ$.  

\begin{definition}[Convergence action]
An action $\Ga\acts Z$ is called a {\em convergence action} and the image of $\Ga$ in $\Homeo(Z)$ is said to be a {\em convergence group} if one of the following equivalent conditions hold:

(i) The action $\Ga \acts TZ$ is properly discontinuous. 

(ii) For every  sequence $\ga_n\to\infty$ in $\Ga$ 
there exist points $z_{\pm}\in Z$ and a subsequence of 
$(\ga_n)$  which converges  to the constant map $\equiv z_+$ 
uniformly on compacts in $Z-\{z_-\}$. The points $z_+$ and $z_-$ are called the {\em limit point} (or the {\em attractor}) and the {\em exceptional point} (or the {\em repeller}) of this subsequence.\footnote{Of course, it might happen that $z_-=z_+$.}  

A convergence action $\Ga\acts Z$ is said to be {\em uniform} 
if the action $\Ga\acts TZ$ is cocompact. 
\end{definition}

A proof for the equivalence of the definitions (i) and (ii) can be found in 
\cite{Bowditch_config}.

The main example of convergence actions comes from the following fact: 
Every discrete group $\Ga$ of isometries of a proper Gromov hyperbolic geodesic metric space $X$ acts as a convergence group on the Gromov boundary $\geo X$ of $X$. 
Furthermore, every word hyperbolic group $\Ga$ acts on its Gromov boundary $\geo \Ga$ as a uniform convergence group. See e.g. \cite{Tukia1994}. 

Bowditch proved that, vice versa, 
this dynamical behavior characterizes the natural actions of 
word hyperbolic groups 
on their boundaries: 

\begin{thm}[{\cite[Thm.\ 0.1]{Bowditch_char}}]
\label{thm:charhypbow}
Let $\Ga \acts Z$ be a uniform convergence action 
on a nonempty perfect\footnote{Recall that a topological space is called perfect if it has no isolated points.} 
compact metrizable space. Then $\Ga$ is word hyperbolic and $Z$ is equivariantly homeomorphic 
to $\geo\Ga$. 
\end{thm}

For every convergence action $\Ga\acts Z$ one defines the {\em limit set}, the {\em conical limit set} and the domain 
of discontinuity (which is the same as the domain of proper discontinuity).

\begin{definition}
A sequence $(\ga_n)$ in $\Ga$ is said to {\em converge} to a point $z_+\in Z$, $\ga_k\to z_+$, 
if  every subsequence in $(\ga_n)$ contains a further subsequence, which converges to $z_+$ uniformly on compacts in $Z- \{z_-\}$, for some $z_-\in Z$ (which depends on the subsubsequence). 
\end{definition}

\begin{definition} 
[See Section 8 of \cite{Bowditch_char}] 
A sequence $(\ga_n)$ in $\Ga$ which converges to $z_+$ is said to {\em converge conically} to $z$ if for every point 
$\hat{z}\in Z - \{z\}$, the sequence of pairs $\ga_n^{-1}(z, \hat{z})$ is relatively compact in $Z^2 -Diag(Z^2)$. \end{definition}

\begin{definition}
The {\em limit set} $\La(\Ga)\subset Z$ of a convergence action $\Ga\acts Z$ is the subset consisting of limits
 $z$ of sequences $\ga_k\to z$, $\ga_k\in \Ga$.  The   {\em conical limit set} $\La_c(\Ga)\subset Z$ 
  is the subset consisting of conical limits of sequences $\ga_k\in \Ga$.  
\end{definition}

Both $\La(\Ga)$ and $\La_c(\Ga)$ are $\Ga$-invariant; the limit set $\La(\Ga)$ is closed, while the conical limit set $\La_c(\Ga)$, in general, is not closed. The {\em domain of discontinuity} of the action $\Ga\acts Z$ is the complement 
$Z - \La(\Ga)$. The action of $\Ga$ on $\Om(\Ga)$ is properly discontinuous. An action $\Ga\acts Z$ is called {\em elementary} if $\La(\Ga)$ contains at most two points and {\em nonelementary} otherwise. The limit set of every nonelementary convergence action is perfect. 

In the case when $\Ga$ is a regular antipodal subgroup of the isometry group $G$ of a symmetric space $X$ and 
$Z=\La_{ch}(\Ga)\subset \DF X$, we refer to conical limit points for the convergence action $\Ga\acts Z$ as {\em intrinsically conical} in order to distinguish this notion of conicality from the {\em extrinsic notion} described in Sections \ref{sec:conical_convergence} and \ref{sec:conical_convergence2}.

\begin{thm}[{\cite[Thm.\ 8.1]{Bowditch_char}}, \cite{Tukia}]
\label{thm:bowditch-conical}
A convergence action $\Ga \acts Z$ on a perfect compact metric space $Z$ 
is uniform  if and only if $Z=\La_c(\Ga)$, i.e., every point of $Z$ is a conical limit point of $\Ga$. 
\end{thm}

We now fix a metric $d$ on the metrizable space $Z$. 

\begin{definition}
A convergence action $\Ga\acts (Z,d)$ is {\em expanding} if it is expanding at $\La=\La(\Ga)$ in the sense of Section \ref{sec:expanding_actions}. 
\end{definition}

The following theorem is proven in \cite{morse} using a different method: 

\begin{thm}\label{thm:conical}
Each nonelementary expanding convergence action $\Gamma\acts Z$ restricts to a uniform action on $\La$. In particular, if $Z$ is perfect then for every expanding convergence action $\Gamma\acts Z$, all limit points of $\Ga$ are conical. 
\end{thm}
\proof Our argument will use, and will illustrate, a generalization of the the concept of {\em thickening} discussed earlier in the context of group actions on higher rank symmetric spaces and flag manifolds. For each $\la\in \La=\La(\Ga)$ define its {\em thickening} $Th(\la)\subset \La^3$, 
$$
Th(\la)= \{\la\} \times \{\la\} \times \La \cup  \{\la\} \times  \La \times \{\la\} \cup   \La \times \{\la\} \times \{\la\}.
$$
Clearly, $Th(\ga \la)= \ga Th(\la), \la\in \La, \ga\in \Ga$.  
Note that the subsets $Th(\la)$ are pairwise disjoint (i.e. the thickening is {\em slim}). 
Then 
$$
Th(\La):= \bigcup_{\la\in \La} Th(\la)
$$
is the {\em large diagonal} in $\La^3$. Of course,
$$
T\La= \La^3 - Th(\La). 
$$
We have the $\Ga$-equivariant fibration 
$$
\pi: Th(\La)\to \La, \quad Th(\la)\to \{\la\}. 
$$
We equip $\La^3$ with the following product metric induced from the metric $d$ on $Z$:
$$
d^2((z_1, z_2, z_3), (w_1, w_2, w_3))= d^2(z_1, w_1) + d^2(z_2, w_2) + d^2(z_3, w_3).  
$$
 The fact that the action $\Ga \acts \La$ is expanding translates to the statement that the action $\Ga\acts \La^3$ is expanding relative to $\pi: Th(\La)\to \La$. Indeed, for $z=(z_1, z_2, z_3)$, $z_i\in \La$, and $\la'\in \La$,  
$$
d(z, Th(\la'))= \min\left( \sqrt{d^2 (z_1, \la') + d^2(z_2, \la')}, \sqrt{d^2 (z_1, \la') + d^2(z_3, \la')}, 
\sqrt{d^2 (z_2, \la') + d^2(z_3, \la')} \right). 
$$
The expansion condition for the action of $\ga\in \Ga$ on a neighborhood $U\subset \La$ of $\la\in \La$ (containing the points $z_1, z_2, z_3$ and $\la'$) implies that
$$
d(\ga z_i, \ga \la') \ge c d(z_i, \la'), \quad c>1, i=1, 2, 3. 
$$
From this we conclude that
$$
d(\ga z, Th(\ga \la')) \ge c d(z, Th(\la')),
$$
which means relative expansion. Therefore, according to Lemma \ref{lem:exp_coco}, applied to the action $\Ga\acts (\La^3, Th(\La))$, the action $\Ga\acts T\La= \La^3 - Th(\La)$ is cocompact. 
In other words, $\Ga\acts \La$ is a uniform convergence action  \qed

\section{Appendix. Model spaces: Symmetric spaces and buildings}
\label{sec:modsp}

In these lectures we use two classes of buildings: Spherical and euclidean\footnote{Although the latter are only mentioned in passing in sections 2---5.}. Spherical buildings were introduced by Tits in order to generalize {\em incidence geometry} from classical groups to general semisimple Lie groups; they emerged as an important geometric tool for studying geometry of symmetric spaces and Lie groups. Similarly, euclidean buildings were introduced by Bruhat and Tits as a tool for studying algebraic groups over fields with nonarchimedean  valuations, e.g. $p$-adic numbers $\Q_p$. A way to think about buildings is as hybrids of simplicial complexes and manifolds equipped with geometric structures: From manifolds they acquire an atlas, from simplicial complexes they acquire certain discrete features. We first define {\em model spaces} which include both symmetric spaces and buildings and then add more axioms to specialize to buildings. 

Below are basic axioms of buildings. As in the case of geometric structures, one starts with a {\em model space} and 
a group acting on this space. The model space in our setting is a {\em model apartment} $a_{mod}$, which is either a unit sphere (in the case of spherical buildings) or a euclidean space (for euclidean buildings). One also fixes a {\em model Coxeter group} $W$ acting isometrically on $a_{mod}$; this group is generated by  reflections in hyperplanes in $a_{mod}$. In the case of spherical apartment, $W$ is required to be finite; in the euclidean case $W$ is required to have finite linear part. For euclidean apartments, in general, no discreteness of $W$ is assumed. In order to avoid the notation confusion, we will frequently use the notation $W_{aff}$ for the Coxeter group in the euclidean case. The pair $(a_{mod}, W)$ is called a {\em Coxeter complex}. 

\begin{example}
For instance, if $a_{mod}=A$ is an affine space and $W$ is a finite reflection group of isometries of $A$, then one takes $W_{aff}= T \rtimes W$, where $T$ is the full group of translations of $A$. 
\end{example}

Hyperplanes fixed by reflections in $W$  are called {\em walls} in $a_{mod}$. A (closed) half-apartment in $a_{mod}$ is a half-space bounded by a wall. In the spherical case, $a_{mod}$ is divided into fundamental domains of $W$, called {\em chambers}. The {\em model chamber} $\sigma_{mod}$ is the quotient $a_{mod}/W$, it can be identified with one of the chambers in $a_{mod}$. The quotient projection $\theta: a_{mod}\to  \sigma_{mod}$ is the {\em type map}. In the spherical case one uses the notation $\angle$ for the angular distance on $a_{mod}$. 

Now we can state the main axioms of model spaces:

\medskip
{\bf Axiom 1.} A model space is a metric space which is either a CAT(1) space (spherical case) or a CAT(0) space (euclidean case).

{\bf Axiom 2.} A model space is a metric space $X$ equipped with an atlas where charts are isometric embeddings $a_{mod}\to X$, such that transition maps are restrictions of elements of the model Coxeter group. 

\medskip
Note that images of charts (called {\em apartments} in the case of buildings) are not required to be open in $X$ (unlike the case of geometric structures on manifolds). Using this axioms one can transport various (invariant) notions from the model apartment to the model space; in particular, one defines {\em walls} in $X$ as images of walls in the model apartment. 

\medskip
{\bf Axiom 3.} For any two points $x, y$ of a model space there is a chart whose image contains both $x$ and $y$.

Note that the model apartment clearly satisfies the first three axioms. One frequently adds one more general axiom in order to distinguish model apartments from ``more interesting'' model spaces:

\medskip
{\bf Axiom 4.} A model space $X$ is {\em thick} if any wall in $X$ equals the intersection of three half-apartments. 

\begin{example}
1. A metric tree is a model space modeled on $(A,W)$, where $A$ is the line and $W$ is the full group of euclidean isometries of $A$.  

2. Each symmetric space $X$ is a model space. Let $G$ denote the connected component of the identity in the isometry group of $X$.  The model apartment is a maximal flat $F$ in $X$ and the model Coxeter group acting on $F$ is the image $W_{aff}$ in $Isom(F)$ of the stabilizer $G_F$ of $F$ in $G$.     
\end{example}

In order to differentiate between symmetric spaces and buildings, one introduces one more {\em angle discreteness} axiom. This axiom is void in the spherical case and we, therefore, restrict now to the case of euclidean model spaces with the Coxeter group $W_{aff}=T\rtimes W$, where $W$ is a finite Coxeter group. Let $\sigma_{mod}$ denote the model chamber for action of $W$ on the unit sphere in the affine space $A=a_{mod}$. We let $\Delta\subset A$ denote a model euclidean Weyl chamber, the cone over $\sigma_{mod}$. We then have the $\Delta$-distance function
$$
d_\Delta(x,y)\in \Delta
$$
defined for all points $x,y\in X$: Pick a chart $\phi: A\to X$; $\phi(x')=x, \phi(y')=y$, then consider the vector $v$ in $A$ represented by the directed segment $x'y'$ and project $v$ to $\Delta=A/W$. The result is $d_\Delta(x,y)$. 
For each nondegenerate geodesic segment $xy$ in $X$, we define its 
$\sigma_{mod}$-{\em direction} $\theta(xy)$ as the unit vector in the direction of $d_\Delta(x,y)$.

For each $x\in X$ one has the {\em space of directions} $\Sigma_x X$ of $X$ at $x$, which is the space of germs of nondegenerate geodesic segments emanating from $x$. On this space we have the (metric) notion of {\em angle} denoted $\angle$. We also have the {\em type map} $\theta_x: \Sigma_x X\to \sigma_{mod}$, sending each $xy$ to 
its direction $\theta(xy)$.

\medskip
{\bf Axiom 5: Angle discreteness}. Let $v_1, v_2$ be elements of $\Si_xX$. Then we require that $\angle(v_1,v_2)$ 
belongs to the finite set of angles
$$
\angle(\theta(v_1), w \theta(v_2)), \quad w\in W. 
$$

A Riemannian manifold (of dimension $\ge 2$) cannot satisfy this axiom. 

\begin{definition}
A (thick) spherical building is a model space, modeled on a spherical Coxeter complex and satisfying Axioms 1---4. A (thick) euclidean building is a model space  modeled on a euclidean Coxeter complex and satisfying Axioms 1---5.
\end{definition}

For instance, metric trees and their products are examples of euclidean buildings. A building is said to be {\em discrete} if the model Coxeter group is discrete. Below is an example of a discrete euclidean building $X$ on which $PSL(3,\Q_p)$ acts isometrically. We will only describe the underlying simplicial complex and not the rest of the structure. Vertices of $X$ are equivalence classes of $\Z_p$-lattices $\La$ in $\Q_p^3$, where two lattices are equivalent iff they differ by a $\Q_p$-scaling. Edges of $X$ represented by pairs of lattices: 
$$
\La\subset \La', \quad |\La':\La|=p. 
$$
A 2-simplex in $X$ is a chain of proper inclusions of lattices 
$$
\La_0\supset \La_1 \supset \La_2 \supset p\La_0
$$
where each inclusion is necessarily of the index $p$ (note that $\La_0, p\La_0$ represent the same vertex). 

It turns out that for each point $x$ in a building (spherical or euclidean), the space of directions $\Si_xX$ has a natural structure of a spherical building. In the case of a vertex in a discrete building, the building $\Si_xX$ (treated as a simplicial complex)  is identified with the link of $x$.

Buildings enter naturally into the theory of symmetric spaces via two {\em asymptotic} constructions:

1. The visual boundary of a symmetric space $X$ (also, of a euclidean building) is a spherical building, see e.g. \cite{Eberlein}. 

2.  Every {\em asymptotic cone} of a symmetric space $X$  
is again a euclidean building modeled on the same Coxeter complex, see \cite{qirigid}. 

\medskip
{\bf The $SL(3,\R)$ example.}  Below we supplement our earlier discussion of visual boundaries with 
the detailed example of the symmetric space $X$ of the group $G=SL(3,\R)$. Our treatment  follows \cite[Appendix 5]{BGS}. This symmetric space $X$ is identified with the space $P(3,\R)$ of conformal structures on $\R^3$, more precisely, positive definite bilinear forms $b$ on $\R^3$ up to scalar. After fixing the standard euclidean bilinear form $q_0$ on $\R^3$, such structures can be identified with symmetric positive-definite $3\times 3$ matrices $A$ of the unit determinant:
$$
b(u, v)= u^T A v. 
$$
The group $G$ acts on bilinear forms $b$ 
by change of variables:
$$
g^*b=b', b'(u, v)= b(g(u), g(v)), 
$$
in terms of matrices, the action is given by
$$
g^* A= g^T A g, \quad g\in SL(3,\R).   
$$
In matrix terms, the Riemannian metric on $X$ at the identity matrix is given by
$$
( a, b)= tr(ab), 
$$
where $a, b$ are symmetric traceless $3\times 3$ matrices. 

A maximal flat $F$ in $X$ is given by diagonal matrices with positive diagonal entries $x_1, x_2, x_3$ 
and unit determinant (the corresponding quadratic forms have principal axes equal to the coordinate lines in $\R^3$). 
The isometry of $F$ to the euclidean plane is given by
$$
Diag(x_1, x_2, x_3)\mapsto (\log(x_1), \log(x_2), \log(x_3))\in \R_0^3=$$
$$
 \{(y_1,y_2,y_3)\in\R^3: y_1+ y_2+y_3=0\}. 
$$
The action of $W=S_3$ on $F$ is by permuting the diagonal entries. The walls are given by the equations 
$$
x_i=x_j, 1\le i\ne j\le 3. 
$$
The positive chamber in $F$ 
is defined by the inequalities
$$
x_1\ge x_2 \ge x_3>0 $$
and in $\R^3_0$ by the inequalities
$$
y_1\ge y_2 \ge y_3. 
$$
(We will think of the positive chamber $\Delta$ as sitting in $\R^3_0$ since the metric of the symmetric space equals the euclidean metric in this setting.) 
If $q\in X$ is a general quadratic form, then  the segment $q_0 q$ is regular iff all three principal axes of the ellipsoid of $q$ are distinct. The $\Delta$-distance $d_\Delta(q_0, q)$ is the vector
$$
(\log(\la_1), \log(\la_2), \log(\la_3)),
$$
where $\la_i$'s are the eigenvalues of the matrix $A$ of $q$, which are arranged in the descending order. In terms of a matrix $g\in SL(3,\R)$,
$$
d_\Delta(q_0, g(q_0))= (\log(\mu_1), \log(\mu_2), \log(\mu_3))
$$
where $\mu_i$'s are the singular values of the matrix $g$, again arranged in the descending order. A sequence $(g_i)$ is  regular iff
$$
\lim_{i\to\infty} \frac{\mu_1(g_i)}{\mu_2(g_i)}= \lim_{i\to\infty} \frac{\mu_2(g_i)}{\mu_3(g_i)}=\infty. 
$$

\medskip
We now describe the visual boundary of $X$ which we will identify with the space of asymptotic classes of geodesic rays $\rho(t), t\ge 0$, emanating from $q_0$. The initial velocity of such a ray is the matrix with the eigenvalues $a\ge b\ge c, a+b+c=0, a^2+b^2+c^2=1$; the asymptotic behavior of the ray is determined by the parameter    
$$
r= \frac{b-c}{a-b}, \quad r\in [0,\infty].   
$$
Let $v_a, v_b, v_c\in \R^3$ denote the unit eigenvectors corresponding to the eigenvalues $a, b, c$. Define the strip
$$
s_\rho= \R v_a + [-r,r]v_b. 
$$

There are three possibilities: 

1. $r=0$, i.e., $b=c$, the strip degenerates to the line $\R v_a$. We associate with the ray $\rho$ the line $\<v_a\>$. Note that the vectors $v_b, v_c$ are not uniquely defined. This corresponds to the fact that there is no unique maximal flat through the ray $\rho$. 

2. $r=\infty$, i.e., $a=b$, the strip is the plane $\<v_a, v_b\>$ spanned by $v_a, v_b$. We associate with the ray $\rho$ the plane $\<v_a, v_b\>$. Again, the vectors $v_a, v_b$ are not uniquely defined (only the vector $v_c$ and the plane $\<v_a, v_b\>= v_c^\perp$ are well-defined). This again corresponds to nonuniqueness  of a maximal flat through the ray $\rho$.

3. $r\in (0,\infty)$, equivalently, the direction of $\rho$ is regular. The corresponding flag is 
$$
(L\subset P)= (\<v_a\>\subset \<v_a, v_b\>).$$
 Note that in this case $a>b>c$ and the one-dimensional eigenspaces of $A$ are unique.  

If we fix $v_a, v_b$, the set of resulting geodesic rays corresponds to the set of all strips in $\<v_a, v_b\>$ interpolating between the line $\<v_a\>$ and the plane $\<v_a, v_b\>$. This is the Weyl chamber corresponding to the flag $L\subset P$. 

If we equip the set of strips with Gromov-Hausdorff topology, then the resulting set is homeomorphic to the set of rays emanating from $q_0$, i.e., to the visual boundary of $P(3,\R)$. The subspaces corresponding to the sets of singular directions are homeomorphic to the spaces of lines and, resp. planes, in $\R^3$.  Thus we see that the Tits boundary of $P(3,\R)$ is naturally homeomorphic to the incidence graph of the projective plane. 

{
\section{Appendix. Manifolds with corners} \label{sec:corners}

For simplicity we consider here only the concept of manifolds with corners and 
good orbifolds with corners since only they appear in the context of this paper. 
We refer to \cite{Joyce} for the definition of orbifolds with corners in general.

The concept of manifolds with corners generalizes the notion of manifolds with boundary. 
Recall that the latter are defined via atlases with values in the euclidean spaces and half-spaces. 
Let $I$ denote the closed interval $[0,1]$ and $I^n$ the $n$-dimensional cube. The cube $I^n$ is a {\em stratified space} 
where the (open) $k$-dimensional stratum $S_k(I^n)$ is the union of open $k$-dimensional faces of $I^n$, i.e. 
$k$-dimensional subcubes which are products of several copies of the open intervals $(0,1)$ and 
singletons from the set $\{0, 1\}$. Dimensions of the strata range from $0$ to $n$. We let $H(I^n)$ denote 
the pseudogroup of homeomorphisms between open subsets of $I^n$ which preserve the stratification of $I^n$, i.e. 
map points of $S_k(I^n)$ to points of $S_k(I^n)$ for every $k=0,...,n$.

\begin{definition}
An $n$-dimensional topological manifold with corners is a 2nd countable Hausdorff topological space $X$ equipped with 
a certain {\em atlas}, which is a maximal system of homeomorphisms (``charts'') $\phi_\alpha: U_\alpha\to V_\alpha$, from open subsets $U_\al\subset I^n$ to open subsets $V_\al\subset X$. It is required that the (partially defined) transition maps 
$$
g_{\al,\be}= \phi_\beta^{-1}\circ \phi_\al
$$
preserve the stratification of $I^n$: $g_{\al,\be}\in H(I^n)$ for all  $\phi_\beta, \phi_\al$. 
Thus, every manifold with corners is stratified as
$$
S_0(X)\sqcup S_1(X) \sqcup \ldots \sqcup S_n(X),
$$
where the {\em strata} $S_k(X)$ consist of the points $x\in  X$ which are mapped to 
$S_k(I^n)$ under the maps $\phi_\al^{-1}$. 
\end{definition}

\begin{example}
If $X$ is a manifold with corners and $Y\subset X$ is an open subset, 
then one obtains a pull-back structure of the manifold with corners from $X$ to $Y$, where the  charts for $Y$ 
are suitable restrictions of the charts for $X$. Similarly, one defines the pull-back of the manifold with corners structure 
via a local homeomorphism $f: Y\to X$, where $Y$ is Hausdorff and 2nd countable and $X$ is a manifold with corners. 
\end{example}

We note that, in particular, every  manifold with corners is automatically a  manifold with boundary,  
where $S_n(X)=\inte(X)$ and the union of the rest of the strata is the boundary of $X$. Conversely, every $n$-manifold with boundary $X$ has a natural structure of a manifold with corners where all strata of dimension $< n-1$ are empty. 

Unlike the topological boundary, the strata $S_k(X)$ are not uniquely determined by the topology of $X$. 

\begin{example}
The cube $X=I^n$ has a natural manifold with corners structure which is the maximal atlas containing the identity map $I^n\to I^n$ and $S_k(X)=S_k(I^n)$, $k=0,...,n$. On the other hand, $I^n$ is homeomorphic to the closed ball $B^n$ 
which we treat as a manifold with boundary. Hence, $S_k(B^n)=\emptyset$ for $k<n-1$. 
\end{example}

As in the case of manifolds, one defines manifolds with corners in other categories, e.g. smooth (with smooth transition maps) or real-analytic (with real-analytic transition maps). Here we recall that a smooth (respectively, real analytic) 
map of an open subset of $I^n\subset \R^n$ to $\R^k$ is a map which admits a smooth (respectively, real analytic) extension to an open subset of $\R^n$. A homeomorphism $f: X\to Y$ between $n$-dimensional manifolds with corners  is 
an {\em isomorphism} if for every pair of charts $\phi_\al: U_\al\subset I^n \to X$,  
$\phi_\be: U_\be\subset I^n\to Y$ the composition $\phi_\be^{-1} \circ f \circ \phi_\al$ belongs to the pseudogroup 
$H(I^n)$. Similarly, one defines an isomorphism between smooth (respectively, real analytic) manifolds with corners 
as an isomorphism of the underlying topological manifolds with corners which gives rise to local diffeomorphisms 
(respectively, real-analytic diffeomorphisms) $\phi_\be^{-1} \circ f \circ \phi_\al$. We let $Aut(X)$ denote the group 
of automorphisms of the topological (respectively, smooth, or real-analytic) manifold with corners $X$. 

\begin{definition}
Let $X$ be a real analytic manifold with corners and $\Gamma < Aut(X)$ be a subgroup which acts properly discontinuously on 
$X$. Then the quotient $X/\Gamma$ is called a 
{\em good real-analytic orbifold with corners}. Here $X/\Gamma$  is a topological space equipped with the collection of 
{\em orbi-charts} which are quotient maps of open subsets $U_\al\subset I^n$ by finite groups of real-analytic automorphisms of $U_\al$ obtained by restrictions of subgroups of $\Gamma$.  
\end{definition}

We note that (analogously to an ordinary orbifold) the orbi-charts of a good orbifold with corners satisfy certain compatibility conditions which are used to define an orbifold with corners in full generality, see \cite{Joyce}.

}

\medskip 

\noindent M.K.: Department of Mathematics,  
University of California, Davis,  
CA 95616, USA\\
email: kapovich@math.ucdavis.edu\\

and

\medskip 
\noindent Korea Institute for Advanced Study,\\ 
207-43 Cheongnyangri-dong, Dongdaemun-gu,\\ 
Seoul, South Korea\\

\noindent B.L.: Mathematisches Institut, 
Universit\"at M\"unchen, 
Theresienstr. 39,  
D-80333 M\"unchen, Germany,  
email: b.l@lmu.de


\begin{thebibliography}{BLP05}

\bibitem[Ah]{Ahlfors}
L. Ahlfors, {\em Fundamental polyhedrons and limit point sets of Kleinian groups} 
Proc. Nat. Acad. Sci. U.S.A. Vol. 55 (1966) pp. 251--254.

\bibitem[Al]{Albuquerque}
P.\ Albuquerque, 
{\em Patterson-Sullivan theory in higher rank symmetric spaces}, 
Geom. Funct. Anal. Vol. {9} (1999), no. 1, pp. 1-28. 

\bibitem[Ba]{Ballmann}
W. Ballmann,
``Lectures on spaces of nonpositive curvature". With an appendix by Misha Brin. 
DMV Seminar, vol.\ 25, Birkh\"auser Verlag, Basel, 1995.


\bibitem[BGS]{BGS}
W. Ballmann, M. Gromov and V. Schroeder, ``Manifolds of nonpositive curvature'', Birkh\"auser Verlag, 1985. 


\bibitem[BM]{Beardon-Maskit}
A. Beardon, B. Maskit, 
{\em Limit points of Kleinian groups and finite sided fundamental polyhedra}, 
Acta Math. Vol. 132 (1974), pp. 1--12. 

\bibitem[Ben]{Benoist}
Y.\ Benoist, 
{\em Propri\'et\'es asymptotiques des groupes lin\'eaires}, Geom. Funct. Anal. Vol. {7} (1997), no. 1, pp. 1--47. 


\bibitem[Bi]{Bishop}
C. Bishop, {\em On a theorem of Beardon and Maskit}, 
Ann. Acad. Sci. Fenn. Math. Vol. 21 (1996), no. 2, pp. 383--388.

\bibitem[BB]{BB}
A. Bj\"{o}rner, F. Brenti, ``Combinatorics of Coxeter groups.'' 
Graduate Texts in Mathematics, 231. Springer, New York, 2005.


\bibitem[BCSS]{BCSS}
L. Blum, F. Cucker, M. Shub, S. Smale, {\em Complexity and Real Computation}, Springer Verlag, 1997. 

\bibitem[BJ]{Borel-Ji}
A. Borel, L. Ji, ``Compactifications of Symmetric and Locally Symmetric Spaces'', 
Springer Verlag, 2006. 


\bibitem[BT]{Borel-Tits} A. Borel, J. Tits, {\em Compl\'ements \'a l'article: ``Groupes r\'eductifs''}. Math. Publ. IHES, No. 41 (1972), pp. 253--276. 


\bibitem[Bou]{Bourbaki}
N.\ Bourbaki, 
{\em Lie Groups and Lie Algebras}, ch.\ 4-6, Springer Verlag, 2002.

\bibitem[Bo1]{Bowditch93}
B.\ Bowditch, 
{\em Geometrical finiteness for hyperbolic groups}, J. Funct. Anal. Vol. 113 (1993), no. 2, pp. 245--317. 


\bibitem[Bo2]{Bowditch_gf}
B.\ Bowditch, 
{\em Geometrical finiteness with variable negative curvature},
Duke Math. J. Vol. {77} (1995) pp. 229-274.


\bibitem[Bo3]{Bowditch-stab}
B.\ Bowditch, 
{\em Spaces of geometrically finite representations}, 
Ann. Acad. Sci. Fenn. Math. Vol. 23 (1998), no. 2, pp. 389--414.


\bibitem[Bo4]{Bowditch_char}
B.\ Bowditch,
{\em A topological characterisation of hyperbolic groups},
J. Amer. Math. Soc. Vol. {11} (1998), no. 3, pp. 643--667.


\bibitem[Bo5]{Bowditch_config}
B.\ Bowditch, 
{\em Convergence groups and configuration spaces}, 
in ``Geometric group theory down under'' (Canberra, 1996), 
pp.\ 23-54, de Gruyter, Berlin, 1999. 


\bibitem[BrHa]{BH} M. Bridson, A. Haefliger, 
``Metric spaces of non-positive curvature'',  
Grundlehren der Mathematischen Wissenschaften, 319. Springer-Verlag, Berlin, 1999.



\bibitem[Ch]{Chevalley}
C. Chevalley,  {\em Sur les d\'ecompositions cellulaires des espaces $G/B$.} With a foreword by Armand Borel. Proc. Sympos. Pure Math., 56, Part 1, Algebraic groups and their generalizations: classical methods (University Park, PA, 1991), pp. 1--23, Amer. Math. Soc., Providence, RI, 1994. 


\bibitem[CDP]{CDP}
M. Coornaert, T. Delzant, A. Papadopoulos, ``Les groupes hyperboliques de Gromov,'' Lecture notes in Mathematics No. 1441, Springer Verlag (1990). 

\bibitem[CP]{CP}
M.\ Coornaert, A.\ Papadopoulos,
{\em Horofunctions and symbolic dynamics on Gromov hyperbolic groups}, 
Glasg.\ Math.\ J., 43 (2001), no.\ 3, pp. 425-456. 

\bibitem[D]{Dolgachev}
I. Dolgachev, {\em Lectures on invariant theory}, 
LMS Lecture Notes Series, 296, Cambridge University Press, 2003. 

\bibitem[Eb]{Eberlein}
P. Eberlein, ``Geometry of nonpositively curved manifolds'', University of Chicago Press, 1997. 

\bibitem[FHS]{FHS}
M. Farber, J.-C. Hausmann, and D. Sch\"utz, 
{\em On the conjecture of Kevin Walker}, J. Topol. Anal. Vol. 1 (2009), no. 1, pp. 65--86. 



\bibitem[Gi1]{Gilman1} 
J. Gilman, {\em Two-generator discrete subgroups of $PSL(2,\R)$},  
Mem. Amer. Math. Soc., Vol. 117 (1995) no. 561. 

\bibitem[Gi2]{Gilman2} 
J. Gilman, {\em Algorithms, Complexity and Discreteness Criteria in $PSL(2,\C)$}, 
Journal d'Analyse Math\'ematique, Vol. 73 (1997) pp. 91--114.


\bibitem[GM]{Gilman-Maskit}
J. Gilman, B. Maskit, 
{\em An algorithm for 2-generator Fuchsian groups}, 
Michigan Math. J., Vol. 38 (1991), no. 1, pp. 13--32. 

\bibitem[Gr1]{Gromov_hypmfs}
M. Gromov,
{\em Hyperbolic manifolds, groups and actions}, 
Riemann surfaces and related topics: 
Proceedings of the 1978 Stony Brook Conference (State Univ. New York, Stony Brook, N.Y., 1978), pp. 183--213,
Ann. of Math. Stud., vol. 97, Princeton Univ. Press, Princeton, N.J., 1981. 

\bibitem[Gr2]{Gromov_hypgps}
M. Gromov,
{\em Hyperbolic groups}, 
in: "Essays in group theory", Math. Sci. Res. Inst. Publ. 8, Springer, New York (1987), pp. 75--263.

\bibitem[GW]{GW}
O.\ Guichard, A.\ Wienhard, 
{\em Anosov representations: Domains of discontinuity and applications}, 
Invent. Math. Vol. 190 (2012) no. 2, pp. 357--438. 

\bibitem[GGKW1]{GGKW1}
F. Gu\'eritaud, O.\ Guichard, F.\ Kassel, A.\ Wienhard, 
{\em Anosov representations and proper actions}, 
Preprint arXiv:1502.03811, February 2015. 

\bibitem[GGKW2]{GGKW2}
F. Gu\'eritaud, O.\ Guichard, F.\ Kassel, A.\ Wienhard, 
{\em Tameness of Riemannian locally symmetric spaces arising from Anosov representations},
Preprint arXiv:1508.04759, August 2015.

\bibitem[Ha]{Hamilton}
E. Hamilton, {\em Geometric finiteness for hyperbolic orbifolds}, Topology, Vol. 37 (1998), pp. 635--657.

\bibitem[He]{Helgason}
S. Helgason, ``Differential geometry, Lie groups and symmetric spaces'',  AMS series Graduate Studies in Mathematics, 2001. 

\bibitem[HP]{HP}
M. Heusener, J. Porti, {\em Representations of knot groups into $SL_n(\C)$ and twisted Alexander polynomials}, Pacific J. Math. 277 (2015), no. 2, pp. 313--354.


\bibitem[Hi]{Hitchin}
N.\ Hitchin,
{\em Lie groups and Teichm\"uller space}, 
Topology, Vol. 31 (1992), no. 3, pp. 449--473. 

\bibitem[Iz]{Izeki}
H. Izeki, {\em Quasiconformal stability of Kleinian groups and an embedding of a space of flat conformal structures}, Conform. Geom. Dyn., 4 (2000), pp. 108--119.

{ \bibitem[JS]{JS}
L. Ji, A-S. Schilling, {\em Polyhedral horofunction compactification as polyhedral ball}, e-reprint, 
arXiv:1607.00564. }

\bibitem[Jo]{Joyce} 
D. Joyce, {\em D-manifolds, d-orbifolds and derived differential geometry: a detailed summary}, 
e-print, arXiv:1208.4948.

\bibitem[Kai]{Kaimanovich}
V. Kaimanovich, {\em Lyapunov exponents, symmetric spaces and a multiplicative ergodic theorem for semisimple Lie groups}, J. Soviet Math. Vol. 47 (1989), no. 2, pp. 2387--2398.  

\bibitem[Ka1]{Kapovich00}
M. Kapovich, ``Hyperbolic Manifolds and Discrete Groups'', Birkh\"auser's series ``Progress in Mathematics'', 2000.

\bibitem[Ka2]{Kapovich2007}
M.\ Kapovich, {\em Kleinian groups in higher dimensions}, In ``Geometry and Dynamics of Groups and Spaces. In memory of Alexander Reznikov'', M. Kapranov et al (eds). Birkh\"auser, Progress in Mathematics, Vol. 265, 2007, pp. 485--562.  

\bibitem[Ka3]{Kapovich2015}
M.\ Kapovich, {\em Discreteness is undecidable},  International Journal of Algebra and Computation,  
26 (2016) pp. 467--472. 

\bibitem[KLP1a]{coco13}
M.\ Kapovich, B.\ Leeb, J.\ Porti, 
{\em Dynamics at infinity 
of regular discrete subgroups of isometries 
of higher rank symmetric spaces}, arXiv e-print, arXiv:1306.3837v1, June 2013. 

\bibitem[KLP1b]{coco15}
M.\ Kapovich, B.\ Leeb, J.\ Porti,
{\em Dynamics on flag manifolds:
domains of proper discontinuity and cocompactness},
arXiv e-print, 2nd version of \cite{coco13}, arXiv:1306.3837v2, 
October 2015.


\bibitem[KLP2]{morse}
M.\ Kapovich, B.\ Leeb, J.\ Porti, 
{\em Morse actions of discrete groups on symmetric spaces}, arXiv e-print,  arXiv:1403.7671,  March 2014. 

\bibitem[KLP3]{mlem}
M.\ Kapovich, B.\ Leeb, J.\ Porti, 
{\em A Morse Lemma for quasigeodesics in symmetric spaces and euclidean buildings},  arXiv e-print, 
arXiv:1411.4176, November 2014. 


\bibitem[KLP4]{anolec}
M.\ Kapovich, B.\ Leeb, J.\ Porti, 
{\em Anosov subgroups: Dynamical and geometric characterizations}, Preprint, January 2016. 

\bibitem[KLP5]{anosov}
M.\ Kapovich, B.\ Leeb, J.\ Porti,
{\em Some recent results on Anosov representations},
Transformation Groups 21 (2016), no. 4, 1105--1121. 


\bibitem[KL1]{bordif}
M.\ Kapovich, B.\ Leeb, {\em Finsler bordifications of symmetric and certain locally symmetric spaces},  arXiv e-print, arXiv:1505.03593, 2015. 

\bibitem[KL2]{KL-undistorted}
M.\ Kapovich, B.\ Leeb, {\em Examples of group actions on products of simplicial trees}, in preparation. 


\bibitem[KLM]{KLM}
M.\ Kapovich, B.\ Leeb, J. J.\ Millson, 
{\em Convex functions on symmetric spaces, side lengths of polygons 
and the stability inequalities for weighted configurations at infinity}, 
Journal of Differential Geometry, Vol. 81 (2009), pp. 297--354. 

\bibitem[KM]{KM}
M. Kapovich, J. J. Millson, {\em On the moduli space of polygons in the euclidean plane}, 
J. Differential Geom. Vol. 42 (1995), no. 1, pp. 133--164. 

\bibitem[KR]{KR}
M. Kapovich, K. Rafi, {\em Teichm\"uller space is semihyperbolic}, in preparation. 

\bibitem[Kar]{Karlsson}
A.\ Karlsson, 
{\em On the dynamics of isometries}, 
Geom. Topol. Vol. 9 (2005), pp. 2359--2394. 

\bibitem[KL98]{qirigid}
B.\ Kleiner, B.\ Leeb, 
{\em Rigidity of quasi-isometries for symmetric spaces 
and euclidean buildings}, Inst.\ Hautes \'Etudes Sci. Publ.\ Math.\ No.\ 86 (1997) pp. 115--197.

\bibitem[KL06]{convcoco}
B.\ Kleiner, B.\ Leeb, 
{\em Rigidity of invariant convex sets in symmetric spaces}, 
Invent. Math. Vol. 163, No. 3, (2006) pp. 657--676. 

\bibitem[La]{Labourie}
F.\ Labourie, 
{\em Anosov flows, surface groups and curves in projective space}, 
Invent. Math. Vol. 165, No. 1, (2006) p. 51--114. 

\bibitem[Man]{Manning}
J. Manning, {\em Algorithmic detection and description of hyperbolic structures on closed 3-manifolds with solvable word problem}, Geom. Topol. 6 (2002), pp. 1--25. 

\bibitem[Mar]{Marden}
A. Marden, {\em Geometry of finitely generated Kleinian groups}, Ann. of Math., 99 (1974), pp. 383--496.

\bibitem[Mas]{Maskit}
B. Maskit, ``Kleinian groups'', vol. 287 of Grundlehren der Math. Wissenschaften, Springer, 1987.

\bibitem[Mi]{Mineyev}
I.\ Mineyev, 
{\em Flows and joins of metric spaces}, 
Geom. Topol. 9 (2005), pp. 403--482. 

\bibitem[Mi1]{Mitchell} S. Mitchell, 
{\em Quillen's theorem on buildings and the loops on a symmetric space}, 
Enseign. Math. (2) 34 (1988), no. 1-2, pp. 123--166.

\bibitem[Mi2]{Mitchell2008}
S. Mitchell, {\em Parabolic orbits in flag varieties}, Preprint, 2008, 
\newline
http://www.math.washington.edu/~mitchell/Notes/flag.pdf. 

\bibitem[Mo]{Mostow}
G.D.\ Mostow, 
{\em Strong rigidity of locally symmetric spaces},
Annals of Mathematics Studies, No. 78, 
Princeton University Press, 1973.

\bibitem[Mu]{Mumford}
D.\ Mumford, J. Fogarty, F. Kirwan, {\em Geometric Invariant Theory}, Springer Verlag, 1994. 

\bibitem[N]{Newstead}
P. E. Newstead, {\em Geometric Invariant Theory}, CIMAT Lectures, 2006. 

\bibitem[PS]{PapaSwen}
P.\ Papasoglu, E.\ Swenson, 
{\em Boundaries and JSJ decompositions of CAT(0)-groups}, 
Geom. Funct. Anal. 19 (2009), no. 2, pp. 559--590. 

\bibitem[Park]{Parker}
J. Parker, ``Hyperbolic spaces'', Jyv\"askyl\"a
Lectures in Mathematics, Vol. 2, 2008.

\bibitem[Parr]{Parreau}
A.\ Parreau, 
{\em La distance vectorielle dans les immeubles affines et les espaces sym\'et\-ri\-ques},
in preparation. 

\bibitem[Q]{Quint}
J.-F.\ Quint,
{\em Groupes convexes cocompacts en rang sup\'erieur},
Geom.\ Dedicata 113 (2005), pp. 1--19. 

\bibitem[Rag]{Raghunathan}
M. Raghunathan, ``Discrete subgroups of Lie groups'', Ergebnisse der Mathematik und ihrer Grenzgebiete, Band 68. Springer-Verlag, New York-Heidelberg, Springer Verlag, 1972. 

\bibitem[Rat]{Ratcliffe}
J. Ratcliffe, ``Foundations of hyperbolic manifolds'', Springer Verlag, 1994.


\bibitem[Ri]{Riley}
R. Riley, {\em Applications of a Computer Implementation
of Poincar\'e's Theorem on Fundamental Polyhedra}, Math of Computation, Vol. 40 (1983) pp. 607--632. 


\bibitem[S1]{Schutz}
D. Sch\"utz, {\em The isomorphism problem for planar polygon spaces}, 
J. Topol. Vol. 3 (2010), no. 3, pp. 713--742. 

\bibitem[S2]{Schutz1}
D. Sch\"utz, {\em Intersection homology of linkage spaces in odd-dimensional euclidean space},  
Algebr. Geom. Topol. Vol. 16 (2016), no. 1, pp. 483--508. 

\bibitem[S3]{Schutz2}
D. Sch\"utz, {\em Intersection homology of linkage spaces},  
J. Topol. Anal. Vol. 8 (2016), no. 1, pp. 25--58. 

\bibitem[Su]{Sullivan}
D.\ Sullivan, 
{\em Quasiconformal homeomorphisms and dynamics. II. Structural stability
implies hyperbolicity for Kleinian groups}, 
Acta Math. Vol. {155} (1985), no. 3-4, pp. 243--260.

\bibitem[Th]{Thurston}
W. Thurston, ``The geometry and topology of three-manifolds'', Princeton lecture notes (1978--1981).


\bibitem[Ti]{Tits}
J. Tits, 
{\em Free subgroups in linear groups}, 
J. Algebra, Vol. {\bf 20} (1972) pp. 250--270.


\bibitem[Tu94]{Tukia1994}
P. Tukia, 
{\em Convergence groups and Gromov-hyperbolic metric spaces}, 
New Zealand Journal of Math.\  Vol. {23} (1994), pp. 157--187. 


\bibitem[Tu98]{Tukia}
P. Tukia, 
{\em Conical limit points and uniform convergence groups}, 
J.\ Reine Angew.\ Math.\  Vol. {501} (1998), pp. 71--98. 



\bibitem[Y]{Yue}
C. Yue, {\em Dimension and rigidity of quasifuchsian representations}, Ann. of Math. (2), Vol. 143 (1996), pp. 331--355.


\end{thebibliography}
\end{document}